\documentclass[review]{elsarticle}

\usepackage{hyperref, verbatim,soul}
\usepackage{amssymb}
\usepackage{amsthm}
\usepackage{float}
\usepackage{graphicx}
\usepackage{calc}
\usepackage{xcolor}

\usepackage[linesnumbered,noend,ruled]{algorithm2e}
\SetAlCapSty{}
\usepackage{algpseudocode}
\SetKw{Break}{break}
\SetKw{Continue}{continue}

\theoremstyle{definition}
\newtheorem{thm}{Theorem}[section]
\newtheorem{dfn}[thm]{Definition}
\newtheorem{pro}[thm]{Problem}

\newtheorem{cor}[thm]{Corollary}
\newtheorem{lem}[thm]{Lemma}
\newtheorem{exa}[thm]{Example}

\usepackage{chngcntr}
\counterwithin{table}{section}

\newcommand{\R}{\mathbb R}
\newcommand{\Z}{\mathbb Z}

\newcommand{\al}{\alpha}
\newcommand{\be}{\beta}
\newcommand{\ga}{\gamma}
\newcommand{\Ga}{\Gamma}
\newcommand{\de}{\delta}
\newcommand{\la}{\lambda}
\newcommand{\La}{\Lambda}
\newcommand{\ti}{\tilde}
\newcommand{\si}{\sigma}
\newcommand{\ep}{\varepsilon}
\newcommand{\ph}{\varphi}

\newcommand{\om}{\omega}

\newcommand{\IT}{\mathrm{IT}}

\newcommand{\Or}{\mathrm{O}}
\newcommand{\SO}{\mathrm{SO}}

\newcommand{\GC}{\mathrm{GC}}

\newcommand{\AMD}{\mathrm{AMD}}
\newcommand{\PDD}{\mathrm{PDD}}
\newcommand{\EMD}{\mathrm{EMD}}
\newcommand{\BT}{\mathrm{BT}}
\newcommand{\sym}{\mathrm{Sym}}

\newcommand{\vol}{\mathrm{vol}}
\newcommand{\bs}{\hfill $\blacksquare$}
\newcommand{\lra}{\leftrightarrow}
\newcommand{\bd}{\partial}
\newcommand{\vl}{\,:\,}
\newcommand{\angstrom}{\textup{\AA}}
\newcommand{\ar}[1]{\langle #1 \rangle}

\journal{Pattern Recognition}

\begin{document}

\begin{frontmatter}

\title{Recognition of near-duplicate periodic patterns by continuous metrics with approximation guarantees}

\author[address]{Olga~D~Anosova}
\author[address]{Daniel~E~Widdowson}
\author[address]{Vitaliy~A~Kurlin\corref{corresponding}}
\address[address]{Department of Computer Science, University of Liverpool, Liverpool L69 3BX, UK}
\cortext[corresponding]{
The corresponding author was supported by the Royal Society APEX fellowship (APX/R1/231152), Royal Academy of Engineering Fellowship 
(IF2122/186) at the Cambridge Crystallographic Data Center, and the EPSRC New Horizons grant 
(EP/X018474/1).}
\ead{vkurlin@liverpool.ac.uk, http://kurlin.org}

\begin{abstract}
This paper rigorously solves the challenging problem of recognizing periodic patterns under rigid motion in Euclidean geometry.
The 3-dimensional case is practically important for justifying the novelty of solid crystalline materials (periodic crystals) and for patenting medical drugs in a solid tablet form.  
\smallskip

Past descriptors based on finite subsets fail when a unit cell of a periodic pattern discontinuously changes under almost any perturbation of atoms, which is inevitable due to noise and atomic vibrations.
The major problem is not only to find complete invariants (descriptors with \emph{no false negatives} and \emph{no false positives} for all periodic patterns) but to design efficient algorithms for distance metrics on these invariants that should continuously behave under noise.
\smallskip

The proposed continuous metrics solve this problem in any Euclidean dimension and are algorithmically approximated with small error factors in times that are explicitly bounded in the size and complexity of a given pattern. 
\smallskip

The proved Lipschitz continuity allows us to confirm all near-duplicates filtered by simpler invariants in major databases of experimental and simulated crystals.  
This practical detection of noisy duplicates will stop the artificial generation of `new' materials  from slight perturbations of known crystals.  
Several such duplicates are under investigation by five journals for data integrity.
\end{abstract}


\end{frontmatter}


\section{Continuous metric problem for periodic point sets and crystals}
\label{sec:intro}

In Euclidean geometry, periodic sets of points model all periodic crystals since any atom has a physically meaningful nucleus represented by an atomic center \cite{feynman1971feynman}.
This approach is more fundamental than using graphs with chemical bonds, which are not real sticks and only abstractly representing inter-atomic interactions depending on various thresholds for distances and angles. 
\medskip

\begin{figure}[h!]
\includegraphics[width=\textwidth]{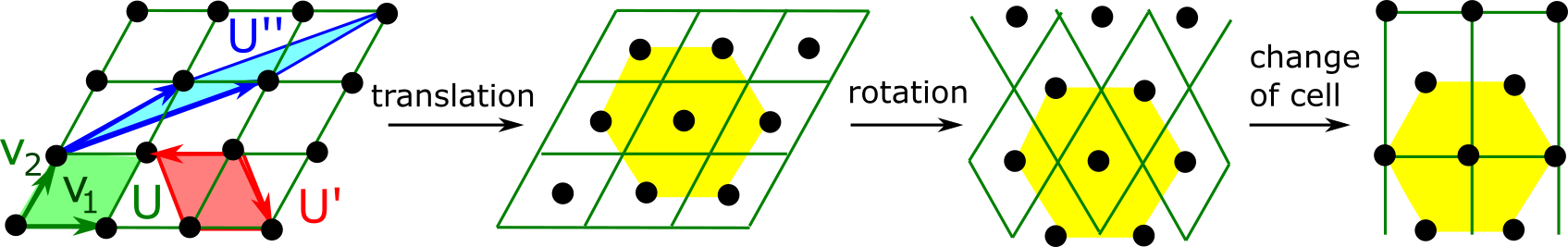}
\vspace*{-6mm}
\caption{
\textbf{Left}: three (of infinitely many) primitive cells $U,U',U''$ of the same minimal area for the hexagonal lattice $\La$. 
\textbf{Other images} show periodic sets $\La+M$ with different cells and motifs, which are all isometric to $\La$ whose hexagonal Voronoi domain is highlighted in yellow.}
\label{fig:hexagonal_lattices}
\end{figure}

A \emph{lattice} $\La\subset\R^n$ is the infinite set of all integer linear combinations $\sum\limits_{i=1}^n c_i v_i$ of a basis $v_1,\dots,v_n$ of Euclidean space $\R^n$.
Any basis defines a parallelepiped $U$ called a primitive \emph{unit cell} of $\La$.
The first picture in Fig.~\ref{fig:hexagonal_lattices} shows different unit cells in red, green, and blue, which generate the same hexagonal lattice.
\medskip

A \emph{periodic point set} $S\subset\R^n$ is a finite union of lattice translates $\La+p$ obtained from $\La$ by shifting the origin to a point $p$ from a finite \emph{motif} $M\subset U$.

This paper is motivated by the growing crisis of artificial data in crystallography \cite{anosova2024importance}
because a slight modification of a known material can be claimed as `new'.
Indeed, the fundamental question ``same or different'' \cite{sacchi2020same} was not rigorously answered for periodic crystals.
In the mathematical language, what periodic crystals can we consider equivalent, i.e. in the same class under an equivalence relation?
One classical equivalence between crystals is by symmetry, e.g. crystallographic space groups were classified into 230 types (if mirror images are distinguished) already in the 19th century by Fedorov 
and Schonflies. 

Now, all experimental databases, e.g. the Cambridge Structural Database (CSD) of 1.3+ million real materials \cite{widdowson2024continuous} need much stronger classifications than by 230 space groups or by chemical compositions.
In fact, many crystals such as diamond and graphite consist of the same elements but have vastly different properties due to essential differences in their geometric structures.
\medskip

Structures of periodic crystals are experimentally determined in a rigid form.
Hence their most practical equivalence is \emph{rigid motion}, which is a composition of translations and rotations in $\R^n$.
The slightly weaker equivalence is \emph{isometry} (any distance-preserving transformation), including mirror reflections.
\medskip

\begin{figure}[h]
\includegraphics[width=\textwidth]{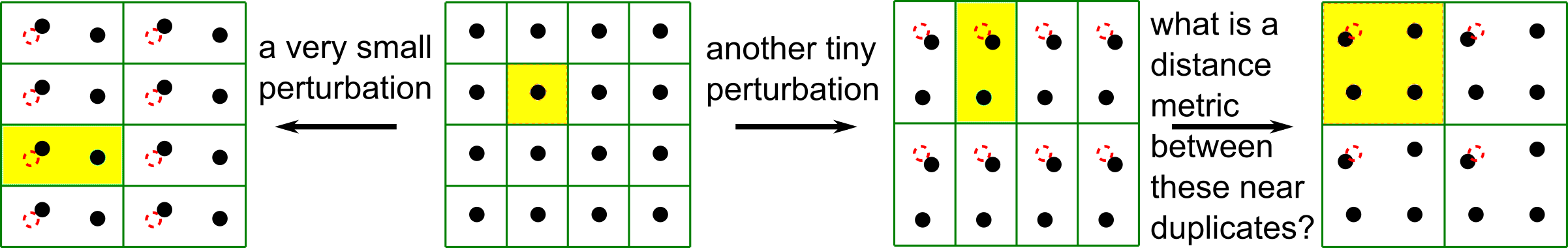}
\vspace*{-6mm}
\caption{Past descriptors based on a primitive cell cannot continuously quantify a distance between near-duplicates.
For example, under almost any perturbation, symmetry groups break down and a primitive cell volume discontinuously changes up to any integer factor.}
\label{fig:square_lattice_perturbations}
\end{figure}

An \emph{isometry class} (or pattern or or \emph{orbit} under the action of all Euclidean isometries) consists of all sets that are isometric to each other.
All isometry classes of periodic point sets form a continuously infinite space.
Indeed, almost any perturbation of points such as atomic displacements caused by noise in data produces a non-isometric crystal, which might have an arbitrarily scaled-up primitive cell, an enlarged motif, and a different symmetry group as in Fig.~\ref{fig:square_lattice_perturbations}.

Perturbations in Fig.~\ref{fig:square_lattice_perturbations} can be applied to any periodic crystal and mistakenly used for re-formatting old crystals as `new' by additionally replacing atomic types with similar ones \cite{bimler2022better}.
The first cases of geometric near-duplicates were exposed for the CSD  \cite[section~7]{widdowson2022average}, 43 materials claimed by A-lab \cite{widdowson2025geographic}, and on a much larger scale in Google's GNoME database \cite{cheetham2024artificial,widdowson2025pointwise}, see \cite[Table~1]{anosova2024importance}.
To avoid machine learning on skewed data, the filtering of near-duplicates is necessary for any database of simulated or experimental objects, e.g. point clouds, that can have infinitely many representations in different coordinate systems. 
\smallskip

A pseudo-symmetry approach, for a threshold $\ep>0$, calls periodic sets \emph{equivalent} if their cell parameters and atomic coordinates differ by at most $\ep$ \cite{zwart2008surprises}.
Then any sets can be joined by a long enough chain of $\ep$-perturbations \cite[Prop~2.1]{kurlin2024mathematics}.
If we allow any threshold $\ep>0$, the transitivity axiom (if $A\sim B\sim C$, then $A\sim C$) implies that all periodic point sets in $\R^n$ become \emph{equivalent}.
\smallskip

A mathematical approach to noisy data is to quantify perturbations by a distance metric satisfying all axioms in Definition~\ref{dfn:metric} below and taking small positive values on pairs of sets in Fig.\ref{fig:square_lattice_perturbations}, which is formalized in Problem~\ref{pro:metric}(a) .

\begin{dfn}[metric]
\label{dfn:metric}
A \emph{metric} on isometry classes of periodic sets of unordered points in $\R^n$ is a real-valued function $d$ satisfying these axioms:

\noindent
(\ref{dfn:metric}a) 
$d(S,Q)=0$ if and only if sets $S,Q$ are isometric (denoted by $S\simeq Q$);

\noindent
(\ref{dfn:metric}b) 
symmetry: $d(S,Q)=d(Q,S)$ for any periodic point sets $S,Q$ in $\R^n$;

\noindent
(\ref{dfn:metric}c) 
triangle inequality: $d(S,Q)+d(Q,T)\geq d(S,T)$ for any $S,Q,T$.
\bs 
\end{dfn}

Without the first axiom in (\ref{dfn:metric}a), even the zero function $d(S,Q)=0$ satisfies Definition~\ref{dfn:metric}a.
The atomic vibrations \cite[chapter~1]{feynman1971feynman} motivate a metric whose continuity is quantified via a maximum displacement of atoms in (\ref{pro:metric}a) below. 
 
\begin{pro}[continuous metric on periodic sets]
\label{pro:metric}
Find a metric $d$ on periodic point sets in $\R^n$ such that all the metric axioms of Definition~\ref{dfn:metric} hold and
\smallskip

\noindent
(\ref{pro:metric}a) 
$d$ is Lipschitz continuous : 
there is a constant $\la>0$ such that, for any sufficiently small $\ep>0$, if $Q$ is obtained from any periodic set $S\subset\R^n$ by perturbing each point of $S$ within its $\ep$-neighborhod, then $d(S,Q)\leq \la\ep$;
\smallskip

\noindent
(\ref{pro:metric}b) 
$d(S,Q)$ is computed or approximated (up to an explicit error factor) in a time that has a polynomially upper bound in the sizes of motifs of $S,Q$.
\bs  
\end{pro}

Problem~\ref{pro:metric} can be widened to any real data (instead of crystals) and equivalences (instead of isometry).
Condition~(\ref{pro:metric}a) goes beyond a complete classification of periodic point sets modulo isometry.
Indeed, any metric $d$ satisfying (\ref{pro:metric}a) detects all non-isometric sets $S\not\simeq Q$ by checking if $d(S,Q)\neq 0$.
Conversely, detecting an isometry $S\simeq Q$ gives only a discontinuous metric $d$ , e.g.
$d(S,Q)=1$ for any non-isometric $S\not\simeq Q$ and $d(S,Q)=0$ for any $S\simeq Q$.
\smallskip

For finite sets under isometry,
the persistent homology turned out to be weaker than anticipated \cite{smith2024generic}.
In this case, Problem~\ref{pro:metric} was solved by easier and faster invariants \cite{kurlin2024polynomial,widdowson2023recognizing}.
In the periodic case, Problem~\ref{pro:metric} was solved in dimension $n=1$ \cite{kurlin2025complete} and for lattices in $\R^2$ \cite{kurlin2024mathematics,bright2023geographic,bright2023continuous} but was open for $n>2$.
\medskip

Accuracies such as precision, recall, and F1-score make sense for finitely many classes with (usually manual) labels.
However, experimental noise (or thermal vibrations of atoms) always produce slightly different objects whose deviations should be quantified by a continuous distance metric.
Hence the real ground truth in many applications is not one of finitely many labels but an experimental structure within a continuous space of all potential objects.
Using crystals as an example, Problem~\ref{pro:metric} states the necessary conditions 
towards continuous machine learning for any data including real (non-discrete) values.
\medskip

This paper solves Problem~\ref{pro:metric} by defining a continuous metric on the complete invariant \emph{isoset} from \cite{anosova2021isometry,anosova2021introduction}.
The first step introduces a boundary tolerant metric $\BT$ on local clusters around points of a periodic set $S$, which continuously changes when points cross a cluster boundary.
This discontinuity at the boundary can be formally resolved by an extra factor, which smoothly goes down to 0 depending on an extra parameter. 
Without using extra parameters, the new metric $\BT$ will be exactly expressed in terms of simpler distances.  
\medskip

The second step uses the Earth Mover's Distance \cite{rubner2000earth} to extend $\BT$ to complete invariants \cite{anosova2021isometry} that are weighted distributions of local clusters under rotations.
The resulting metric on periodic sets in $\R^n$ is approximated with a factor $\eta$, e.g. $\eta\approx 4$ in $\R^3$, in a time depending polynomially on the input size.
\medskip

\begin{figure}[h!]
\includegraphics[height=19mm]{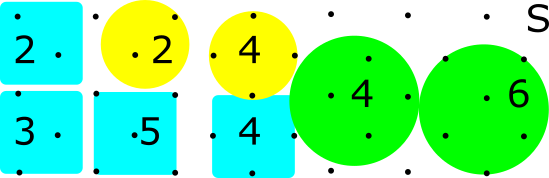}
\hspace*{2mm}
\includegraphics[height=19mm]{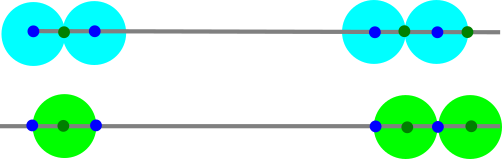}
\caption{
\textbf{Left}: for any lattice $S$ and a fixed size of a box or a ball, one can choose many non-isometric finite subsets of different sizes.
\textbf{Right}: the blue set $S$ and green set $Q$ in the line $\R$ have a small Hausdorff distance 
$d_H=\ep$ but are not related by a small perturbation. }
\label{fig:past_methods}
\end{figure}

The third step proves the metric axioms and continuity $d(S,Q)\leq 2\ep$, which also has practical importance.
Indeed, if $d(S,Q)$ is approximated by a value $d$ with a factor $\eta$, we get the lower bound $\ep\geq\frac{d}{2\eta}$ for the maximum displacement $\ep$ of points.
Such a lower bound is impossible to guarantee by analyzing only finite subsets, which can be very different in identical periodic sets, see Fig.~\ref{fig:past_methods}.

\section{Past work on distances and invariants of periodic point sets}
\label{sec:past}

This section clarifies that all past descriptors of periodic crystals are either discontinuous under perturbations as in Fig.~\ref{fig:square_lattice_perturbations} or were not proved to be complete under rigid motion.
Problem~\ref{pro:metric} was open for periodic sets for $n>2$.
\medskip

One can try comparing periodic point sets by finding an isometry of $\R^n$ that makes them as close as possible \cite{chisholm2005compack}.
This approximate matching is much easier for finite sets.
Hence it is very tempting to restrict any periodic point set to a large rectangular box  or a cube with identified opposite sides 
(a fixed 3D torus). 
However, differently located boxes or balls of any fixed size can contain non-isometric finite sets as shown in Fig.~\ref{fig:past_methods}~(left) for the square lattice.
Then extra justifications are needed to show that a comparison of periodic sets by their finite subsets does not depend on the choices of these finite subsets. 

\begin{dfn}[Hausdorff distance $d_H$, bottleneck distance $d_B$]
\label{dfn:Hausdorff+bottleneck}
\textbf{(a)}
For any sets $S,Q$ in a metric space,  $d_{\vec H}(S,Q)=\sup\limits_{p\in S} \inf\limits_{q\in Q} d(p,q)$ is
 the \emph{directed Hausdorff} distance. 
The \emph{Hausdorff distance} is $d_H(S,Q)=\max\{d_{\vec H}(S,Q),d_{\vec H}(Q,S)\}$.
\medskip

\noindent
\textbf{(b)}
The \emph{bottleneck distance} $d_B(S,Q)=\inf\limits_{g:S\to Q} \sup\limits_{p\in S}d(p,g(p))$ for sets $S,Q$ of the same cardinality is minimized over bijections $g$ and maximized over $p\in S$.
\bs
\end{dfn}

Fig.~\ref{fig:past_methods}~(right) shows the sets $S,Q$ consisting of blue and green points, respectively, where all green points of $Q$ are covered by small closed blue balls centered at all points of $S$ in the top right picture, and vice versa. 
Hence a small Hausdorff distance $d_H(S,Q)$ doesn't guarantee that the sets $S,Q$ are related by a small perturbation of points.  
A non-bijective matching of points is inappropriate for real atoms that cannot disappear and reappear from thin air.
Hence the bottleneck distance $d_B$ is more suitable for measuring atomic displacements than $d_H$.
\cite[Example~2.1]{widdowson2022resolving} shows that the 1-dimensional lattices $\Z$ and $(1+\de)\Z$ have $d_B=+\infty$ for any $\de>0$.
If any lattices have equal density (or unit cell volume), they have a finite bottleneck distance $d_B$ by \cite[Theorem~1(iii)]{duneau1991bounded}. 
\medskip

If we consider only periodic point sets $S,Q\subset\R^n$ with the same density (or unit cells of the same volume), the bottleneck distance $d_B(S,Q)$ becomes a well-defined wobbling distance \cite{carstens1999geometrical}, which is still discontinuous under perturbations by \cite[Example~2.2]{widdowson2022resolving},
see the related results for non-periodic sets in \cite{senechal1996quasicrystals}
\medskip

Another approach to comparing crystals is by Voronoi diagrams, 
which can be defined for periodic sets but remain combinatorially unstable as for finite sets.
Under almost any perturbation of basis vectors in $\R^2$, a rectangular lattice becomes generic with a hexagonal Voronoi domain.
Hence combinatorial descriptors of Voronoi domains discontinuously change under perturbations of non-generic sets as in Fig.~\ref{fig:square_lattice_perturbations}.
Geometric descriptors such as the area or volume can be continuously compared by the Hausdorff distance 
and helped define two continuous metrics between lattices in $\R^n$ \cite{mosca2020voronoi}, though their implementation sampled finitely many rotations without approximation guarantees.
\medskip

Other comparisons of periodic sets use a manually chosen number of neighbors  or a cut-off radius \cite{chisholm2005compack}. 
A reduction to a finite subset cannot provide a complete and continuous invariant of periodic sets because, under tiny perturbations, a \emph{primitive} (minimal by volume) cell can become larger than any bounded subset of a fixed size, see Fig.~\ref{fig:square_lattice_perturbations}.
One can guarantee the continuity under perturbations by extra smoothing at a fixed cut-off radius so that non-matched points covertly cross a fixed boundary, 
e.g. \cite{kawano2021classification} starts from a Gromov-Hausdorff distance between finite sets of any sizes and adds terms converging to 0 at the boundary.
The continuity was shown for three motions \cite[Fig.~3,4,5]{kawano2021classification} but the triangle inequality needs a proof, else clustering may not be trustworthy \cite{rass2024metricizing}.
\medskip

Crystallographers often compared periodic crystals by using reduced or conventional cells. 
In $\R^2$, a cell with basis vectors $\vec v_1,\vec v_2$ is \emph{reduced} if $|\vec v_1|\leq|\vec v_2|$ and $-\frac{1}{2}\vec v_1^2\leq \vec v_1\cdot\vec  v_2\leq 0$.
The vectors $\vec v_1=(2a,0)$ and $\vec v_2^{\pm}=(-a,\pm b)$ for $b\geq a\sqrt{3}$ and both signs $\pm$ are reduced and define isometric lattices related by reflection.
This ambiguity of bases can be resolved by an additional condition $\det(\vec v_1,\vec v_2)>0$, which creates the inevitable discontinuity, see more details in \cite[Fig.~4]{kurlin2024mathematics}. 
In $\R^3$, the most widely used reduced cell is Niggli's cell, 
which has a minimum volume and all angles as close to $90^\circ$ as possible.
Niggli's cell was known to be experimentally discontinuous since 1965 \cite{lawton1965reduced} or even earlier 
due to Fig.~\ref{fig:square_lattice_perturbations}. 
\medskip

In $\R^3$, all generic periodic sets are distinguished by density functions \cite{edels2021}, which can be computed at discrete values of a continuous radius $t\in\R$ \cite{smith2022practical}.
A metric between density functions was defined in terms of suprema over infinitely many $t\in\R$, so the metric was approximated without guarantees.
The density functions \cite{anosova2023density} coincide for the periodic sets $S_{15}=X+Y+15\Z$, $Q_{15}=X-Y+15\Z$, where $X = \{0, 4, 9\}$, $Y = \{0, 1, 3\}$ \cite[Example~11]{anosova2022density}.
This pair and all generic periodic sets are distinguished by a faster Pointwise Distance Distribution (PDD) \cite[Theorem~4.4]{widdowson2022resolving} whose averages $\AMD$ are incomplete by \cite[Example~3.3]{widdowson2022resolving}. The PDD cannot distinguish any mirror images.
\medskip

A distance between invariant values can be a metric on isometry classes only if the underlying invariant is complete under isometry.
Otherwise, non-isometric sets can have identical invariant values with a distance of 0.  
Hence a complete classification should take into account a potential high complexity of periodic sets.
Inspired by \cite{delone1976local,dolbilin1998multiregular}, the isometry classification of periodic sets was reduced \cite{anosova2021isometry} to only rotations of local clusters whose radius can be determined from $S$.
\medskip

Section~\ref{sec:isosets} reminds us of a complete invariant isoset from \cite{anosova2021isometry}.
Section~\ref{sec:metrics} introduces a Lipschitz continuous metric (Definition~\ref{dfn:EMD_isosets} and Theorem~\ref{thm:continuity}), whose polynomial time bounds (Corollaries~\ref{cor:compare_isosets}, \ref{cor:approximate_EMD}) are proved in section~\ref{sec:algorithms}.
Section~\ref{sec:lower_bound} proves a lower bound (Theorem~\ref{thm:lower_bound}) for the new metric via faster $\PDD$s.
Section~\ref{sec:discussion} discusses the significance of the continuous metric for detecting near-duplicates in major crystal databases and for upholding scientific integrity.

\section{Isometry classification of periodic point sets by complete invariants}
\label{sec:isosets}

This section reviews the complete invariant \cite{anosova2021isometry}
based on local clusters and their symmetry groups, which were previously studied in \cite{delone1976local,dolbilin1998multiregular}.

\begin{dfn}[\emph{global} clusters and \emph{$m$-regular} periodic sets]
\label{dfn:m-regular}
For any point $p$ in a periodic set $S\subset\R^n$, the \emph{global cluster} is $C(S,p)=\{\vec q - \vec p \vl q\in S\}$.
For any $p,q\in\R^n$, let the set $\Or(\R^n;p,q)$ consist of all isometries of $\R^n$ that map $p$ to $q$.
Global clusters $C(S,p)$ and $C(S,q)$ are called \emph{isometric} if there is $f\in\Or(\R^n;p,q)$ such that $f(S)=S$. 
A periodic point set $S\subset\R^n$ is called \emph{$m$-regular} if all global clusters of $S$ form exactly $m\geq 1$ isometry classes.
\bs
\end{dfn}

For any point $p\in S$, its global cluster is a view of $S$ from the position of a point $p$.
We view all astronomical stars in the universe $S$ from our planet Earth at $p$.
Any lattice is 1-regular since all its global clusters are related by translations.
Though global clusters $C(S,p),C(S,q)$ at any different points $p,q\in S$ contain the same set $S$, they may not match under the translation shifting $p$ to $q$.
The global clusters are infinite, hence distinguishing them under isometry is not easier than original periodic sets.
However, the $m$-regularity of a periodic set can be checked in terms of finite local $\al$-clusters below.

\begin{dfn}[\emph{local $\al$-clusters} $C(S,p;\al)$ and \emph{symmetry groups} $\sym(S,p;\al)$]
\label{dfn:local_cluster}
For a point $p$ in a periodic point set $S\subset\R^n$ and any $\al\geq 0$, the local \emph{$\al$-cluster} $C(S,p;\al)$ is the set of all vectors $\vec q - \vec p$ such that $q\in S$ and $|\vec q - \vec p|\leq\al$.
Let the group $\Or(\R^n;p)$ consist of all isometries that fix $p$.
If $p=0$ is the origin, $\Or(\R^n;0)$ is the usual orthogonal group.
The \emph{symmetry} group $\sym(S,p;\al)$ consists of all isometries $f\in\Or(\R^n;p)$ that map $C(S,p;\al)$ to itself so that $f(p)=p$.
\bs
\end{dfn}

For any periodic set $S$, if $\al$ is smaller than the minimum distance between all points of $S$, then any $\al$-cluster $C(S,p;\al)$ is one point $\{p\}$.
Its symmetry group  consists of all isometries fixing the center $p$, so $\sym(S,p;\al)=\Or(\R^n;p)$.
When $\al$ is increasing, the $\al$-clusters $C(S,p;\al)$ become larger and there can be fewer (not more) isometries $f\in\Or(\R^n;p)$ that bijectively map $C(S,p;\al)$ to itself.
So the group $\sym(S,p;\al)$ can become smaller (not larger) and eventually stabilizes (stops changing), which is formalized in Definition~\ref{dfn:stable_radius}.
This stabilization uses the bridge length extending the idea of a longest edge in a Minimum Spanning Tree to a periodic set $S$, which does not easily reduce to the finite case \cite{mcmanus2025computing}.

\begin{dfn}[\emph{bridge length} $\be(S)$]
\label{dfn:bridge_length}
For a periodic point set $S\subset\R^n$, the \emph{bridge length} is a minimum distance $\be(S)>0$ such that any $p,q\in S$ can be connected by a sequence of points $p_0=p,p_1,\dots,p_k=q$ such that 
any two successive points $p_{i-1},p_{i}$ are close so that $|\vec p_{i-1} - \vec p_{i}|\leq\be(S)$ for $i=1,\dots,k$.
\bs
\end{dfn}

\cite[Theorem~1.3]{dolbilin1998multiregular} described how a family of clusters determines a periodic point set under isometry.
These results motivated the \emph{isotree}, \emph{stable} radius, and \emph{isoset} in Definitions~\ref{dfn:isotree}, \ref{dfn:stable_radius},~\ref{dfn:isoset}, respectively, leading to the isometry classification of periodic point sets via isosets in Theorem~\ref{thm:isoset_complete}.
The \emph{isotree} in Definition~\ref{dfn:isotree} is inspired by a clustering dendrogram because points of $S$ split into isometry classes of $\al$-clusters at variable radii $\al$, not at a fixed $\al$.

\begin{dfn}[\emph{isotree} $\IT(S)$ of \emph{$\al$-partitions}]
\label{dfn:isotree}
Fix a periodic point set $S\subset\R^n$.
Points $p,q\in S$ are \emph{$\al$-equivalent} if their $\al$-clusters $C(S,p;\al)$ and $C(S,q;\al)$ can be related by an isometry that matches their centers.
The \emph{isometry class} $[C(S,p;\al)]$ consists of all $\al$-clusters isometric to $C(S,p;\al)$.
The \emph{$\al$-partition} $P(S;\al)$ is the splitting of $S$ into $\al$-equivalence classes of points.
Call a value $\al$ \emph{singular} if $P(S;\al)\neq P(S;\al-\ep)$ for any small enough $\ep>0$.
Represent each $\al$-equivalence class by a vertex of the \emph{isotree} $\IT(S)$.
The top vertex of $\IT(S)$ represents the $0$-equivalence class coinciding with $S$.
For any successive singular values $\al<\al'$, connect the vertices representing any classes $A\in P(S;\al)$ and $A'\in P(S;\al')$ such that $A'\subset A$ by an edge of the length $\al'-\al$ in $\IT(S)$.
\bs
\end{dfn}

\vspace*{-2mm}
\begin{figure}[h]
\includegraphics[width=\textwidth]{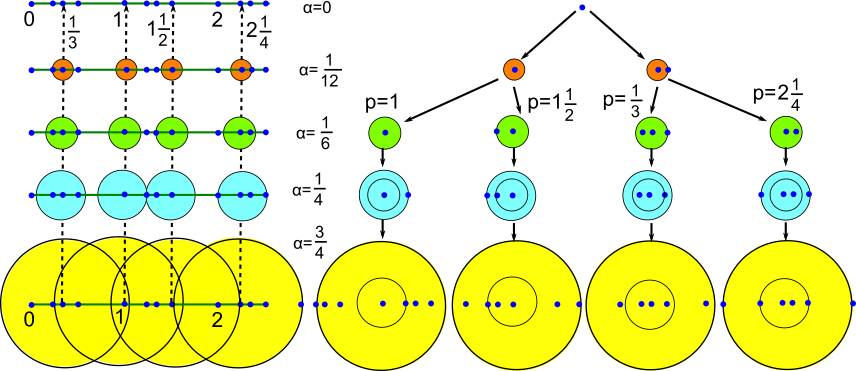}
\vspace*{-8mm}
\caption{\textbf{Left}: the 1-dimensional set $S_4=\{0,\frac{1}{4},\frac{1}{3},\frac{1}{2}\}+\Z$ has four points in the unit cell $[0,1)$ and is 4-regular by Definition~\ref{dfn:m-regular}.
\textbf{Right}: the colored disks show $\al$-clusters in the line $\R$ with radii $\al=0,\frac{1}{12},\frac{1}{6},\frac{1}{4},\frac{3}{4}$ and represent points in the isotree $\IT(S_4)$ from Definition~\ref{dfn:isotree}.
}
\label{fig:4-regular_set_isotree}
\end{figure}

For any periodic point set $S\subset\R^n$, the root vertex of $\IT(S)$ at $\al=0$ is the single class $S$, because any 0-cluster $C(S,p;0)$ of a point $p\in S$ consists only of its center $p$.
When the radius $\al$ is increasing, $\al$-clusters $C(S,p;\al)$ include more points and hence may not be isometric.
In other words, any $\al$-equivalence class from $P(S;\al)$ may split into two or more classes, which cannot merge at any larger $\al'$.
Branched vertices of $\IT(S)$ correspond to the values of $\al$ when an $\al$-equivalence class is split into subclasses for $\al'$ slightly larger than $\al$.
So the number $|P(S;\al)|$ of $\al$-equivalence is non-decreasing in $\al$, see Fig.~\ref{fig:4-regular_set_isotree}.
\medskip

The $\al$-clusters of the 1-dimensional periodic point set $S_4\subset\R$ in Fig.~\ref{fig:4-regular_set_isotree} are intervals in $\R$, shown as disks for better visibility.
In Fig.~\ref{fig:4-regular_set_isotree}, this class persists until $\al=\frac{1}{12}$, when all points $p\in S_4$ are split into two classes: one represented by 1-point cluster $\{p\}$ for $p\in\{0,\frac{1}{2}\}+\Z$, and another represented by 2-point clusters $\{p,p+\frac{1}{12}\}$, $p\in\{\frac{1}{4},\frac{1}{3}\}+\Z$.
The periodic set $S_4$ has four $\al$-equivalence classes for any radius $\al\geq\frac{1}{6}$.
For any point $p\in\Z\subset S_4$, the symmetry group $\sym(S_4,p;\al)=\Z_2$ is generated by the reflection in $p$ for $\al\in[0,\frac{1}{4})$.
For all $p\in S_4$, the symmetry group $\sym(S_4,p;\al)$ is trivial for any $\al\geq\frac{1}{4}$.
For any periodic point set $S\subset\R^n$, the $\al$-partitions of $S$ stabilize in the sense below.

\begin{dfn}[the \emph{minimum stable radius} $\al(S)$]
\label{dfn:stable_radius}
Let $S\subset\R^n$ be a periodic point, $\be\geq\be(S)$ be an upper bound of the bridge length $\be(S)$ from Definition~\ref{dfn:bridge_length}.
A radius $\al\geq\be$ is called \emph{stable} if the following conditions hold:
\medskip

\noindent
(\ref{dfn:stable_radius}a) 
the $\al$-partition $P(S;\al)$ equals the $(\al-\be)$-partition $P(S;\al-\be)$;
\medskip

\noindent
(\ref{dfn:stable_radius}b)
the groups stabilize so that
$\sym(S,p;\al)=\sym(S,p;\al-\be)$ for any $p\in S$, i.e.
any isometry $f\in\sym(S,p;\al-\be)$ preserves the larger cluster $C(S,p;\al)$.
\medskip

\noindent
A minimum value of a stable radius $\al$ satisfying (\ref{dfn:stable_radius}ab) for $\be=\be(S)$ from Definition~\ref{dfn:bridge_length} is called \emph{the minimum stable radius} and denoted by $\al(S)$. 
\bs
\end{dfn}

Due to the upper bounds in Lemma~\ref{lem:upper_bounds}(b,c), the minimum stable radius $\al(S)\geq 0$ exists and is achieved because $P(S;\al)$ and $\sym(S,p;\al)$ are continuous on the right (unchanged when $\al$ increases by a sufficiently small value). 
\medskip

Any $m$-regular periodic point set $S\subset\R^n$ has at most $m$ $\al$-equivalence classes, so the isotree $\IT(S)$ stabilizes with maximum $m$ branches.
Though (\ref{dfn:stable_radius}b) is stated for all points $p\in S$ for simplicity, it suffices to check condition~(\ref{dfn:stable_radius}b) for points only from a finite motif $M$ of $S$ due to periodicity.
\medskip

All stable radii of $S$ form the interval $[\al(S),+\infty)$ by
Lemma~\ref{lem:stable_radius} in the appendix.
The periodic set $S_4$ in Fig.~\ref{fig:4-regular_set_isotree} has $\be(S_4)=\frac{1}{2}$ and $\al(S)=\frac{3}{4}$ since the $\al$-partition and symmetry groups $\sym(S_4,p;\al)$ are stable for $\frac{1}{4}\leq\al\leq\frac{3}{4}$.
\medskip

Condition (\ref{dfn:stable_radius}b) doesn't follow from condition (\ref{dfn:stable_radius}a) due to the following example. 
Let $\La$ be the 2D lattice with the basis $(1,0)$ and $(0,\be)$ for $\be>1$.
Then $\be$ is the bridge length of $\La$. 
Condition (\ref{dfn:stable_radius}a) is satisfied for any $\al\geq 0$, because all points of any lattice are equivalent under translations. However, condition (\ref{dfn:stable_radius}b) fails for any $\al<\be+1$.
Indeed, the $\al$-cluster of the origin $(0,0)$ contains five points $(0,0),(\pm 1,0),(0,\pm\be)$, whose symmetries are generated by the two reflections in the axes $x,y$, but the $(\al-\be)$-cluster of the origin $(0,0)$ consists of its center  and has the symmetry group $\Or(\R^2)$. 
It is possible that condition (\ref{dfn:stable_radius}b) might imply condition (\ref{dfn:stable_radius}a), but in practice it makes sense to verify (\ref{dfn:stable_radius}b) only after checking much simpler condition (\ref{dfn:stable_radius}a). Both conditions are essentially used in the proof of Isometry Classification Theorem~\ref{thm:isoset_complete}. 
\medskip

Conditions (\ref{dfn:stable_radius}ab) appeared in \cite{dolbilin1998multiregular} with different notations $\rho,\rho+t$.
Since many applied papers use $\rho$ for the physical density and have many types of bond distances, we replaced $t$ and $\rho+t$ with the bridge length $\be$ and radius $\al$, respectively, as for growing $\al$-shapes in Topological Data Analysis \cite{smith2024generic}. 
\medskip

Recall that the \emph{covering radius} $R(S)$ of a periodic point set $S\subset\R^n$ is the minimum radius $R>0$ such that $\bigcup\limits_{p\in S}\bar B(S;R)=\R^n$, or the largest radius of an open ball in the complement $\R^n\setminus S$.
For $m$-regular point sets in $\R^n$, an upper bound of $\al(S)$ can be extracted from \cite[Theorem 1.3]{dolbilin1998multiregular} whose proof motivated a stronger bound in Lemma~\ref{lem:upper_bounds}(c), see comparisons in Example~\ref{exa:upper_bounds}(c). 
\medskip

A periodic point set $S$ is \emph{locally antipodal} if the local cluster $C(S,p;2R(S))$ is centrally symmetric for any point $p\in S$, i.e. bijectively maps to itself under $\vec q\mapsto 2\vec p-\vec q$, $q\in\R^n$.
\cite[Theorem~1]{dolbilin2016uniqueness} says that all locally antipodal Delone sets, hence all periodic sets $S$, are globally antipodal, i.e. $S$ is preserved under the isometry $\vec q\mapsto 2\vec p-\vec q$ for any fixed $p\in S$, e.g. any lattice is antipodal.

\begin{lem}[upper bounds for a stable radius $\al(S)$ and bridge length $\be(S)$]
\label{lem:upper_bounds}
\textbf{(a)}
Let $S\subset\R^n$ be a periodic point set with a unit cell $U$, which has the longest edge $b$ and longest diagonal $d$.
Set $r(U)=\max\{b,\frac{d}{2}\}$.
Then the bridge length $\be(S)$ from Definition~\ref{dfn:bridge_length} 
has the upper bound $\min\{2R(S),r(U)\}\geq\be(S)$.
\medskip

\noindent
\textbf{(b)}
For any antipodal periodic set $S\subset\R^n$ whose covering radius is $R(S)$, the minimum stable radius has the upper bound $2R(S)+\be(S)>\al(S)$.  
\medskip

\noindent
\textbf{(c)}
Let $S\subset\R^n$ be any periodic point set with the bridge length $\be$.
For any point $p\in S$ and a radius $\al_0\geq 2R(S)$, the order $|\sym(S,p;\al_0)|$ of the group $\sym(S,p;\al_0)$ should be finite.
Let $p_1,\dots,p_m\in S$ be all points of an asymmetric unit of $S$.
Set $L=\left[\sum\limits_{i=1}^{m}\big(\log_2|\sym(S,p_i;\al_0)|-\log_2|\sym(S,p_i)|\big)\right]$. 
Then the minimum stable radius $\al(S)$ from Definition~\ref{dfn:stable_radius} 
has the upper bound $\al_0+(L+m)\be\geq\al(S)$.
If $\al_0=2R(S)$, then $(L+m+1)2R(S)\geq\al(S)$. 
\bs
\end{lem}
\begin{proof}
\textbf{(a)}
The lemma in \cite[section~2]{delone1976local} proved that, in any Delone set $S$ with the covering radius $R(S)$, any two points $p,q\in S$ can be connected by a finite sequence of points $p_0=p,p_1,\dots,p_k=q$ such that $|\vec p_{i-1} - \vec p_{i}|\leq 2R(S)$ for $i=1,\dots,k$.
In particular, any periodic point set $S$ has the upper bound $2R(S)\geq\be(S)$.
It remains to prove the second upper bound $r(U)\geq \be(S)$.
\smallskip

For a point $p\in S$, shift the unit cell $U$ so that $p$ becomes the origin of $\R^n$ and a vertex of $U$, so the lattice $\La$ can be considered a subset of the periodic point set $S$.
Any points of $\La$ can be connected by a sequence of lattice points such that any successive points have a distance not greater than the longest edge-length $b$ of $U$.
Any point of a motif $M\subset U$ of $S$ is at most $r(U)$ away from a vertex of $U$, where $d$ is the length of the longest diagonal of $U$.
Any points of $S$ can be connected by a sequence whose successive points are at most $r(U)=\max\{b,\frac{d}{2}\}$ away from each other, so $\be(S)\leq r(U)$ by Definition~\ref{dfn:bridge_length}.
\medskip

\noindent
\textbf{(b)}
We will prove that the conditions of Definition~\ref{dfn:stable_radius} hold for $\al=2R(S)+\be(S)$ and $\be=\be(S)$.
To prove condition~\ref{dfn:stable_radius}(a), we check below that any $2R(S)$-equivalent points $p,q\in S$ are $\al$-equivalent for any $\al>2R(S)$.
The $2R(S)$-equivalence means that there is an isometry $f\in\Or(\R^n;p,q)$ such that $f(C(S,p;2R(S)))=C(S,q;2R(S))$.
Set $Q=f(S)$. 
Then $f(C(S,p;2R(S)))=C(f(S),f(p);2R(S))$ means that $C(S,q;2R(S))=C(Q,q;2R(S))$.
\cite[Theorem~3]{dolbilin2016uniqueness} implies that if antipodal periodic point sets $S,Q\subset\R^n$ have a common point $q$ with $C(S,q;2R(S))=C(Q,q;2R(S))$, then $S=Q$.
In our case, $f(S)=S$ implies that $f$ makes the points $p$ and $q=f(p)$ $\al$-equivalent for any $\al>2R(S)$.
Condition~\ref{dfn:stable_radius}(b) says that 
any isometry $f\in\sym(S,p;2R(S))$ should belong to $\sym(S,p;\al)$ for any point $p\in S$ and radius $\al>2R(S)$.
Indeed, \cite[Theorem~3]{dolbilin2016uniqueness} implies that $Q=f(S)$ and $S$ should coincide, so $f$ isometrically maps any cluster $C(S,p;\al)$ to itself, hence $f\in\sym(S,p;\al)$.
\medskip

\noindent
\textbf{(c)}
Lemma~\ref{lem:finite_symmetry}, which was briefly proved in \cite[p.~20]{delone1976local}, says that the symmetry group $\sym(S,p;2R(S))$ is finite.
For any initial radius $\al_0\geq 2R(S)$, we aim to find a radius $\al=\al_0+k\be$ such that both conditions~\ref{dfn:stable_radius}(a,b) hold for a suitable index $k=1,2,3,\dots$ whose upper bound we will determine below.
\medskip

If condition~\ref{dfn:stable_radius}(a) fails for some $\al=\al_0+k\be$, the number $|P(S;\al_0+(k-1)\be)|$ of $\al$-equivalence classes increases at least by one when $\al_0+(k-1)\be$ increases to $\al_0+k\be$.
Since an asymmetric unit of $S$ consists of $m\geq 1$ points, there are at most $m-1$ incremental values $0=k_0\leq k_1\leq \dots\leq k_{m-1}$ when $1\leq|P(S;\al_0+(k_i-1)\be)|<|P(S;\al_0+k_i\be)|\leq m$ for $i=1,\dots,m-1$.
\smallskip

In a degenerate case, if all points of $S$ are $(\al_0+(k-1)\be)$-equivalent, this single class can split into the maximum $m>1$ classes of $(\al_0+k\be)$-equivalence, then $k_1=\dots=k_{m-1}=k\geq 1$.
For any successive incremental values $k_{i-1}<k_{i}$, the number $|P(S;\al_0+k\be)|$ of $(\al_0+k\be)$-equivalence classes is constant for $k=k_{i-1}+1,\dots,k_{i}$, so condition~\ref{dfn:stable_radius}(a) holds for every radius $\al=\al_0+k\be$.
\smallskip

By reordering the points $p_1,\dots,p_m$ from an asymmetric unit of $S$, we can assume that $p_1,\dots,p_{i}$ represent $i$ classes of $(\al+k_{i-1}\be)$-equivalence for any fixed $i=1,\dots,m$.
Set $L(k)=\sum\limits_{i=1}^{m}\log_2|\sym(S,p_i;\al_0+k\be)|$.
When $k$ increases, any group $\sym(S,p_i;\al_0+k\be)$ can become only smaller, not larger, so $L(k)$ is non-increasing. 
If $L(k-1)=L(k)$ for any $0<k\neq k_1,\dots,k_{m-1}$, both conditions~\ref{dfn:stable_radius}(a,b) hold, so $\al_0+k\be$ is a stable radius.
We will find an upper bound for a minimum value of such $k$.
If condition~\ref{dfn:stable_radius}(b) fails for all radii $\al=\al_0+k\be$ with $k=k_{i-1}+1,\dots,k_{i}$, then at least one of the groups $\sym(S,p;\al_k+j\be)$ for $p\in\{p_1,\dots,p_{i}\}$ is a proper subgroup of $\sym(S,p;\al_0+(k-1)\be)$.
The order of a proper subgroup is at most a half of the order of the group, so $$\log_2|\sym(S,p;\al_0+k\be)|\leq \log_2|\sym(S,p;\al_0+(k-1)\be)|-1, k=k_{i-1}+1,\dots,k_{i}.$$
Hence the sum $L(k)$ decreases at least by 1 for any failure of condition~\ref{dfn:stable_radius}(b)
from $L(0)=\sum\limits_{i=1}^{m}\log_2|\sym(S,p_i;\al_0)|$
to $L(+\infty)=\sum\limits_{i=1}^{m}\log_2|\sym(S,p_i)|$, where $\sym(S,p_i)$ is the symmetry group of the global cluster $C(S;p_i)$.
Adding $m-1$ potential failures of condition~\ref{dfn:stable_radius}(a) for $\al_0+k_i\be$ with $i=1,\dots,m-1$, the radius $\al_0+k\be$ cannot be stable 
for a maximum $L+m-1$ values of $k$, where $$L=[L(0)-L(+\infty)]=\left[\sum\limits_{i=1}^{m}\big(\log_2|\sym(S,p_i;\al_0)|-\log_2|\sym(S,p_i)|\big)\right].$$
Then any $\al=\al_0+k\be$ with $k\geq L+m$ is stable, so 
$\al(S)\leq \al_0+(L+m)\be$.
To get $\al(S)\leq (L+m+1)2R(S)$, set $\al_0=2R(S)$ and use $\be\leq 2R(S)$ from~(a).
\end{proof}

The upper bound in Lemma~\ref{lem:upper_bounds}(a) holds for any unit cell of $S$.
If a cell is non-reduced and too long, its reduced form can have smaller bounds for $\be(S)$.   

\begin{exa}[upper bounds for $\al(S)$ and $\be(S)$]
\label{exa:upper_bounds}
Let $\La(b)\subset\R^n$ be a lattice whose unit cell is a rectangular box with the longest edge $b\geq 1$.
\smallskip

\noindent
\textbf{(a)}
In Lemma~\ref{lem:upper_bounds}(a), the upper bound $b\geq\be(S)$ is tight because $\be(\La(b))=b$.
\smallskip

\noindent
\textbf{(b)}
In Lemma~\ref{lem:upper_bounds}(b), the ratio $(2R(S)+\be(S))/\al(S)\geq 1$ tends to $1$ as $b\to+\infty$ for any fixed $n$.
Indeed, a cluster $C(\La(b),0;\al)$ is $n$-dimensional only for $\al\geq b$, so the group $\sym(\La(b),0;\al)$ stabilizes at $\al=b$, hence $\al(S)=b+\be(\La(b))=2b$ is the minimum stable radius. 
The covering radius $R(\La(b))$ is half of the longest diagonal of the rectangular cell $U$.
If $b\to+\infty$ and all other sizes of $U$ remain fixed, the ratio $(2R(\La(b))+\be(\La(b)))/\al(S)$ tends to $1$ for any fixed $n$.
\medskip

\noindent
\textbf{(c)}
Lemma~\ref{lem:upper_bounds}(c) was motivated by \cite[Theorem~1.3]{dolbilin1998multiregular}, which implies the upper bound $\be(S)+2m(n^2+1)\log_2(2+R(S)/r(S))>\al(S)$ for $m$-regular point sets.
Let $\La\subset\R^2$ be a lattice whose unit cell is a rhombus with sides 1.
Then $m=1$, $n=2$, $r(\La)=0.5$, $\be(\La)=1$, and $\al(\La)=2$.
If $\La$ deforms from a square lattice to a hexagonal lattice, the covering radius $R(\La)$ varies in the range $[\frac{1}{\sqrt{3}},\frac{1}{\sqrt{2}}]$.
The past bound above gives the estimate $1+2(2^2+1)\log_2(2+\frac{2}{\sqrt{3}})\approx 17.6>\al(\La)=2$.  
For any lattice $\La$ in this family, the symmetry group $\sym(\La,0)=\sym(\La,0;1)$ stabilizes at $\al_0=1$. 
Lemma~\ref{lem:upper_bounds}(c) for $\al_0=1$ gives $L=\log_2(2)-\log_2(2)=0$, so the upper bound $\al_0+(L+m)\be(S)\geq\al(S)$ is tight: $2\geq\al(\La)$. 
In practice, if $L$ is large because some local clusters $C(S;p;\al_0)$ have too many symmetries, one can increase the radius $\al_0$ to reduce $L$ for a better bound of $\al(S)$. 
\bs
\end{exa}

Definition~\ref{dfn:isoset} reminds of the \emph{isoset}, which was initially introduced in \cite[Definition~9]{anosova2021isometry}. 
We also cover the case of rigid motion 
and prove Completeness Theorem~\ref{thm:isoset_complete} in the appendix in more detail than in \cite[Theorem~9]{anosova2021isometry}.
  
\begin{dfn}[isoset $I(S;\al)$ at a radius $\al\geq 0$]
\label{dfn:isoset}
Let a periodic point set $S\subset\R^n$ have a motif $M$ of $m$ points.
Split all points $p\in M$ into $\al$-equivalence classes.
Each $\al$-equivalence class of (say) $k$ points in $M$ can be associated with the \emph{isometry class} $\si=[C(S,p;\al)]$ of an $\al$-cluster centered at some $p\in M$.
The \emph{weight} of $\si$ is $w=k/m$.
The \emph{isoset} $I(S;\al)$ is the unordered set of all isometry classes $(\si;w)$ with weights $w$ for all points $p$ in the motif $M$.
If we replace isometry with rigid motion, we get the \emph{oriented} isoset $I^o(S;\al)$.
\bs
\end{dfn}

All points $p$ of a lattice $\La\subset\R^n$ from one $\al$-equivalence class for any radius $\al\geq 0$ because all $\al$-clusters $C(\La,p;\al)$ are isometrically equivalent to each other by translations. 
Hence the isoset $I(\La;\al)$ is one isometry class of weight 1 for $\al\geq 0$, see examples in Fig.~\ref{fig:square_vs_hexagon}.
All isometry classes $\si$ in $I(S;\al)$ are in a 1-1 correspondence with all $\al$-equivalence classes in the $\al$-partition $P(S;\al)$ from Definition~\ref{dfn:isotree}.
So $I(S;\al)$ without weights can be viewed as a set of points in the isotree $\IT(S)$ at the radius $\al$.
The size of the isoset $I(S;\al)$ equals the number $|P(S;\al)|$ of $\al$-equivalence classes in the $\al$-partition. 
Formally, $I(S;\al)$ depends on $\al$ because $\al$-clusters grow in $\al$.
To distinguish any $S,Q\subset\R^n$ under isometry, we will compare their isosets at a maximum stable radius of $S,Q$.

\begin{exa}[isosets of simple lattices]
\label{exa:isosets_lattices}
\textbf{(a)}
Any lattice $\La\subset\R^n$ is 1-regular by Definition~\ref{dfn:m-regular} and can be assumed to contain the origin $0$ of $\R^n$.
Then the isoset $I(\La;\al)$  consists of a single isometry class of a cluster $C(\La,0;\al)$.
So the isotree $\IT(\La)$ is a linear path, which is horizontally drawn for the hexagonal and square lattices $\La_6,\La_4$ in Fig.~\ref{fig:lattice_isotree}.
If both $\La_6,\La_4$ have a minimum inter-point distance 1, then the bridge length from Definition~\ref{dfn:bridge_length} is $\be=1$.
\medskip

\begin{figure}[ht]
\includegraphics[height=16mm]{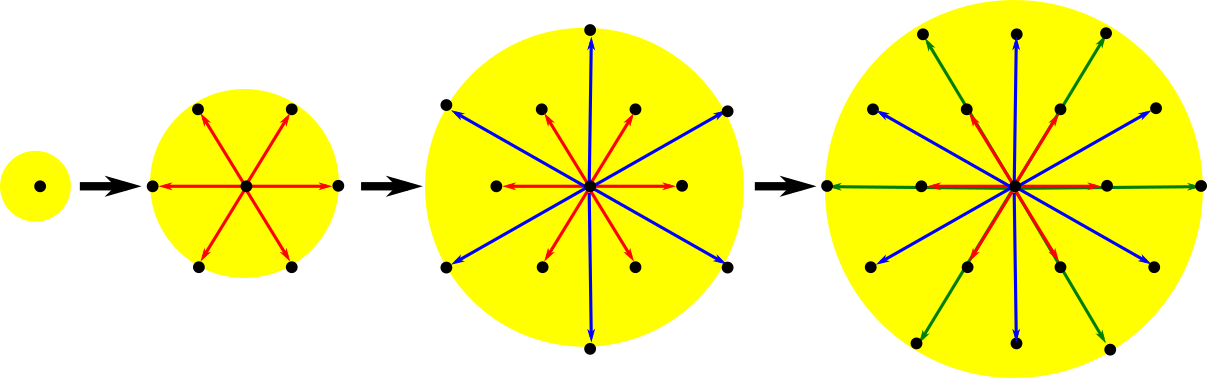}
\hspace*{1mm}
\includegraphics[height=16mm]{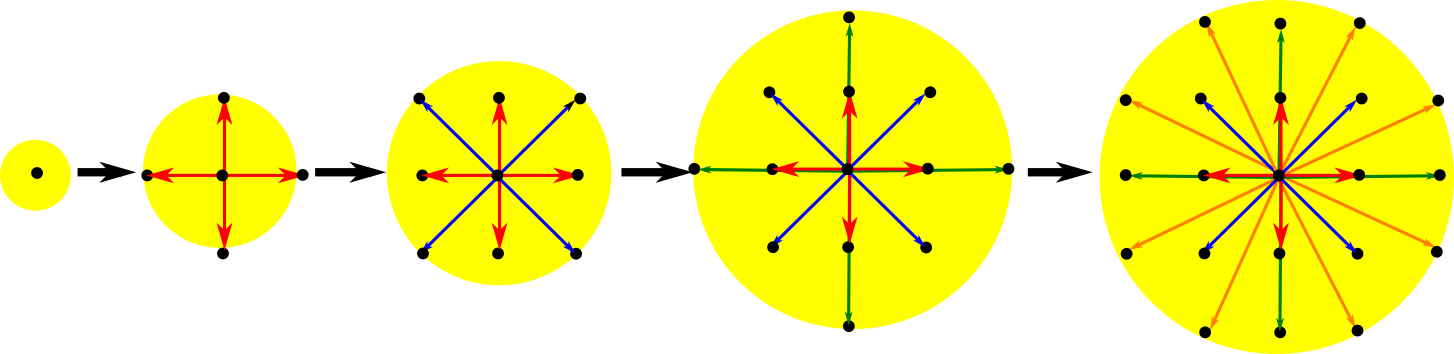}
\vspace*{-4mm}
\caption{
The isotree of any lattice $\La$ is $[0,+\infty)$ is a line $\R$ parametrized by the radius $\al$.
\textbf{Left}: the isotree of the hexagonal lattice $\La_6$.
\textbf{Right}: the isotree of the square lattice $\La_4$.
}
\label{fig:lattice_isotree}
\end{figure}

\noindent
\textbf{(b)}
For the hexagonal lattice $\La_6\subset\R^2$, 
$C(\La_6,(0,0);\al)$ includes points $p\neq(0,0)$ only for $\al\geq 1$.
The cluster $C(\La_6,(0,0);1)=\{(0,0),(\pm 1,0),(\pm\frac{1}{2},\pm\frac{\sqrt{3}}{2})\}$ appears in the 2nd step of  Fig.~\ref{fig:lattice_isotree}~(left).
The symmetry group $\sym(\La_6,(0,0);\al)$ becomes the dihedral group $D_6$ (all symmetries of a regular hexagon) for $\al\geq 1$.
Hence any $\al\geq\be+1=2$ is stable.
The isoset $I(\La_6;1)$ is the isometry class of the cluster $C(\La_6,(0,0);1)$ of six vertices of the regular hexagon and its center.
\medskip

\noindent
\textbf{(c)}
For the square lattice $\La_4\subset\R^2$, $C(\La_4,(0,0);\al)$ has points $p\neq(0,0)$ only for $\al\geq 1$.
 $C(\La_4,(0,0);2)=\{(0,0),(\pm 1,0),(0,\pm 1),(\pm\sqrt{2},\pm\sqrt{2}),(\pm 2,0),(0,\pm 2)\}$ includes the origin $(0,0)$ with its 12 neighbors in the 4th step of Fig.~\ref{fig:lattice_isotree}~(right).
The group $\sym(\La_4,(0,0);\al)$ becomes the dihedral group $D_4$ (all symmetries of a square) for $\al\geq 1$.
So any $\al\geq\be+1=2$ is stable.
The isoset $I(\La_4;1)$ is the isometry class of $C(\La_4,(0,0);1)$ of four vertices of the square and its center.
\bs
\end{exa}

An equality $\si=\xi$ between isometry classes of clusters means that some (hence any) clusters $C(S,p;\al)$ and $C(Q,q;\al)$ representing $\si,\xi$, respectively, are related by $f\in\Or(\R^n;p,q)$, which will be algorithmically tested in Corollary~\ref{cor:compare_isosets}.

\begin{thm}[isometry classification of periodic point sets]
\label{thm:isoset_complete}
For any periodic point sets $S,Q\subset\R^n$, let $\al$ be a common stable radius satisfying Definition~\ref{dfn:stable_radius} for an upper bound $\be\geq\be(S),\be(Q)$.
Then $S,Q$ are isometric (related by rigid motion, respectively) if and only if there is a bijection $\ph:I(S;\al)\to I(Q;\al)$ (between oriented isosets, respectively) that preserves all their weights.
\bs
\end{thm}

Theorem~\ref{thm:isoset_complete} was inspired by \cite[Theorem~1.3]{dolbilin1998multiregular} saying that, for a multi-regular point set $X$, ``the only Delone sets $Y$ all of whose $\rho$-stars are isometric to $\rho$-stars of $X$ are sets globally isometric to $X$''.
After renaming $\rho$-stars as $\al$-clusters, we collected their isometry classes (with weights) into the \emph{isoset} to rephrase \cite[Theorem~1.3]{dolbilin1998multiregular} as a classification of all periodic point sets by isosets. 
The $\al$-equivalence and isoset in Definition~\ref{dfn:isoset} can be refined by labels such as chemical elements, which keeps Theorem~\ref{thm:isoset_complete} valid for labeled points.
\smallskip

When comparing sets from a finite database, it suffices to build their isosets only up to a common upper bound of a stable radius $\al$ in Lemma~\ref{lem:upper_bounds}(c).

\section{Continuous metrics on isometry classes of periodic sets in $\R^n$}
\label{sec:metrics}

This section proves the continuity of the isoset $I(S;\al)$ in Theorem~\ref{thm:continuity} by using the Earth Mover's Distance (EMD) from Definition~\ref{dfn:EMD_isosets}.
For a point $p\in\R^n$ and a radius $\ep$, 
the {\em closed} ball $\bar B(p;\ep)=\{q\in\R^n \vl |\vec q - \vec p|\leq\ep\}$ has as its the boundary $(n-1)$-dimensional sphere $\bd\bar B(p;\ep)\subset\R^{n}$.
The \emph{$\ep$-offset} of any set $C\subset\R^n$ is the Minkowski sum $C+\bar B(0;\ep)=\{\vec p+\vec q \vl p\in C, q\in \bar B(0;\ep)\}$.
\medskip

Then the directed \emph{Hausdorff} distance from Definition~\ref{dfn:Hausdorff+bottleneck}(a) $d_{\vec H}(C,D)$ 
is the minimum radius $\ep\geq 0$ such that $C\subseteq D+\bar B(0;\ep)$.
Definition~\ref{dfn:tolerant_metric} introduces the crucial new metric, which will be explicitly computed in Lemma~\ref{lem:max-min_formula}. 

\begin{dfn}[boundary tolerant metric $\BT$ on isometry classes of clusters] 
\label{dfn:tolerant_metric}
For a radius $\al$ and periodic point sets $S,Q\subset\R^n$, let clusters $C(S,p;\al),C(Q,q;\al)$ represent isometry classes $\si\in I(S;\al),\xi\in I(Q;\al)$, respectively.
The \emph{boundary tolerant} metric $\BT(\si,\xi)$ 
 is defined as the minimum $\ep\geq 0$ such that 
\smallskip

\noindent
(\ref{dfn:tolerant_metric}a) 
$C(Q,q;\al-\ep)\subseteq f(C(S,p;\al))+\bar B(0;\ep)$ for some $f\in\Or(\R^n;p,q)$, and 
\smallskip

\noindent
(\ref{dfn:tolerant_metric}b) 
$C(S,p;\al-\ep)\subseteq g(C(Q,q;\al))+\bar B(0;\ep)$ for some $g\in\Or(\R^n;q,p)$.
\bs
\end{dfn}

In Definition~{\ref{dfn:tolerant_metric}}, if one cluster consists of only its centre, e.g. $C(S,p;\al)=\{p\}$, then the boundary tolerant metric is $\BT=\max\{|\vec s-\vec q| \mid s\in C(Q,q;\al)\}$.

\begin{lem}[correctness of $\BT$]
\label{lem:tolerant_metric}
The metric $\BT(\si,\xi)$ in Definition~\ref{dfn:tolerant_metric} 
is independent of cluster representatives and
satisfies the metric axioms below: 
\medskip

\noindent
(\ref{lem:tolerant_metric}a)
$\BT(\si,\xi)=0$ if and only if $\si=\xi$ as isometry classes
 of $\al$-clusters;
\medskip

\noindent
(\ref{lem:tolerant_metric}b)
symmetry : $\BT(\si,\xi)=\BT(\xi,\si)$ for any isometry classes of $\al$-clusters;
\medskip

\noindent
(\ref{lem:tolerant_metric}c)
triangle inequality : $\BT(\si,\zeta)\leq \BT(\si,\xi)+\BT(\xi,\zeta)$ for any $\si,\xi,\zeta$.
\bs
\end{lem}

\begin{exa}[square lattice vs hexagonal]
\label{exa:square_vs_hexagon}
The isoset $I(\La;\al)$ of any lattice $\La\subset\R^n$ containing the origin $0$ consists of a single isometry class $[C(\La,0;\al)]$, see Example~\ref{exa:isosets_lattices}.
For the square (hexagonal) lattice with minimum inter-point distance 1 in Fig.~\ref{fig:square_vs_hexagon}, the cluster $C(\La,0;\al)$ consists of only 0 for $\al<1$ and includes four (six) nearest neighbors of 0 for $\al\geq 1$.
Hence $\sym(\La,0;\al)$ stabilizes as the symmetry group of the square (regular hexagon) for $\al\geq 1$.
The lattices have the minimum stable radius $\al(\La)=2$ and $\be(\La)=1$ by Example~\ref{exa:upper_bounds}(c).

\begin{figure}[h!]
\includegraphics[width=\linewidth]{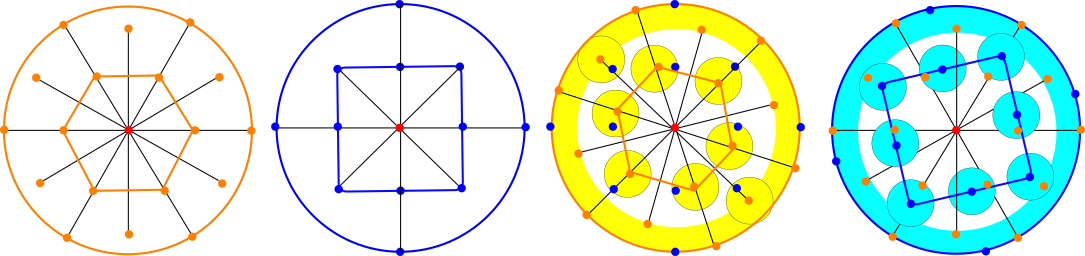}
\caption{
Example~\ref{exa:square_vs_hexagon} computes the metric $\BT$ from Definition~\ref{dfn:tolerant_metric} for the isometry classes of the $2$-clusters
in the square and hexagonal lattices $\La_4,\La_6$.
\textbf{1st}: the 2-cluster $C(\La_6,0;2)$ with its boundary circle $\bd\bar B(0;2)$;
\textbf{2nd}: the 2-cluster $C(\La_4,0;2)$ with its boundary circle $\bd\bar B(0;2)$;
\textbf{3rd}: for $\ep=\sqrt{2}-1\approx 0.41$, the cluster $C(\La_4,0;2)$ is covered by the yellow $\ep$-offset of $C(\La_6,0;2)\cup\bd\bar B(0;2)$ rotated through $15^\circ$ clockwise.
\textbf{4th}: $C(\La_6,0;2)$ is covered by the blue $\ep$-offset of $C(\La_4,0;2)\cup\bd\bar B(0;2)$ rotated through $15^{\circ}$ anticlockwise, so $\BT=\sqrt{2}-1$.}
\label{fig:square_vs_hexagon}
\end{figure}

\noindent
Fig.~\ref{fig:square_vs_hexagon} illustrates the computations whose extra details are in Example~\ref{exa:square_vs_hexagon_details}.
\bs
\end{exa}

Non-isometric periodic sets $S,Q$ such as perturbations in Fig.~\ref{fig:square_lattice_perturbations} can have isosets of different numbers of isometry classes.
A distance between these weighted distributions of different sizes can be measured by EMD below.

\begin{dfn}[Earth Mover's Distance on isosets]
\label{dfn:EMD_isosets}
Let periodic point sets $S,Q\subset\R^n$ have a common stable radius $\al$ and isosets $I(S;\al)=\{(\si_i,w_i)\}$ and $I(Q;\al)=\{(\xi_j,v_j)\}$, where $i=1,\dots,m(S)$ and $j=1,\dots,m(Q)$.
The \emph{Earth Mover's Distance} \cite{rubner2000earth} is 
$\EMD( I(S;\al), I(Q;\al) )=\sum\limits_{i=1}^{m(S)} \sum\limits_{j=1}^{m(Q)} f_{ij} \BT(\si_i,\xi_j)$ minimized over \emph{flows} 
$f_{ij}\in[0,1]$ subject to 
$\sum\limits_{j=1}^{m(Q)} f_{ij}\leq w_i$ for $i=1,\dots,m(S)$,
$\sum\limits_{i=1}^{m(S)} f_{ij}\leq v_j$ for $j=1,\dots,m(Q)$, and
$\sum\limits_{i=1}^{m(S)}\sum\limits_{j=1}^{m(Q)} f_{ij}=1$.
\bs
\end{dfn}

\begin{lem}[EMD is a metric on isosets]
\label{lem:EMD_metric}
The Earth Mover's Distance from Definition~\ref{dfn:EMD_isosets} satisfies the metric axioms for all $\al$ and periodic sets $S,Q,T$.
\medskip

\noindent
(\ref{lem:EMD_metric}a)
$\EMD( I(S;\al), I(Q;\al) )=0$ if and only if $I(S;\al)=I(Q;\al)$;
\medskip

\noindent
(\ref{lem:EMD_metric}b)
$\EMD( I(S;\al), I(Q;\al) )=\EMD( I(Q;\al), I(S;\al) )$;
\medskip

\noindent
(\ref{lem:EMD_metric}b)
$\EMD( I(S;\al), I(Q;\al) )+\EMD( I(Q;\al), I(T;\al) )\geq \EMD( I(S;\al), I(T;\al) )$.
\end{lem}

\begin{exa}[EMD for lattices with $d_B=+\infty$]
\label{exa:EMD_cluster}
\cite[Example~2.1]{widdowson2022resolving} 
showed that the lattices $S=\Z$ and $Q=(1+\de)\Z$ have the bottleneck distance $d_B(S,Q)=+\infty$ for any $\de>0$.
We show that $S,Q$ have Earth Mover's Distance $\EMD=2\de$ at their common stable radius $\al=2+2\de$. 
The bridge lengths are $\be(S)=1$ and $\be(Q)=1+\de$.
The $\al$-cluster $C(S,0;\al)$ contains non-zero points for $\al\geq 1$, e.g. $C(S,0;1)=\{0,\pm 1\}$.
The symmetry group $\sym(S,0;\al)=\Z_2$ includes a non-trivial reflection with respect to 0 for all $\al\geq 1$, so the stable radius of $S$ is any $\al\geq \be+1=2$.
Similarly, $Q$ has $\be(Q)=1+\de$ and stable radii $\al\geq 2(1+\de)$.
The Earth Mover's Distance between $I(S;\al)$ and $I(Q;\al)$ at the common stable radius $\al=2+2\de$ equals the metric $\BT$ between the only $\al$-clusters $C(S,0;\al)=\{0,\pm 1,\pm 2\}$ and $C(Q,0;\al)=\{0,\pm (1+\de),\pm2(1+\de)\}$. 
\medskip

By Definition~\ref{dfn:tolerant_metric} we look for a minimum $\ep>0$ such that the cluster $C(S,0;\al-\ep)$ is covered by $\ep$-offsets of $\pm (1+\de),\pm2(1+\de)$ and vice versa.
If $\ep<2\de$, the points $\pm 2\in C(S,0;\al-\ep)$ cannot be $\ep$-close to $\pm (1+\de),\pm1(+\de)$, but $\ep=2\de$ is large enough.
The cluster $C(Q,0;\al-2\de)=\{0,\pm (1+\de)\}$ is covered by the $2\de$-offset of $C(S,0;\al)=\{0,\pm 1,\pm 2\}$, so $\EMD(I(S;\al),I(Q;\al))
=2\de$.
\bs 
\end{exa}

\begin{dfn}[\emph{packing radius}]
\label{dfn:packing_radius}
For a discrete set $Q\subset\R^n$, the \emph{packing radius} $r(Q)$
is the minimum half-distance between any points of $Q$.
Also, $r(Q)$ is the maximum radius $r$ such that the open balls $B(p;r)$ are disjoint for all $p\in Q$. 
\bs
\end{dfn}

Lemma~\ref{lem:common_lattice} is proved in the appendix and is needed for Theorem~\ref{thm:continuity}.

\begin{lem}
\label{lem:common_lattice}
Let periodic point sets $S,Q\subset\R^n$ have bottleneck distance $d_B(S,Q)<r(Q)$, where $r(Q)$ is the packing radius.
Then $S,Q$ have a common lattice $\La$ with a unit cell $U$ such that $S=\La+(U\cap S)$ and $Q=\La+(U\cap Q)$.
\bs
\end{lem}

For rigid motion instead of general isometry, Definition~\ref{dfn:tolerant_metric} of a boundary tolerant metric $\BT$ is updated to $\BT^o$ by considering only orientation-preserving isometries from $\SO(\R^n;p,q)$, which also makes the continuity below valid for oriented isosets $I^o(S;\al)$ under $\EMD$ using $\BT^o$ instead of $\BT$ in Definition~\ref{dfn:EMD_isosets}. 

\begin{thm}[continuity of isosets under perturbations]
\label{thm:continuity}
Let periodic point sets $S,Q\subset\R^n$ have a bottleneck distance $d_B(S,Q)<r(Q)$, where $r(Q)$ is the packing radius in Definition~\ref{dfn:packing_radius}. 
Then the isosets $I(S;\al),I(Q;\al)$ are close in the Earth Mover's Distance: 
$\EMD( I(S;\al), I(Q;\al) )\leq 2d_B(S,Q)$ for $\al\geq 0$. 
\bs
\end{thm}
\begin{proof}
By Lemma~\ref{lem:common_lattice} the given periodic point sets $S,Q$ have a common unit cell $U$.
Let $g:S\to Q$ be a bijection such that $|\vec p-g(\vec p)|\leq\ep=d_B(S,Q)=\inf\limits_{g:S\to Q} \sup\limits_{p\in S}|\vec p-g(\vec p)|$ for all points $p\in S$.
Since the bottleneck distance 
$\ep<r(Q)$ is small, the bijective image $g(p)$ of any point $p\in S$ is a unique $\ep$-close point of $Q$ and vice versa.
Hence we can assume that the common unit cell $U\subset\R^n$ contains the same number (say, $m$) points from $S$ and $Q$.
Expand the initial $m(S)$ isometry classes $(\si_i,w_i)\in I(S;\al)$
to $m$ isometry classes (with equal weights $\frac{1}{m}$) represented by clusters $C(S,p;\al)$ for $m$ points $p\in S\cap U$.
If the $i$-th initial isometry class had a weight $w_i=\frac{k_i}{m}$, $i=1,\dots,m(S)$, the expanded isoset contains $k_i$ equal isometry classes of weight $\frac{1}{m}$.
For example, the 1-regular set $S_1$ in Fig.~\ref{fig:alpha-clusters} has the isoset consisting of a single class $[C(S_1,p;\al)]$, which is expanded to four identical classes of weight $\frac{1}{4}$ for the four points in the motif.
The isoset $I(Q;\al)$ is similarly expanded to $m$ isometry classes of weight $\frac{1}{m}$. 
\medskip

For any point $p\in S\cap U$, the image $g(p)\in Q$ has a unique point $h(p)\in Q\cap U$ such that $h(p)$ is equivalent to $g(p)$ modulo the lattice of $Q$.
Then the $\al$-clusters of $g(p)$ and $h(p)$ in $Q$ are isometric for any $\al\geq 0$. 
The bijection $p\mapsto h(p)$ between the expanded motifs of $S,Q$ induces the bijection between the expanded sets of $m$ isometry classes.
Each correspondence $\si_l\mapsto\xi_l$ in the latter bijection can be visualized as the flow $f_{ll}=\frac{1}{m}$ for $l=1,\dots,m$, so $\sum\limits_{l=1}^m f_{ll}=1$.
\medskip

To show that the Earth Mover's Distance (EMD) between any initial isoset and its expansion is 0, we collapse all identical isometry classes in the expanded isosets, but keep the arrows with the flows above.
Only if both tail and head of two (or more) arrows are identical, we collapse these arrows into one arrow that gets the total weight.
All equal weights $\frac{1}{m}$ correctly add up at heads and tails of final arrows to the initial weights $w_i,v_j$ of isometry classes.
So the total sum of flows is $\sum\limits_{i=1}^{m} \sum\limits_{j=1}^{m} f_{ij}=1$ as required by Definition~\ref{dfn:EMD_isosets}.
It suffices to consider below the EMD only for the expanded isosets of exactly $m$ classes.
\medskip

We will estimate the boundary tolerant metric between isometry classes $\si_l,\xi_l$ whose centers $p$ and $g(p)$ are $\ep$-close within the common unit cell $U$.
For any fixed point $p\in S\cap U$, shift $S$ by the vector $g(\vec p)-\vec p$.
This shift makes $p\in S$ and $g(p)\in Q$ identical and keeps all pairs $q,g(q)$ for $q\in C(S,p;\al)$ within $2\ep$ of each other.
Using the identity map $f$ 
in Definition~\ref{dfn:tolerant_metric}, we get the upper bound $\BT([C(S,p;\al)],[C(Q,g(p);\al)])\leq 2\ep$.
Then $\EMD(I(S;\al),I(Q;\al))\leq
\sum\limits_{l=1}^{m} f_{ll} \BT([C(S,p;\al)],[C(Q,g(p);\al)])\leq
2\ep\sum\limits_{l=1}^{m} f_{ll}=
2\ep$ as required.
\end{proof}

Corollary~\ref{cor:metric}a justifies that the EMD satisfies all metric axioms for periodic point sets that have a stable radius $\al$.
Corollary~\ref{cor:metric}b avoids this dependence on $\al$ and scales any periodic point set $S$ to the minimum stable radius $\al(S)=1$.

\begin{cor}
\label{cor:metric}
\textbf{(a)}
For $\al>0$, $\EMD( I(S;\al), I(Q;\al) )$ is a metric on the space of isometry classes of all periodic point sets with a stable radius $\al$ in $\R^n$.
\medskip

\noindent
\textbf{(b)}
For a periodic point set $S\subset\R^n$, let $S/r(S)\subset\R^n$ denote $S$ after uniformly dividing all vectors by the packing radius $r(S)$.
Then  $|r(S)-r(Q)|+\EMD( I(S/r(S);1), I(Q/r(Q);1) )$ is a metric on all periodic point sets. 
\bs
\end{cor}
\begin{proof}
\textbf{(a)}
Lemma~\ref{lem:EMD_metric} proved the metric axioms for the $\EMD$ on isosets.
The equality $I(S;\al)=I(Q;\al)$ is equivalent to isometry $S\simeq Q$ by Theorem~\ref{thm:isoset_complete}.
\medskip

\noindent
\textbf{(b)}
By part (a), $\EMD( I(S/r(S);1), I(Q/r(Q);1) )$ satisfies the symmetry and triangle inequality, which are preserved by adding the Euclidean distance $d=|r(S)-r(Q)|$ between the packing radii.
The equality $\EMD=0$ means that $S/r(S)\simeq Q/r(Q)$ are isometric.
Hence $S,Q$ are isometric up to a uniform factor.
Adding the distance $d=|r(S)-r(Q)|$ guarantees that the sum becomes zero only if $r(S)=r(Q)$, so the given sets $S,Q$ should be truly isometric.
\end{proof}

The metric $\EMD( I(S;\al), I(Q;\al) )$ is measured in the same units as atomic coordinates, say in angstroms: 1$\angstrom=10^{-10}$m, and hence is physically meaningful.
By Theorem~\ref{thm:continuity}, a small value $\de=\EMD( I(S;\al), I(Q;\al) )$ means that atoms of $S$ should be perturbed by at least $0.5\de$ on average for a complete match with $Q$.   
Since crystals are practically compared within a finite dataset, we can take any common upper bound of $\al(S)$ from Lemma~\ref{lem:upper_bounds}, also in Corollary~\ref{cor:metric}(b).

\section{Algorithms to test isometry and to approximate metrics on isosets}
\label{sec:algorithms}

This section describes time complexities for computing the complete invariant isoset (Theorem~\ref{thm:compute_isoset}), comparing isosets (Corollary~\ref{cor:compare_isosets}), approximating 
the boundary tolerant metric $\BT$ and 
Earth Mover's Distance on isosets (Corollary~\ref{cor:approximate_EMD}).  
All estimates will use the geometric complexity $\GC(S)$ below.

\begin{dfn}[geometric complexity $\GC$]
\label{dfn:geom_complexity}
Let a periodic point set $S\subset\R^n$ have an asymmetric unit of $m$ points in a  cell $U$ of volume $\vol[U]$.
Let $L$ be the symmetry characteristic for $\al_0=2R(S)$ in Lemma~\ref{lem:upper_bounds}(c), where $R(S)$ is the covering radius.
The \emph{geometric complexity} is 
$\GC(S)=\frac{(10(L+m+2)R(S)/n)^n}{2\vol[U]}$.
\bs
\end{dfn}

Let $V_n=\frac{\pi^{n/2}}{\Ga(\frac{n}{2}+1)}$ be the volume of the unit ball in $\R^n$, where the Gamma function $\Ga$ has $\Ga(k)=(k-1)!$ and $\Ga(\frac{k}{2}+1)=\sqrt{\pi}(k-\frac{1}{2})(k-\frac{3}{2})\cdots\frac{1}{2}$ for any integer $k\geq 1$.
Set $\nu(U,\al,n)=\frac{(\al+d)^n V_n}{\vol[U]}$, where 
$d=\sup\limits_{p,q\in U}|\vec p-\vec q|$ is a longest diagonal of a unit cell $U$.
All complexities assume the real Random-Access Machine (RAM) model and a fixed dimension $n$ of Euclidean space $\R^n$.
\medskip

The main input size of a periodic set is the number $m$ of motif points because the length of a standard Crystallographic Information File is linear in $m$.
For a fixed dimension $n$, the big $O$ notation $O(m^n)$ in all complexities means a function $t(m)$ such that $t(m)\leq C m^n$ for a fixed constant $C$ independent of $m$.
We will include all other parameters depending on a periodic point set $S$.

\begin{lem}[a local cluster]
\label{lem:compute_cluster}
Let a periodic point set $S\subset\R^n$ have $m$ points in a unit cell $U$.
For any stable radius $\al\geq 0$ and $p\in M=S\cap U$, the cluster $C(S,p;\al)$ has at most $k=\nu m$ points and can be found in time $\nu O(m)$, where $\nu\leq \GC(S)$ for geometric complexity $\GC(S)$ from Definition~\ref{dfn:geom_complexity}.
\bs
\end{lem}

\begin{thm}[computing an isoset]
\label{thm:compute_isoset}
For any periodic point set $S\subset\R^n$ given by a motif $M$ of $m$ points in a unit cell $U$, the isoset $I(S;\al)$ at a stable radius $\al$ can be found in time $O(m^2 k^{\lceil n/3\rceil}\log k)$, where $k=\nu m$ 
for $\nu\leq \GC(S)$.
\bs
\end{thm}
\begin{proof}
Lemma~\ref{lem:compute_cluster} computes the $\al$-clusters of $m$ points $p\in M$ in time $O(k)$.
To verify a congruence (isometry) of finite sets $A,B\subset\R^n$, the algorithm from \cite{brass2000testing} first moves the centers of mass of $A,B$ to $0\in\R^n$.
We instead move the centers of given clusters $A,B$ to the origin and then follow \cite{brass2000testing} to check if the shifted clusters are related by an isometry $f\in\Or(\R^n;0)$ in time $O(k^{\lceil n/3\rceil}\log k)$.
The isoset $I(S;\al)$ is obtained after identifying isometric clusters for $m$ points through $O(m^2)$ pairwise comparisons.
The total time is $O(m^2 k^{\lceil n/3\rceil}\log k)$.
\end{proof}


\begin{cor}[comparing isosets]
\label{cor:compare_isosets}
There is an algorithm to check if any periodic point sets $S,Q\subset\R^n$ with motifs of at most $m$ points are isometric in total time $O(m^2 k^{\lceil n/3\rceil}\log k)$, where $k=\nu m$ for $\nu\leq\max\{\GC(S),\GC(Q)\}$.
\bs
\end{cor}
\begin{proof}
Theorem~\ref{thm:compute_isoset} finds $I(S;\al),I(Q;\al)$ with a common stable radius in time $O(m^2 k^{\lceil n/3 \rceil}\log k)$, where each cluster has $k=\nu m$ points by Lemma~\ref{lem:compute_cluster}.
Any classes from $I(S;\al),I(Q;\al)$ are compared \cite{brass2000testing} in time $O(k^{\lceil n/3 \rceil}\log k)$.
Then $O(m^2)$ comparisons suffice to check if there is a bijection $I(S;\al)\lra I(Q;\al)$.  
\end{proof}

\begin{dfn}[directed distances $d_{\vec R}$ and $d_{\vec M}$]
\label{dfn:directed_distances}
\textbf{(a)}
For any sets $C,D\subset\R^n$, the directed \emph{rotationally invariant} distance  $d_{\vec R}(C,D)=\min\limits_{f\in\Or(\R^n)}d_{\vec H}(C,f(D))$ is minimized over all maps $f\in\Or(\R^n;0)$, which fix the origin $0\in\R^n$. 
\smallskip

\noindent
\textbf{(b)}
For any finite sets $C,D\subset\R^n$, order all points $p_1\dots,p_k\in C$ by increasing distance to the origin $0$.
The \emph{radius} of $C$ is $R(C)= \max\limits_{p\in C}|p|$.
Define the directed \emph{max-min distance} $d_{\vec M}(C,D)=\max\limits_{i=1,\dots,k}\min\{\;\al-|p_i|,\; d_{\vec R}(\{p_1,\dots,p_i\}, D) \;\}$.
\bs
\end{dfn}

If $C'\subset C$, then $d_{\vec R}(C',D)\leq d_{\vec R}(C,D)$.
Let $C,D\subset\bar B (0; \al)$ be finite sets including the origin $0$.
If $C=\{0\}$, then $d_{\vec R}(C,D)=0$  because $C\subset D$, but $d_{\vec R}(D,C)=R(D)$ is the radius of $D$ because $D\subset \{0\}+\bar B(0;\ep)$ only for $\ep\geq R(D)$.
Definition~\ref{dfn:directed_distances}, Lemma~\ref{lem:max-min_formula} and hence all further results work for rigid motion by restricting all maps to the special orthogonal group $\SO(\R^n;0)$.

\begin{lem}[max-min formula for $d_{\vec R}$ via $d_{\vec M}$]
\label{lem:max-min_formula}
For any finite sets $C,D\subset\R^n$, $\al\geq R(C)$, the distance $d_{\vec R}(C\cup\bd\bar B(0;\al),D\cup\bd\bar B(0;\al))$ equals $d_{\vec M}(C,D)$.
\bs
\end{lem}

\begin{exa}[max-min formula]
\label{exa:max-min_formula}
Consider the subcluster $C\subset C(\La_4,0;2)$ of the points $p_1=(1,0)$, $p_2=(1,1)$, $p_3=(1,-1)$, $p_4=(2,0)$ from the square lattice $\La_4$ in Fig.~{\ref{fig:square_vs_hexagon}}.
Let $\al=2$ and $D=C(\La_6,0;2)$ be the 2-cluster of the hexagonal lattice $\La_6$.
Then $d_{\vec R}(p_1,D)=0$ because $p_1$ coincides with $(1,0)\in D$.
Then $d_{\vec R}(\{p_1,p_2\},D)=\sqrt{2}-1$, because the cloud $D$ after the clockwise rotation through $15^\circ$ has the points $(\cos 15^\circ,-\sin 15^\circ)$ and $(\frac{1}{\sqrt{2}},\frac{1}{\sqrt{2}})$ at distances $\sqrt{(\cos 15^\circ -1)^2+\sin^2 15^\circ}\approx 0.26$, $\sqrt{2}-1\approx 0.41$ to $p_1,p_2$, respectively.
Then $d_{\vec R}(\{p_1,p_2,p_3\},D)=\sqrt{2}-1$ because the same rotated image of $D$ has $(\sqrt{\frac{3}{2}}, -\sqrt{\frac{3}{2}})$ at the distance $\sqrt{3}-\sqrt{2}\approx 0.32$ to $p_3$. 
For $i=1$, $\min\{\al-|p_1|,d_{\vec R}(p_1,D)\}=\min\{2-1,0\}=0$.
For $i=2,3$, $\min\{\al-|p_2|,d_{\vec R}(\{p_1,p_2\},D)\}=\min\{\al-|p_3|,d_{\vec R}(\{p_1,p_2,p_3\},D)\}=\min\{2-\sqrt{2},\sqrt{2}-1\}=\sqrt{2}-1$.
For $i=4$, $\min\{\al-|p_4|,d_{\vec R}(C,D)\}=0$ since $\al=2=|p_4|$.
The maximum value is $\sqrt{2}-1$, so Example~\ref{exa:square_vs_hexagon} fits Lemma~\ref{lem:max-min_formula}.
\bs
\end{exa}

Lemma~\ref{lem:rot-inv_distance} extends \cite[section 2.3]{goodrich1999approximate} from  $n=3$ to any dimension $n>1$. 

\begin{lem}[approximating $d_{\vec R}$]
\label{lem:rot-inv_distance}
Let a cloud $C\subset\R^n$ consist of $k=|C|$ points ordered by distances $|p_1|\leq\dots\leq|p_k|$ from the origin 
and $\ar{C}$ denote the number of different vectors $\vec p/|\vec p|$ for $p\in C$. 
For each $j=1,\dots,k$, consider the subcloud $C_j=\{p_1,\dots,p_{j}\}$.
For any cloud $D\subset\R^n$ of $|D|$ points, all distances $d_j=d_{\vec R}(C_j,D)$ from Definition~{\ref{dfn:directed_distances}} for $j=1,\dots,k$ can be approximated by some $d'_j$ in time $O(|C|\ar{C}^{n-1}|D|)$ so that $d_j\leq d'_j\leq \om d_j$, $\omega=1+\frac{1}{2}n(n-1)$.
\bs
\end{lem}

The proof of Lemma~\ref{lem:rot-inv_distance} uses only orientation-preserving isometries from $\SO(\R^n,0)$.
Hence the upper bounds from Lemma~\ref{lem:rot-inv_distance}, Theorem~\ref{thm:approximate_BT}, and Corollary~\ref{cor:approximate_EMD} work for both cases of rigid motion and general isometry in $\R^n$.

\begin{thm}[approximating $\BT$]
\label{thm:approximate_BT}
Let periodic point sets $S,Q\subset\R^n$ have isometry classes $\si,\xi$ represented by clusters $C,D$ of a radius $\al$, respectively.
In the notations of Lemma~{\ref{lem:rot-inv_distance}}, 
$\BT(\si,\xi)$ from Definition~{\ref{dfn:tolerant_metric}} can be approximated with the factor $\omega=1+\frac{1}{2}n(n-1)$ in time $O(|C|(\ar{C}^{n-1}+\ar{D}^{n-1})|D|)$.  
\bs
\end{thm}
\begin{proof}
By Definitions~{\ref{dfn:tolerant_metric}} and {\ref{dfn:directed_distances}}, 
the boundary tolerant metric $\BT(\si,\xi)$ is the maximum of 
$d_{\vec R}(C\cup\bd\bar B(0;\al),D\cup\bd\bar B(0;\al))$ and $d_{\vec R}(D\cup\bd\bar B(0;\al),C\cup\bd\bar B(0;\al))$.
Lemma~{\ref{lem:max-min_formula}} implies that $\BT(\si,\xi)=\max\{d_{\vec M}(C,D),d_{\vec M}(D,C)\}$.
It remains to compute required approximations of the two distances $d_{\vec M}$ above.
\medskip
 
Let $C$ consist of $k$ points ordered by distances $|p_1|\leq\dots\leq|p_k|$ from the origin. 
For each $j=1,\dots,k$, consider the subcloud $C_j=\{p_1,\dots,p_{j}\}$.
We use the approximation $d'_j$ from Lemma~{\ref{lem:rot-inv_distance}} 
to compute $d'=\max\limits_{i=1,\dots,k}\min\{\al-|p_i|,\; d'_i \}$ in the extra time $O(|C|)$.
Now we check that this final approximation $d'$ is between $d_{\vec M}(C,D)=\max\limits_{i=1,\dots,k}\min\{\al-|p_i|,\; d_i \}$ in Lemma~{\ref{lem:max-min_formula}} and $\omega d_{\vec M}(C,D)$.
\medskip

The inequalities $d_j\leq d'_j \leq \omega d_j$ for $j=1,\dots,k$  from Lemma~{\ref{lem:rot-inv_distance}} imply that
$\min\{\al-|p_j|, d_j\}\leq \min\{\al-|p_j|, d'_j\} \leq \min\{\al-|p_j|,\omega d_j \}$.
By fixing an index $j$ maximizing the left-hand side minimum above, we conclude that
$d_{\vec M}(C,D)=$ 
$\min\{\al-|p_j|,\; d_j\}\leq \max\limits_{i=1,\dots,k}\min\{\al-|p_i|,\; d'_i \}=d'$.
By fixing an index $j$ maximizing the middle side minimum above, we get the following upper bound: $d'=\min\{\al-|p_j|,\; d'_j\}\leq \max\limits_{i=1,\dots,k}\min\{\al-|p_i|,\; \om d_i \}\leq \om d_{\vec M}(C,D)$.
The extra \\ time $O(|C|+|D|)$ for the approximations of $d_{\vec M}(C,D)$ and $d_{\vec M}(D,C)$ is dominated by the total time $O(|C|(\ar{C}^{n-1}+\ar{D}^{n-1})|D|)$ from Lemma~{\ref{lem:rot-inv_distance}}.
\end{proof}

\begin{cor}[approximating EMD]
\label{cor:approximate_EMD}
Let $S,Q\subset\R^n$ be periodic point sets whose motifs have at most $m$ points $p$ and $\chi$ different vectors $\vec p/|\vec p|$.
For any $\al>0$, the metric $\EMD(I(S;\al),I(Q;\al))$ can be approximated with the factor $\omega=1+\frac{1}{2}n(n-1)$ in time $O(\nu^2 m^4 \chi^{n-1})$, where $\nu\leq\max\{\GC(S),\GC(Q)\}$.
\bs
\end{cor}
\begin{proof}
Since $S,Q$ have at most $m$ points in their motifs, their isosets at any radius $\al$ have at most $m$ isometry classes.
By Theorem~\ref{thm:approximate_BT} the distance $\BT(\si,\xi)$ between any classes $\si\in I(S;\al)$ and $\xi\in I(Q;\al)$ of $\al$-clusters up to $k$ points can be approximated with the factor $\omega$ in time $O(k^2 \chi^{n-1})$.
The maximum number of points is $k=\nu m$, where $\nu\leq\max\{\GC(S),\GC(Q)\}$ by Lemma~{\ref{lem:compute_cluster}}.
Since Definition~{\ref{dfn:EMD_isosets}} uses normalized distributions, $\omega$ emerges as a multiplicative upper bound in $\EMD(I(S;\al),I(Q;\al))$.
After computing $O(m^2)$ pairwise distances between sets of $m$ clusters, the exact EMD can be found in the extra time $O(m^3\log m)$ \cite{orlin1993faster}, which is dominated by the time $O(m^2\nu^2 m^2 \chi^{n-1})$ for all cluster distances.
So the total time becomes $O(\nu^2 m^4 \chi^{n-1})$.
The EMD can be approximated \cite[section~3]{shirdhonkar2008approximate} with a constant factor in time $O(m)$. 
\end{proof}

Counting directions $\vec p/|p|$ as points ($\chi\leq m$), for dimension $n=3$, the rough bounds for the isoset and its approximate $\EMD'$ in Theorem~\ref{thm:compute_isoset} and Corollary~\ref{cor:approximate_EMD} are $O(m^3\log m)$ and $O(m^6)$, respectively.
Algorithms 1-2 in the appendix describe pseudocodes for
Lemma~\ref{lem:rot-inv_distance}, Theorem~\ref{thm:approximate_BT}, and Corollary~\ref{cor:approximate_EMD}.
 
\section{A lower bound for continuous metrics via simpler invariants}
\label{sec:lower_bound}

Theorem~\ref{thm:lower_bound} gives a lower bound for EMD in terms of the simpler invariant \emph{Pointwise Distance Distribution} \cite{widdowson2022resolving}, 
see Definition~\ref{dfn:PDD} below. 
If $S$ is a lattice or a 1-regular set, then all points are isometrically equivalent, so they  have the same distances to all their neighbors.
In this case, $\PDD(S;k)$ is a single row of $k$ distances, which is the vector $\AMD(S;k)$ of Average Minimum Distances \cite{widdowson2022average}.

\begin{dfn}[Pointwise Distance Distribution $\PDD$]
\label{dfn:PDD}
Let a periodic set $S=\La+M$ have points $p_1,\dots,p_m$ in a unit cell.
For $k\geq 1$, consider the $m\times k$ matrix $D(S;k)$, whose $i$-th row consists of the ordered Euclidean distances $d_{i1}\leq\cdots\leq d_{ik}$ from $p_i$ to its first $k$ nearest neighbors in the full set $S$.
The rows of $D(S;k)$ are \emph{lexicographically} ordered as follows.
A row $(d_{i1},\dots,d_{ik})$ is \emph{smaller} than $(d_{j1},\dots,d_{jk})$ if the first (possibly none) distances coincide: $d_{i1}=d_{j1},\dots,d_{il}=d_{jl}$ for $l\in\{1,\dots,k-1\}$ and the next $(l+1)$-st distances satisfy $d_{i,l+1}<d_{j,l+1}$.
If $w$ rows are identical to each other, these rows are collapsed to one row with the \emph{weight} $w/m$.
Put this weight in the extra first column. 
The final matrix of $k+1$ columns is the \emph{Pointwise Distance Distribution} $\PDD(S;k)$.
The Average Minimum Distance $\AMD(S;k)$ is the vector $(\AMD_1,\dots,\AMD_k)$, where $\AMD_i$ is the weighted average of the $(i+1)$-st column of $\PDD(S;k)$.
\bs
\end{dfn}

\begin{thm}[isometry invariance of $\PDD$]
\label{thm:invariance}
For any finite or periodic set $S\subset\R^n$, 
$\PDD(S;k)$ in Definition~\ref{dfn:PDD} is an isometry invariant of $S$ for $k\geq 1$.
\bs
\end{thm}

Theorem~\ref{thm:invariance} and continuity of $\PDD$ in the metric from Definition~\ref{dfn:EMD_PDD} follows from more general results in \cite{widdowson2022resolving}.
The distance between rows {$\vec R_i(S),\vec R_j(Q)$} of PDD matrices below is measured in the metric $L_{\infty}(\vec p,\vec q)=\max\limits_{i=1,\dots,k}|p_i-q_i|$.

\begin{dfn}[Earth Mover's Distance on Pointwise Distance Distributions] 
\label{dfn:EMD_PDD}
Let finite or periodic sets $S,Q\subset\R^n$ have $\PDD(S;k)$, $\PDD(Q;k)$ with  rows $\vec R_i(S),\vec R_j(Q)$ of weights $w_i(S), w_(Q)$ for $i=1,\dots,m(S)$ and $j=1,\dots,m(Q)$, respectively.
A full \emph{flow} from $\PDD(S;k)$ to $\PDD(Q;k)$ is an $m(S)\times m(Q)$ matrix whose element $f_{ij}\in[0,1]$ is called a partial \emph{flow} from {$\vec R_i(S)$ to $\vec R_j(Q)$}.
The \emph{Earth Mover's Distance} is the minimum value of the \emph{cost} 
{$\EMD(I(S),I(Q))=\sum\limits_{i=1}^{m(S)} \sum\limits_{j=1}^{m(Q)} f_{ij} L_\infty(\vec R_i(S),\vec R_j(Q))$} over \emph{flows} $f_{ij}\in[0,1]$ subject to 
$\sum\limits_{j=1}^{m(Q)} f_{ij}\leq w_i(S)$ for $i=1,\dots,m(S)$,
$\sum\limits_{i=1}^{m(S)} f_{ij}\leq w_j(Q)$ for $j=1,\dots,m(Q)$, 
$\sum\limits_{i=1}^{m(S)}\sum\limits_{j=1}^{m(Q)} f_{ij}=1$.
\bs
\end{dfn}

Lemma~\ref{lem:lower_bound} is a partial case of Theorem~\ref{thm:lower_bound} for 1-regular point sets $S,Q$. 

\begin{lem}[lower bound for the tolerant distance $\BT$]
\label{lem:lower_bound}
Let $S,Q\subset\R^n$ be periodic point sets with a common stable radius $\al$.
Choose any points $p\in S$ and $q\in Q$.
Let the distance between isometry classes of $\al$-clusters $\ep=\BT([C(S,p;\al)],[C(Q,q;\al)])$ be smaller than a minimum half-distance between any point of $S$ and $Q$.
Let $k$ be a minimum number of points in the clusters $C(S,p;\al-\ep)$ and $C(Q,q;\al-\ep)$.
Then the $L_{\infty}$ distance between the rows of the points $p,q$ in $\PDD(S;k),\PDD(Q;k)$, {respectively}, is at most
$\ep$.
\bs
\end{lem}

\begin{thm}[lower bound for $\EMD$]
\label{thm:lower_bound}
Let $S,Q\subset\R^n$ be periodic sets with a common stable radius $\al$.
Let $\ep=\EMD(I(S;\al),I(Q;\al))$ and $k$ be the maximum number of points of $S,Q$ in their $(\al-\ep)$-clusters.
If $\ep$ is less than the half-distance between any points of $S,Q$, then $\EMD(\PDD(S;k),\PDD(Q;k))\leq\ep$.
\bs 
\end{thm}
\begin{proof}
To prove that 
$\EMD(\PDD(S;k),\PDD(Q;k))\leq \EMD(I(S;\al),I(Q;\al)),$ we choose optimal flows $f_{ij}\in[0,1]$, $i=1,\dots,m(S)$ and $j=1,\dots,m(Q)$, that minimize 
$\ep=\EMD(I(S;\al),I(Q;\al))$ in Definition~\ref{dfn:EMD_isosets}.
For any points $p_i\in S$ and $q_j\in Q$, let $\vec R_i(S)$ and $\vec R_j(Q)$ be their rows in $\PDD(S;k)$ and $\PDD(Q;k)$, respectively.
Lemma~\ref{lem:lower_bound} gives
$L_\infty(\vec R_i(S),\vec R_j(Q))\leq \BT([C(S,p_i;\al)],[C(Q,q_j;\al)])$.
These inequalities for all indices $i,j$ and the same flows $f_{ij}$ imply that 
$$\sum\limits_{i=1}^{m(S)} \sum\limits_{j=1}^{m(Q)} f_{ij} L_\infty(\vec R_i(S),\vec R_j(Q))\leq \sum\limits_{i=1}^{m(S)} \sum\limits_{j=1}^{m(Q)} f_{ij} \BT(\si_i,\xi_j)=\ep$$ by the choice of $f_{ij}$ 
The left-hand side of the last inequality can become only smaller when minimizing over $f_{ij}$.
Then $\EMD(\PDD(S;k),\PDD(Q;k))\leq\ep$.
\end{proof}

\section{Real experiments, limitations, significance, and a discussion} 
\label{sec:discussion}

In 1930, future Nobel laureate Linus Pauling noticed the ambiguity of crystal structures obtained by diffraction \cite{pauling1930crystal}.
Such \emph{homometric} crystals with identical diffraction patterns were only manually distinguished until now because even the generically complete $\PDD$s coincide for the Pauling periodic sets $P(\pm u)$ for all $u\in(0,0.25)$, see the real overlaid crystals for $u=0.03$ in Fig.~\ref{fig:Pauling}~(left).

\vspace*{-2mm}
\newcommand{\mh}{28mm}
\begin{figure}[h!]
\includegraphics[height=\mh]{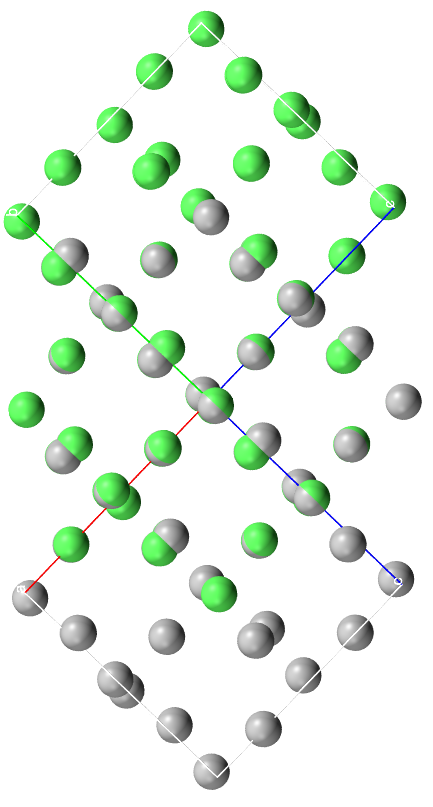}
\hspace*{0mm}
\includegraphics[height=\mh]{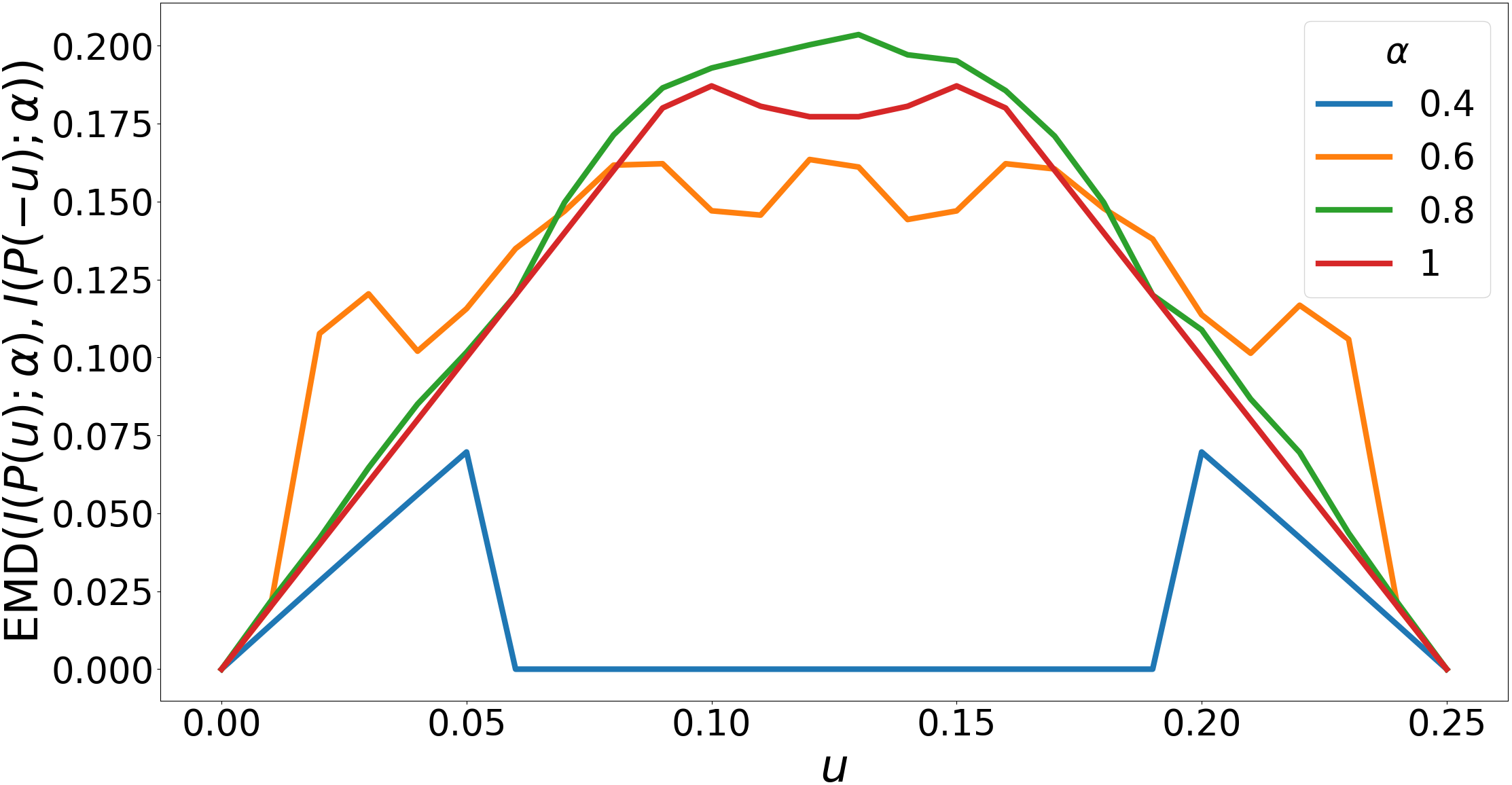}
\hspace*{0mm}
\includegraphics[height=\mh]{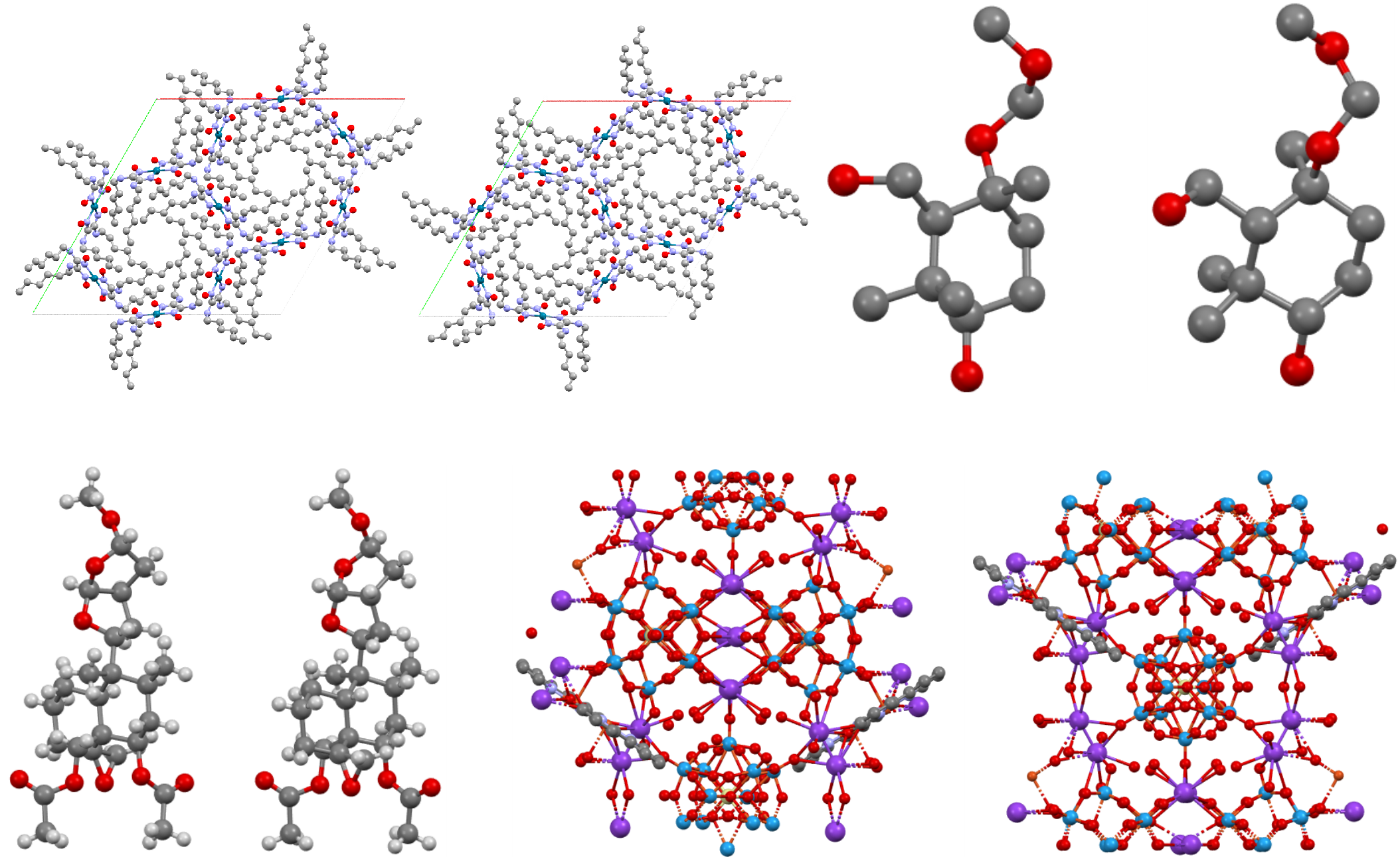}
\caption{\textbf{Left}: a comparison of Pauling's homometric crystals $P(\pm u)$ for $u=0.03$ \cite{pauling1930crystal}, by COMPACK \cite{chisholm2005compack} 
aligning subsets of 48 atoms and outputs RMSD, which fails the triangle inequality. 
The atoms from different $P(\pm 0.03)$ are shown in green and gray.
\textbf{Middle}: the pairs of $P(\pm u)$ have $\EMD'>0$ for all $u\in(0,0.25)$ and $\al>0.4$ (running time 50 ms for $u=0.03$ and $\al=0.5$). 
\textbf{Right}: 
four pairs of mirror images in the CSD are indistinguishable by all past invariants but have approximate $\EMD'>0$ for all radii $\al>1.5\angstrom$ in Fig.~\ref{fig:mirror_images}~(left). 
}
\label{fig:Pauling} 
\end{figure}

\vspace*{-4mm}
\begin{figure}[h!]
\includegraphics[width=0.49\textwidth]{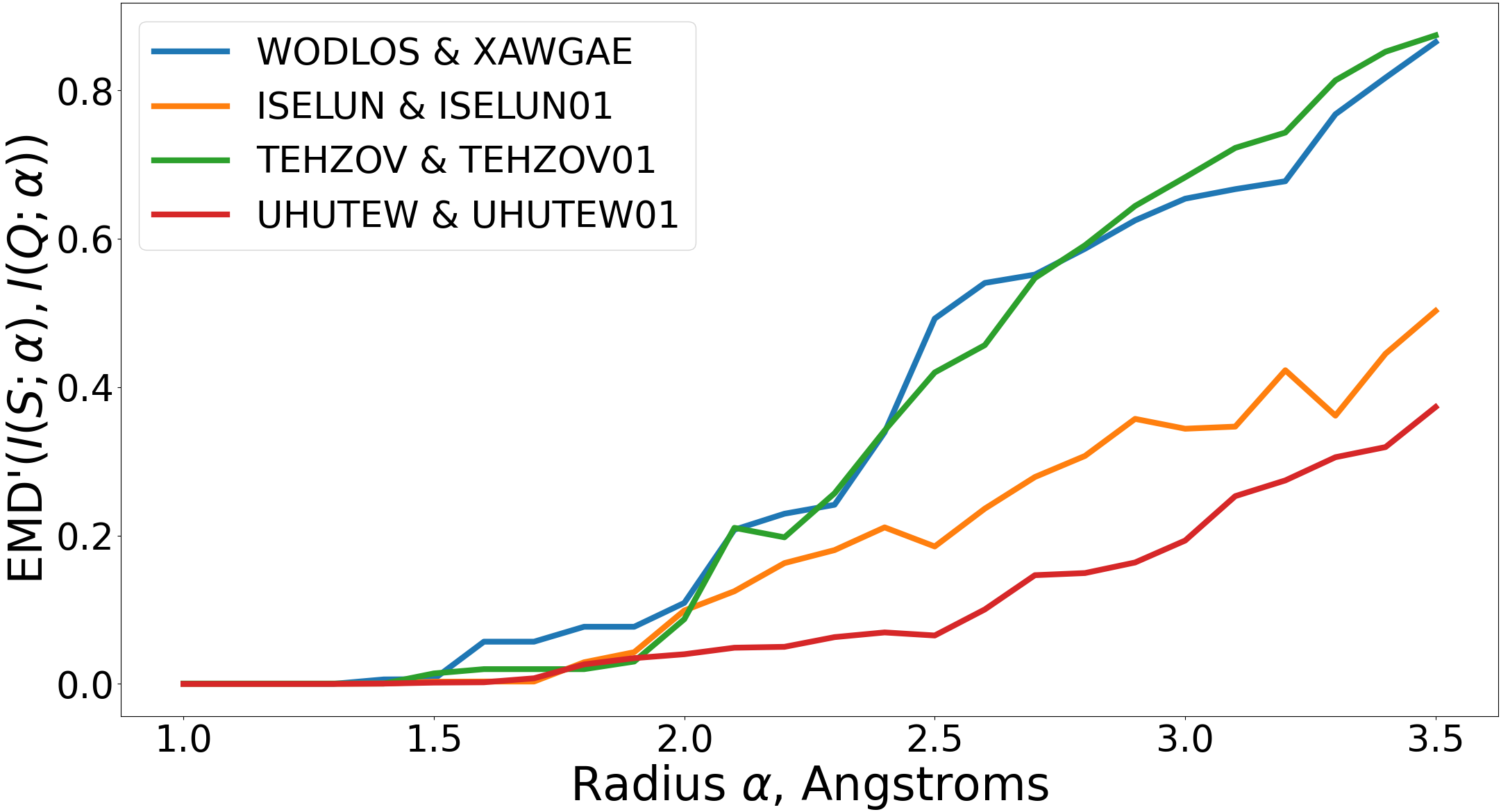}
\includegraphics[width=0.49\textwidth]{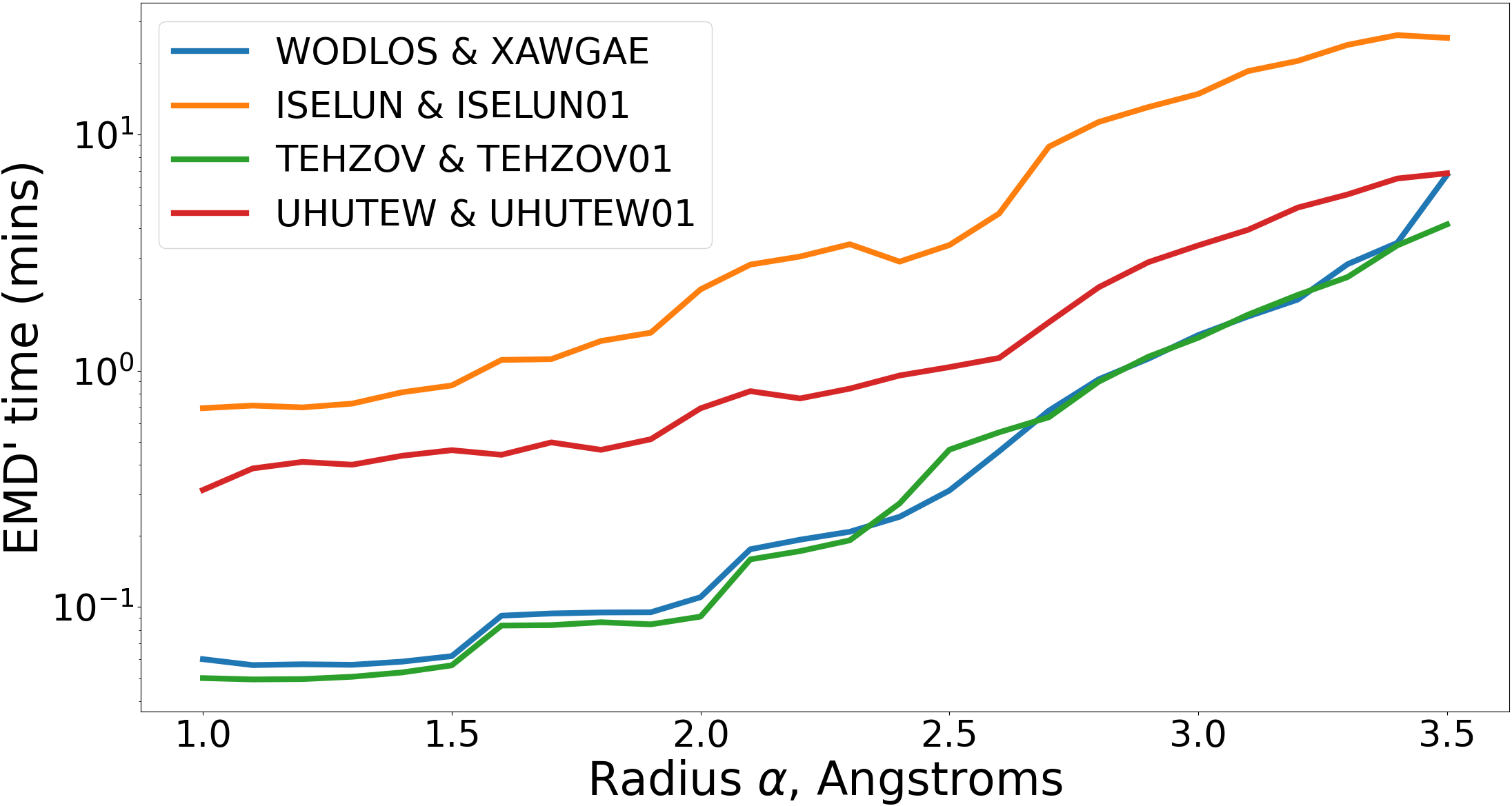}
\caption{
The isosets distinguish all four pairs of mirror images given by their codes in the CSD.
\textbf{Left}: approximate $\EMD'$ for different radii $\al$.
\textbf{Right}: running times on a desktop.
}
\label{fig:mirror_images} 
\end{figure}

\begin{figure}[h!]
\includegraphics[width=0.49\textwidth]{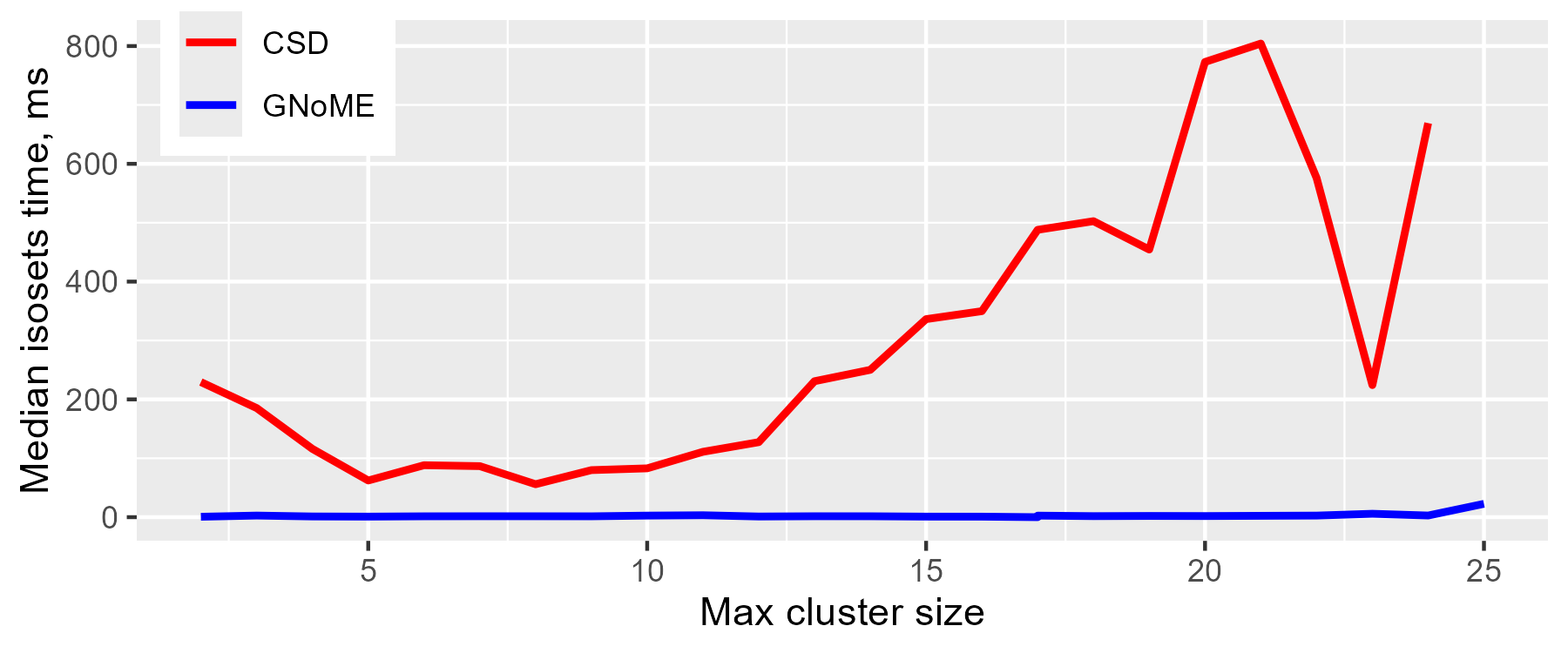}
\includegraphics[width=0.49\textwidth]{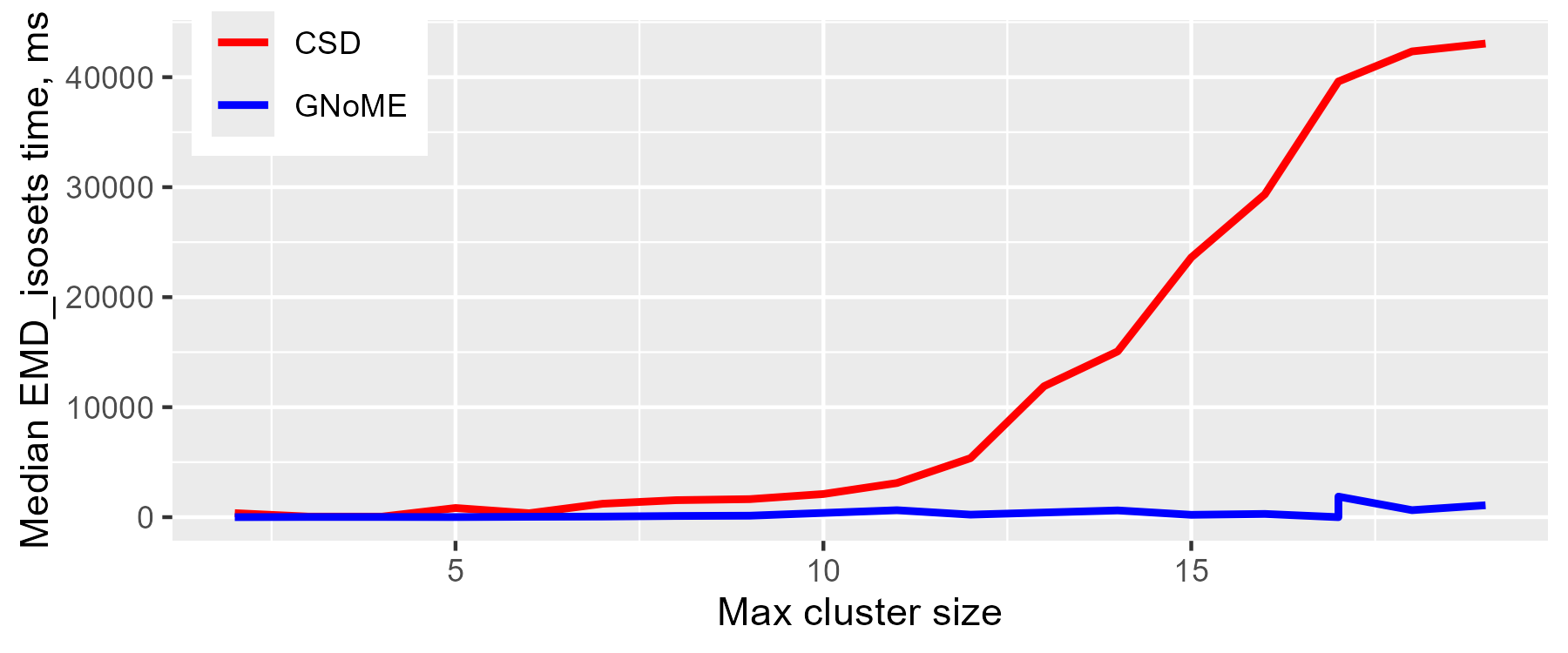}
\caption{
Median run times vs the max size of $\al$-clusters.
\textbf{Left}: invariants.
\textbf{Right}: metrics.
Times are faster for symmetric (inorganic) crystals in the GNoME than in the CSD. 
}
\label{fig:times} 
\end{figure}

The strongest past invariant PDD is based on distances and cannot distinguish mirror images. 
In the CSD, we {found four pairs that have identical PDDs but are mirror images shown in Fig.~{\ref{fig:Pauling}}~(right), distinguished by isosets with $\al\geq 1.5\angstrom$ in Fig.~{\ref{fig:mirror_images}}~(left).
For WODLOS vs XAWGAE and $\al=2$, the total time including isosets and EMD is about 4.3 seconds.
All experiments were on a modest CPU AMD Ryzen 5 5600X, 32GB RAM. 
Fig.~{\ref{fig:times}} summarizes times of new invariants and metrics, see Tables~A.1 and A.2 for details.
The supplementary materials include a Python code with instructions and full tables of metrics and run times for thousands of near-duplicates in the CSD and GNoME.}

\vspace*{-2mm}

The limitations of the $\EMD$ metric on isosets in Definition~\ref{dfn:EMD_isosets} are a slower running  time than for $\AMD,\PDD$ and the approximate (not exact) algorithm in Corollary~\ref{cor:approximate_EMD}, which are outweighed by the following crucial advantages.
\medskip

First, all past invariants could not distinguish infinitely many periodic sets (including all mirror images) under rigid motion, e.g. the real crystals in Fig.~\ref{fig:Pauling}.
The new continuous $\EMD$ fully solved  Problem~\ref{pro:metric}, which remained open at least since 1965 \cite{lawton1965reduced}.  
Second, because the proved error factor in the practical dimension $n=3$ is close to $4$, any near-duplicate crystals that differ by atomic deviations of up to $\ep$ 
have an exact distance $\EMD\leq 2\ep$ by main Theorem~\ref{thm:continuity} and hence an approximate distance up to about $8\ep$ by Corollary~\ref{cor:approximate_EMD}.
\smallskip

Any crystals that can be matched under rigid motion are recognizable since our approximation of $\EMD=0$ is also $0$.
Any approximate value $\de$ of $\EMD$ for real crystals $S,Q$ implies that all atoms of $S$ should be perturbed by at least $\de/8$ on average for a complete match with $Q$.
This result justifies filtering out near-duplicates \cite{zwart2008surprises} in all datasets to avoid machine learning on skewed data.
\smallskip

Future work can use the $\EMD$ to continuously quantify changes in material properties under perturbations of atoms and extend Problem~\ref{pro:metric} to metrics on finite or periodic sets of points under affine and projective transformations.
\smallskip

\textbf{Conclusions}.
Sections~\ref{sec:isosets} and \ref{sec:metrics} prepared the complexity results in section~\ref{sec:algorithms}: algorithms for computing and comparing isosets (Theorems~\ref{thm:compute_isoset}, Corollary~\ref{cor:compare_isosets}),  and approximating the new boundary tolerant metric $\BT$ (Theorem~\ref{thm:approximate_BT}), and EMD on isosets (Corollary~\ref{cor:approximate_EMD}). 
The proofs expressed polynomial bounds in terms of the motif \emph{size} $m=|S|$ of a periodic set $S$ because the input size of a Crystallographic Information File is linear in $m$, e.g. any lattice has $m=1$. 
\smallskip

The factors depending on the dimension and geometric complexity $\GC(S)$ are inevitable due to the curse of dimensionality and the infinite nature of crystals.
In practice, crystal symmetries reduce a motif to a smaller asymmetric part, which usually has less than 20 atoms even for large molecules in the CSD.
The lower bound via faster PDD invariants in Theorem~\ref{thm:lower_bound} justifies applying the algorithm of Corollary~\ref{cor:approximate_EMD} only for a final confirmation of near-duplicates. 
So the isosets finalized the hierarchy of the faster but incomplete invariants.
\smallskip

The main novelty is the boundary tolerant metric in Definition~\ref{dfn:tolerant_metric} that makes the complete invariant isoset Lipschitz continuous (Theorem~\ref{thm:continuity}) without extra parameters that are needed to smooth past descriptors such as powder diffraction patterns and atomic environments with fixed cut-off radii.
Since the isoset is the only Lipschitz continuous invariant whose completeness under isometry was proved for all periodic point sets in $\R^n$, the isoset was used to confirm the duplicates in the CSD and GNoME, see Tables~A.1 and~A.2. 
\smallskip

The resulting \emph{Crystal Isometry Principle} (CRISP) says that all non-isometric periodic crystals should have non-isometric sets of atomic centers and implies that all known and not yet discovered periodic crystals live in a common \emph{Crystal Isometry Space} (CRIS) whose first continuous maps appeared in \cite{bright2023geographic,bright2023continuous}.
\smallskip

We thank all reviewers for their valuable time and helpful comments.

\vspace*{-4mm}
\bibliography{near-duplicate-periodic-patterns}

\begin{thebibliography}{10}
\expandafter\ifx\csname url\endcsname\relax
  \def\url#1{\texttt{#1}}\fi
\expandafter\ifx\csname urlprefix\endcsname\relax\def\urlprefix{URL }\fi
\expandafter\ifx\csname href\endcsname\relax
  \def\href#1#2{#2} \def\path#1{#1}\fi

\bibitem{feynman1971feynman}
R.~Feynman, The Feynman lectures on physics, Vol.~1, 1971.

\bibitem{anosova2024importance}
O.~Anosova, V.~Kurlin, M.~Senechal, The importance of definitions in crystallography, IUCrJ 11 (2024) 453--463.
\newblock \href {https://doi.org/10.1107/S2052252524004056} {\path{doi:10.1107/S2052252524004056}}.

\bibitem{sacchi2020same}
P.~Sacchi, et~al., Same or different -- that is the question: identification of crystal forms from crystal data, CrystEngComm 22 (2020) 7170--7185.

\bibitem{widdowson2024continuous}
D.~Widdowson, V.~Kurlin, Continuous invariant-based maps of the {C}{S}{D}, Crystal Growth and Design 24 (2024) 5627–5636.

\bibitem{bimler2022better}
D.~Bimler, Better living through coordination chemistry: A descriptive study of a prolific papermill that combines crystallography and medicine (2022).
\newblock \href {https://doi.org/10.21203/rs.3.rs-1537438/v1} {\path{doi:10.21203/rs.3.rs-1537438/v1}}.

\bibitem{widdowson2022average}
D.~Widdowson, et~al., Average minimum distances of periodic sets - foundational invariants for mapping periodic crystals, MATCH 87 (2022) 529--559.

\bibitem{widdowson2025geographic}
D.~Widdowson, V.~Kurlin, Geographic-style maps with a local novelty distance help navigate in the materials space, Scientific Rep. 15 (2025) 27588.

\bibitem{cheetham2024artificial}
A.~K. Cheetham, R.~Seshadri, Artificial intelligence driving materials discovery?, Chemistry of Materials 36~(8) (2024) 3490–3495.

\bibitem{widdowson2025pointwise}
D.~Widdowson, V.~Kurlin, Pointwise distance distributions for detecting near-duplicates in large materials databases, SIAM Journal on Applied Mathematics (to appear, arxiv:2108.04798) (2025).

\bibitem{zwart2008surprises}
P.~Zwart, et~al., Surprises and pitfalls arising from (pseudo) symmetry, Acta Cryst. D 64 (2008) 99--107.

\bibitem{kurlin2024mathematics}
V.~Kurlin, Mathematics of 2-dimensional lattices, Foundations of Computational Mathematics 24 (2024) 805–863.

\bibitem{smith2024generic}
P.~Smith, V.~Kurlin, Generic families of finite metric spaces with identical or trivial 1{D} persistence, J Applied Comp. Topology 8 (2024) 839–855.

\bibitem{kurlin2024polynomial}
V.~Kurlin, Polynomial-time algorithms for continuous metrics on atomic clouds of unordered points, MATCH 91 (2024) 79--108.

\bibitem{widdowson2023recognizing}
D.~Widdowson, V.~Kurlin, Recognizing rigid patterns of unlabeled point clouds by complete and continuous isometry invariants with no false negatives and no false positives, in: CVPR, 2023, pp. 1275--1284.

\bibitem{kurlin2025complete}
V.~Kurlin, Complete and continuous invariants of 1-periodic sequences in polynomial time, SIAM J Mathematics of Data Science 7 (2025) 1643--1663.

\bibitem{bright2023geographic}
M.~Bright, A.~Cooper, V.~Kurlin, Geographic-style maps for 2-dimensional lattices, Acta Crystallographica Section A 79 (2023) 1--13.

\bibitem{bright2023continuous}
M.~Bright, A.~Cooper, V.~Kurlin, Continuous chiral distances for 2-dimensional lattices, Chirality 35 (2023) 920--936.

\bibitem{anosova2021isometry}
O.~Anosova, V.~Kurlin, An isometry classification of periodic point sets, in: LNCS (DGMM proceedings), Vol. 12708, 2021, pp. 229--241.

\bibitem{anosova2021introduction}
O.~Anosova, V.~Kurlin, Introduction to periodic geometry and topology, arxiv:2103.02749.

\bibitem{rubner2000earth}
Y.~Rubner, C.~Tomasi, L.~Guibas, The {E}arth {M}over's {D}istance as a metric for image retrieval, International J Computer Vision 40 (2000) 99--121.

\bibitem{chisholm2005compack}
J.~Chisholm, S.~Motherwell, Compack: a program for identifying crystal structure similarity using distances, J. Appl. Cryst. 38 (2005) 228--231.

\bibitem{widdowson2022resolving}
D.~Widdowson, V.~Kurlin, Resolving the data ambiguity for periodic crystals, Advances in Neural Information Processing Systems 35 (2022).

\bibitem{duneau1991bounded}
M.~Duneau, C.~Oguey, Bounded interpolations between lattices, Journal of Physics A: Mathematical and General 24~(2) (1991) 461.

\bibitem{carstens1999geometrical}
H.-G. Carstens, et~al., Geometrical bijections in discrete lattices, Combinatorics, Probability and Computing 8~(1-2) (1999) 109--129.

\bibitem{senechal1996quasicrystals}
M.~Senechal, Quasicrystals and geometry, CUP Archive, 1996.

\bibitem{mosca2020voronoi}
M.~Mosca, V.~Kurlin, Voronoi-based similarity distances between arbitrary crystal lattices, Crystal Research and Technology 55~(5) (2020) 1900197.

\bibitem{kawano2021classification}
S.~Kawano, J.~Mason, Classification of atomic environments via the {G}romov--{W}asserstein distance, Comp. Materials Science 188 (2021) 110144.

\bibitem{rass2024metricizing}
S.~Rass, S.~K{\"o}nig, S.~Ahmad, M.~Goman, Metricizing the euclidean space towards desired distance relations in point clouds, IEEE Transactions on Information Forensics and Security 19 (2024) 7304--7319.

\bibitem{lawton1965reduced}
S.~Lawton, R.~Jacobson, The reduced cell and its crystallographic applications, Tech. rep., Ames Lab, Iowa State University (1965).

\bibitem{edels2021}
H.~Edelsbrunner, T.~Heiss, V.~Kurlin, P.~Smith, M.~Wintraecken, The density fingerprint of a periodic point set, in: SoCG, 2021, pp. 32:1--32:16.

\bibitem{smith2022practical}
P.~Smith, V.~Kurlin, A practical algorithm for degree-k voronoi domains of three-dimensional periodic point sets, in: Lecture Notes in Computer Science (Proceedings of ISVC), Vol. 13599, 2022, pp. 377--391.

\bibitem{anosova2023density}
O.~Anosova, V.~Kurlin, Density functions of periodic sequences of continuous events, Journal of Mathematical Imaging and Vision 65 (2023) 689–701.

\bibitem{anosova2022density}
O.~Anosova, V.~Kurlin, Density functions of periodic sequences, in: LNCS Proceedings of Discrete Geometry and Mathematical Morphology, Vol. 13493, 2022, pp. 395--408.

\bibitem{delone1976local}
B.~Delone, N.~Dolbilin, M.~Shtogrin, R.~Galiulin, A local criterion for regularity of a system of points, in: DAN SSSR, Vol. 227, 1976, pp. 19--21.

\bibitem{dolbilin1998multiregular}
N.~Dolbilin, J.~Lagarias, M.~Senechal, Multiregular point systems, Discrete \& Computational Geometry 20~(4) (1998) 477--498.

\bibitem{mcmanus2025computing}
J.~McManus, V.~Kurlin, Computing the bridge length: the key ingredient in a continuous isometry classification of periodic point sets, Acta Cryst A 81 (2025).

\bibitem{dolbilin2016uniqueness}
N.~Dolbilin, A.~Magazinov, Uniqueness theorem for locally antipodal {D}elaunay sets, Proceedings of the Steklov Institute of Maths 294 (2016) 215--221.

\bibitem{brass2000testing}
P.~Brass, C.~Knauer, Testing the congruence of d-dimensional point sets, in: Proceedings of SoCG, 2000, pp. 310--314.

\bibitem{goodrich1999approximate}
M.~T. Goodrich, J.~S. Mitchell, M.~W. Orletsky, Approximate geometric pattern matching under rigid motions, T-PAMI 21 (1999) 371--379.

\bibitem{orlin1993faster}
J.~B. Orlin, A faster strongly polynomial minimum cost flow algorithm, Operations research 41~(2) (1993) 338--350.

\bibitem{shirdhonkar2008approximate}
S.~Shirdhonkar, D.~Jacobs, Approximate {E}arth {M}over’s {D}istance in linear time, in: Computer Vision and Pattern Recognition, 2008, pp. 1--8.

\bibitem{pauling1930crystal}
L.~Pauling, M.~Shappell, The crystal structure of bixbyite and the c-modification of the sesquioxides, Zeitschrift f{\"u}r Kristallographie-Cryst. Materials 75 (1930) 128--142.

\end{thebibliography}

\renewcommand{\thesection}{\Alph{section}}
\setcounter{section}{0}
\section{Appendix A: detailed proofs of all auxiliary results}
\label{sec:proofs}

This appendix includes more detailed and updated proofs of \cite[Lemmas~7 and 11-13]{anosova2021isometry}.
Fig.~\ref{fig:sections_diagram} visualizes the logical connections between main results.
\medskip

\begin{figure}[h!]
\includegraphics[height=55mm]{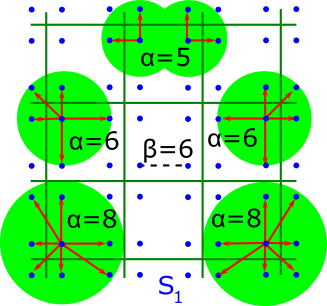}
\hspace*{2mm}
\includegraphics[height=55mm]{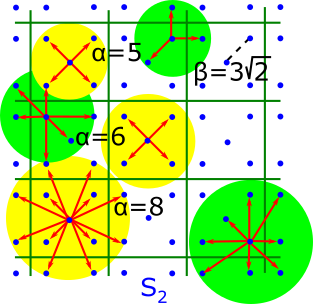}
\caption{
\textbf{Left}: 
in $\R^2$, the periodic point set $S_1$ has the square unit cell $[0,10)^2$ containing the four points $(2,2),(2,8),(8,2),(8,8)$, so $S_1$ isn't a lattice, but is 1-regular by Definition~\ref{dfn:m-regular}, and $\be(S_1)=6$.
All local $\al$-clusters of $S_1$ are isometric, shown by red arrows for $\al=5,6,8$, see Definition~\ref{dfn:local_cluster}.
\textbf{Right}: 
$S_2$ has the extra point $(5,5)$ in the center  of the cell $[0,10)^2$ and is 2-regular with $\be(S_2)=3\sqrt{2}$, so $S_2$ has green and yellow isometry types of $\al$-clusters. }
\label{fig:alpha-clusters}
\end{figure}

Fig.~\ref{fig:alpha-clusters}~(left) shows the 1-regular periodic set $S_1\subset\R^2$ whose all points (close to vertices of square cells) have isometric global clusters related by translations and rotations through $90^\circ,180^\circ,270^\circ$. 
The set $S_2$ has extra points at the centers of all square cells.
The local $\al$-clusters around these centers are not isometric to $\al$-clusters around the points close to cell vertices for any $\al\geq3\sqrt{2}$. 
\medskip

The 1-regular periodic point set $S_1$ in Fig.~\ref{fig:alpha-clusters} for any $p\in S_1$ has the symmetry group $\sym(S_1,p;\al)=\Or(\R^2)$ for $\al\in[0,4)$. 
Then $\sym(S_1,p;\al)$ stabilizes as $\Z^2$ with one reflection for $\al\geq 4$ as soon as $C(S_1,p;\al)$ includes one more point.

\begin{figure}[h!]
\includegraphics[height=17mm]{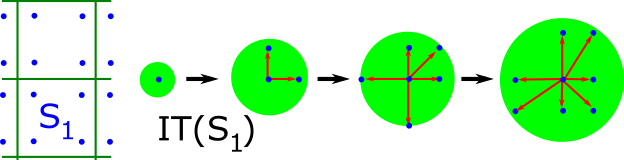}
\hspace*{2mm}
\includegraphics[height=17mm]{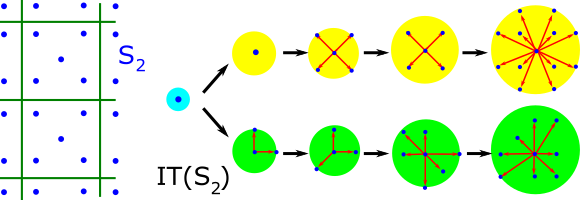}
\caption{
\textbf{Left}:The isotree $\IT(S_1)$ from Definition~\ref{dfn:isotree} of the 1-regular set $S_1$ in Fig.~\ref{fig:alpha-clusters} for any $\al\geq 0$ has one isometry class of $\al$-clusters under rotations.
\textbf{Right}: the isotree $\IT(S_2)$ of the 2-regular set $S_2$ in Fig.~\ref{fig:alpha-clusters} stabilizes with two non-isometric classes of $\al$-clusters for $\al\geq 4$.}
\label{fig:1-regular_set_isotree}
\end{figure}

\begin{figure}[h!]
\includegraphics[width=\textwidth]{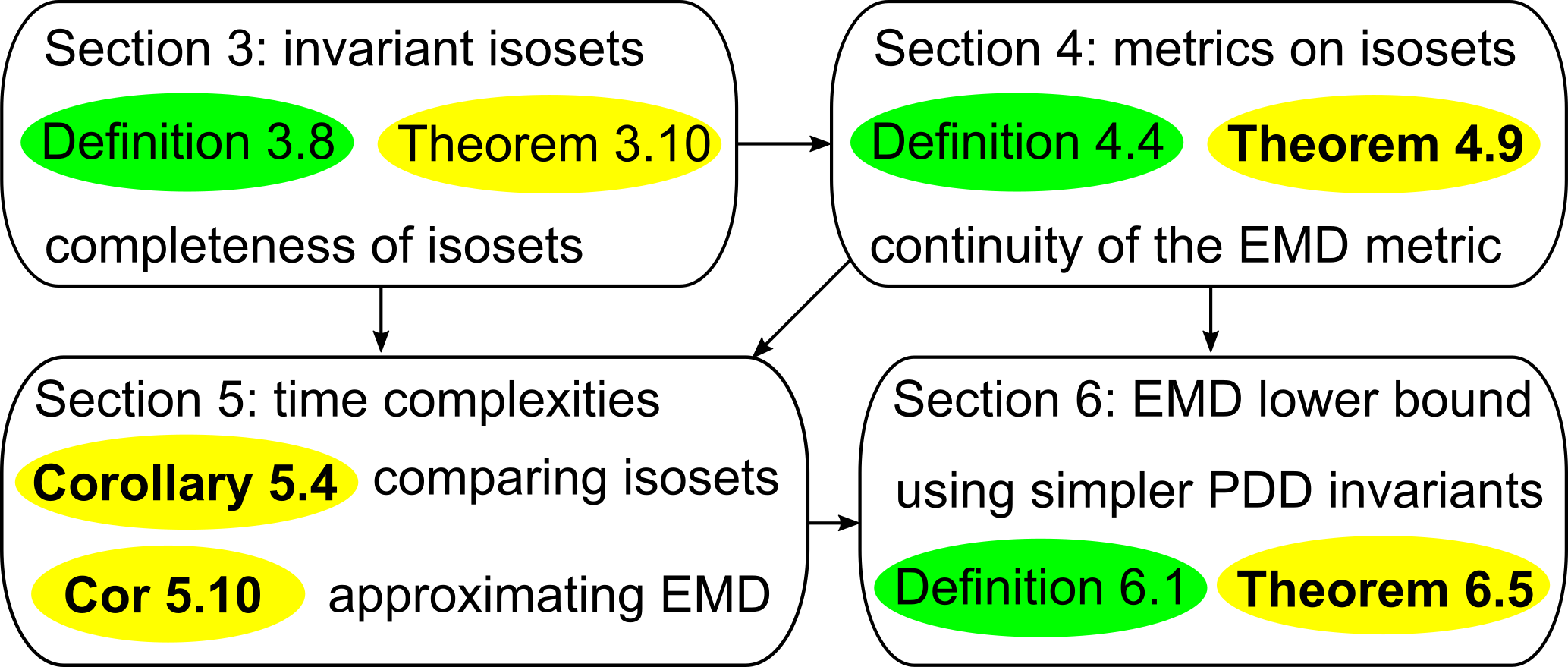}
\vspace*{-8mm}
\caption{Key definitions and main Theorems~\ref{thm:continuity},~\ref{thm:lower_bound} about the continuous metrics on complete invariant isosets with time complexities of algorithms in Corollaries~\ref{cor:compare_isosets},~\ref{cor:approximate_EMD}.}
\label{fig:sections_diagram}
\end{figure}

Fig.~\ref{fig:1-regular_set_isotree} 
illustrates the isosets for the periodic sets $S_1,S_2$ in Fig.~\ref{fig:alpha-clusters}.

\begin{lem}[isotree properties]
\label{lem:isotree} 
The isotree $\IT(S)$  
has the properties below:
\medskip

\noindent
(\ref{lem:isotree}a)
for $\al=0$, the $\al$-partition $P(S;0)$ consists of one class;
\medskip

\noindent
(\ref{lem:isotree}b) 
if $\al<\al'$, then $\sym(S,p;\al')\subseteq\sym(S,p;\al)$ for any point $p\in S$;
\medskip

\noindent
(\ref{lem:isotree}c) 
if $\al<\al'$, the $\al'$-partition $P(S;\al')$ \emph{refines} $P(S;\al)$, i.e. any $\al'$-equivalence class from $P(S;\al')$ is included into an $\al$-equivalence class from the partition $P(S;\al)$.
So the cluster count $|P(S;\al)|$ is non-strictly increasing 
in $\al$.
\bs
\end{lem}
\begin{proof}
(\ref{lem:isotree}a)
Let $\al\geq 0$ be smaller than the minimum distance $2r(S)$ betweens any points of $S$.
Then any cluster $C(S,p;\al)$ is the single-point set $\{p\}$.
All these 1-point clusters are isometric to each other.
So $|P(S;\al)|=1$ for $\al<2r(S)$.
\medskip

\noindent
(\ref{lem:isotree}b)
For any $p\in S$, the inclusion of clusters $C(S,p;\al)\subseteq C(S,p;\al')$ implies that any isometry $f\in\Or(\R^n;p)$ that isometrically maps the larger cluster $C(S,p;\al')$ to itself also maps the smaller cluster $C(S,p;\al)$ to itself.
Hence any element of $\sym(S,p;\al')\subseteq\Or(\R^n;p)$ belongs to $\sym(S,p;\al)$. 
\medskip

\noindent
(\ref{lem:isotree}c)
If points $p,q\in S$ are $\al'$-equivalent at the larger radius $\al'$, i.e. the clusters $C(S,p;\al')$ and $C(S,q;\al')$ are related by an isometry from $\Or(\R^n;p,q)$, then $p,q$ are $\al$-equivalent at the smaller radius $\al$.
Hence any $\al'$-equivalence class of points in $S$ is a subset of an $\al$-equivalence class in $S$.
\end{proof}

The proofs of Lemmas~\ref{lem:local_extension}, \ref{lem:global_extension}, and Theorem~\ref{thm:isoset_complete} follow \cite{delone1976local}, \cite[section~4]{dolbilin1998multiregular}, and extend to the oriented case by taking orientation-preserving isometries. 
Recall that $\Or(\R^n;p,q)$ denotes the set of all isometries of $\R^n$ that map $p$ to $q$.

\begin{lem}[local extension]
\label{lem:local_extension}
Let $S,Q\subset\R^n$ be periodic point sets and $\sym(S,p;\al-\be)=\sym(S,p;\al)$ for some point $p\in S$ and $\al>\be$.
Assume that there is an isometry $g\in\Or(\R^n;p,q)$ such that $g(C(S,p;\al))=C(Q,q;\al)$.
Let $f\in\Or(\R^n;p,q)$ be any isometry such that $f(C(S,p;\al-\be))=C(Q,q;\al-\be)$.
Then $f$ isometrically maps the larger clusters: $f(C(S,p;\al))=C(Q,q;\al)$.
\bs
\end{lem}
\begin{proof}
The composition $h=f^{-1}\circ g$ fixes $p$ and isometrically maps $C(S,p;\al-\be)$ to itself, so $h\in\sym(S,p;\al-\be)$.
The condition $\sym(S,p;\al-\be)=\sym(S,p;\al)$ implies that $h\in\sym(S,p;\al)$, so the isometry $h\in\Or(\R^n;p)$ isometrically maps the larger cluster $C(S,p;\al)$ to itself.
Then the given isometry $f=g\circ h^{-1}$ isometrically maps $C(S,p;\al)$ to $f(C(S,p;\al))=g(C(S,p;\al))=C(Q,q;\al)$. 
\end{proof}

\begin{lem}[global extension]
\label{lem:global_extension}
Let periodic point sets $S,Q\subset\R^n$ have a common stable radius  $\al$ satisfying Definition~\ref{dfn:stable_radius} for an upper bound $\be\geq\be(S),\be(Q)$.
Let $I(S;\al)=I(Q;\al)$ and $p\in S$, $q\in Q$ be any points with an isometry $f\in\Or(\R^n;p,q)$ such that $f(C(S,p;\al))=C(Q,q;\al)$.
Then $f(S)=Q$.
\bs
\end{lem}
\begin{proof}
To show that $f(S)\subset Q$, it suffices to check that the image $f(a)$ of any point $a\in S$ belongs to $Q$.
By Definition~\ref{dfn:bridge_length} the points $p,a\in S$ are connected by a sequence of points $p=a_0,a_1,\dots,a_k=a\in S$ such that the distances $|a_{i-1}- a_{i}|$ between any successive points have the upper bound $\be$ for $i=1,\dots,k$.
\medskip

We will prove that $f(C(S,a_k;\al))=C(Q,f(a_k);\al)$ by induction on $k$, where the base $k=0$ is given.
The induction step below goes from $i$ to $i+1$. 
\medskip

The ball $\bar B(a_i;\al)$ contains the smaller ball $\bar B(a_{i+1};\al-\be)$ around the closely located center $a_{i+1}$.
Indeed, since $|a_{i+1}-a_i|\leq\be$, the triangle inequality for the Euclidean distance implies that any point $a'_{i+1}\in\bar B(a_{i+1};\al-\be)$ with $|a'_{i+1}-a_i|\leq\al-\be$ satisfies $|a'_{i+1}-a_i|\leq |a'_{i+1}-a_{i+1}|+|a_{i+1}-a_i|\leq(\al-\be)+\be=\al$, so 
$\bar B(a_{i+1};\al-\be)\subset\bar B(a_i;\al)$.
Then the inductive assumption  $f(C(S,a_i;\al))=C(Q,f(a_i);\al)$ gives 
$f(C(S,a_{i+1};\al-\be))=f(C(S,a_i;\al))\cap f(\bar B(a_{i+1};\al-\be)) 
=C(Q,f(a_i);\al)\cap \bar B(f(a_{i+1});\al-\be)=C(Q,f(a_{i+1});\al-\be)$.
\medskip

Due to $I(S;\al)=I(Q;\al)$, the isometry class of $C(S,a_{i+1};\al)$ equals an isometry class of $C(Q,b_{i+1};\al)$ for some point $b_{i+1}\in Q$, i.e. there is an isometry $g\in\Or(\R^n;a_{i+1},b_{i+1})$ such that $g(C(S,a_{i+1};\al))=C(Q,b_{i+1};\al)$.
Since $f\circ g^{-1}\in\Or(\R^n;b_{i+1})$ isometrically maps $C(Q,b_{i+1};\al-\be)$ to $C(Q,f(a_{i+1});\al-\be)$, the points $b_{i+1},f(a_{i+1})\in Q$  are in the same $(\al-\be)$-equivalence class of $Q$.
\medskip

By condition~(\ref{dfn:stable_radius}a), the splitting of the periodic point set $Q\subset\R^n$ into $\al$-equivalence classes coincides with its splitting into $(\al-\be)$-equivalence classes. 
Hence the points $b_{i+1},f(a_{i+1})\in Q$ are in the same $\al$-equivalence class of $Q$.
Then $C(Q,f(a_{i+1});\al)$ is isometric to $C(Q,b_{i+1};\al)=g(C(S,a_{i+1};\al))$.\medskip

Now we can apply Lemma~\ref{lem:local_extension} for $p=a_{i+1},q=f(a_{i+1})$ and conclude that the given isometry $f$, which satisfies $f(C(S,a_{i+1};\al-\be))=C(Q,f(a_{i+1});\al-\be)$, isometrically maps the larger clusters: $f(C(S,a_{i+1};\al))=C(Q,f(a_{i+1});\al)$.
The induction step is finished.
The inclusion $f^{-1}(Q)\subset S$ is proved similarly.
\end{proof}

\begin{lem}[all stable radii $\al\geq\al(S)$]
\label{lem:stable_radius}
If $\al$ is a stable radius of a periodic point set $S\subset\R^n$, then so is any  larger radius $\al'>\al$.
Then all stable radii form the interval $[\al(S),+\infty)$, where $\al(S)$ is the minimum stable radius of $S$.
\bs
\end{lem}
\begin{proof}
Due to Lemma~(\ref{lem:isotree}bc), conditions (\ref{dfn:stable_radius}ab) imply that the $\al'$-partition $P(S;\al')$ and the symmetry groups $\sym(S,p;\al')$ remain the same for all $\al'\in[\al-\be(S),\al]$, where $\be(S)$ is the bridle length. 
We need to show that they remain the same for any $\al'>\al$ and will apply Lemma~\ref{lem:global_extension} for $S=Q$ and $\be=\be(S)$.
Let points $p,q\in S$ be $\al$-equivalent, i.e. there is an isometry $f\in\Or(\R^n;p,q)$ such that $f(C(S,p;\al))=C(S,q;\al)$.
By Lemma~\ref{lem:global_extension}, $f$ isometrically maps the full set $S$ to itself.
Then all larger $\al'$-clusters of $p,q$ are matched by $f$, so $p,q$ are $\al'$-equivalent  and $P(S;\al)=P(S,\al')$.
Similarly, any isometry $f\in\sym(S,p;\al)$ by Lemma~\ref{lem:global_extension} for $S=Q$ and $p=q$, isometrically maps the full set $S$ to itself.
Then $\sym(S,p;\al')$ coincides with $\sym(S,p;\al)$ for any $\al'>\al$. 
\end{proof}
 
 \begin{figure}[h!]
\includegraphics[width=\textwidth]{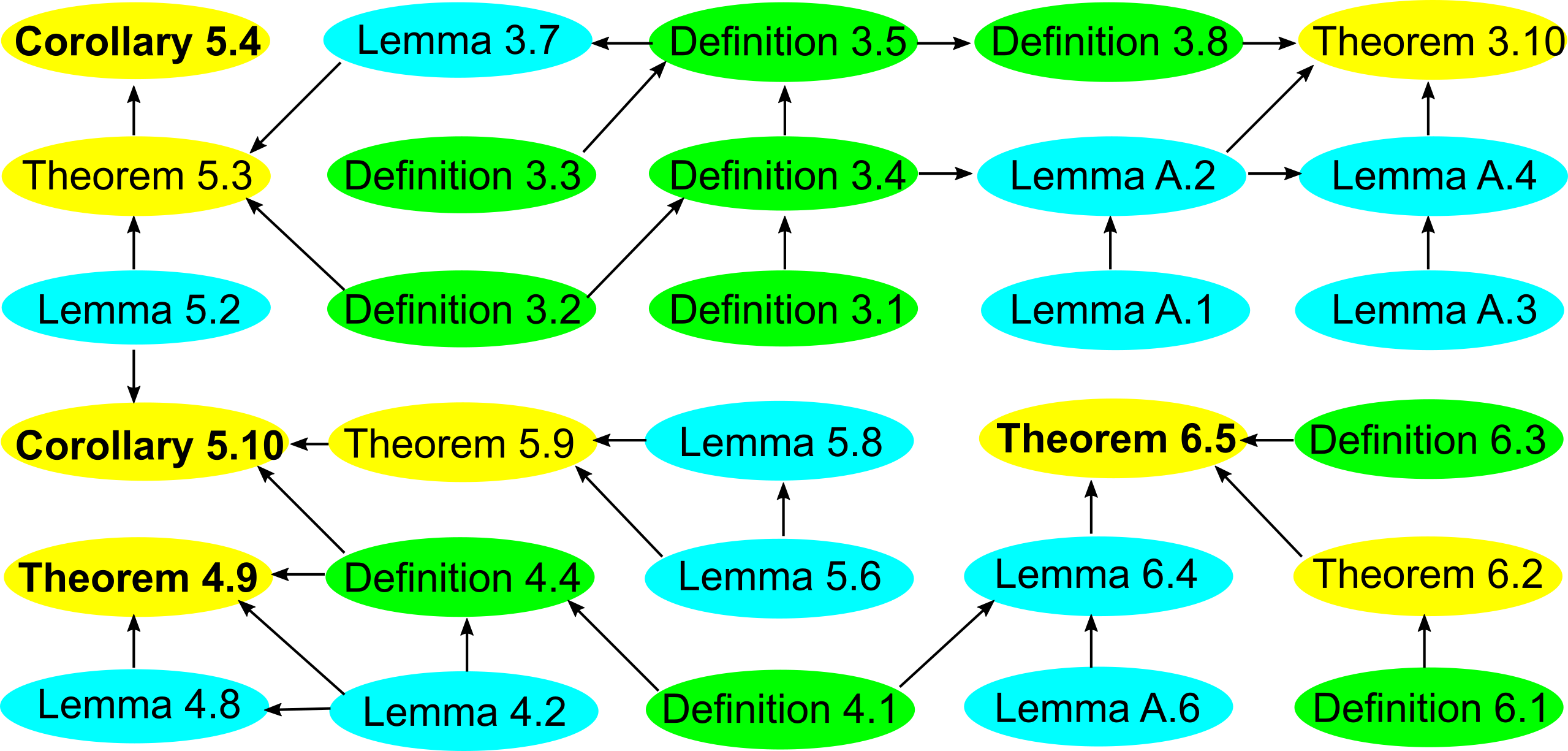}
\vspace*{-6mm}
\caption{Logical steps towards main Theorems~\ref{thm:continuity},~\ref{thm:lower_bound} and Corollaries~\ref{cor:compare_isosets} and~\ref{cor:approximate_EMD}.}
\label{fig:results_diagram}
\end{figure}

\begin{proof}[Proof of Theorem~\ref{thm:isoset_complete}]
The part \emph{only if} $\Rightarrow$.
Let $f$ be an isometry of $\R^n$, which isometrically maps one periodic point set $S$ to another $Q$.
For any point $p$ in a motif $M(S)$ of $S$, the image $f(p)\in Q$ is equivalent to a unique point $g(p)$ in a motif $M(Q)$ of $Q$ modulo a translation along a vector from the lattice of $Q$.
\medskip

Then, for any $p\in M(S)$ and $\al\geq 0$, the clusters $C(S,p;\al)$ and $C(Q,g(p);\al)$ are related by an isometry of $\R^n$.
Hence the bijection $g:M(S)\to M(Q)$ induces a bijection $I(S;\al)\to I(Q;\al)$ between all isometry classes with weights.
\medskip

The part \emph{if} $\Leftarrow$.
Fix a point $p\in S$. 
The cluster $C(S,p;\al)$ represents a class $\si\in I(S;\al)$.
Due to $I(S;\al)=I(Q;\al)$, the class $\si$ equals some $\xi\in I(Q;\al)$.
Hence there is an isometry $f$ of $\R^n$ such that the cluster $f(C(S,p;\al))=C(Q,f(p);\al)$ represents $\xi$. 
By Lemma~\ref{lem:global_extension}, $f$ isometrically maps $S$ to $Q$. 
\end{proof}

\begin{proof}[Proof of Lemma~\ref{lem:tolerant_metric}]
Since the set $\Or(\R^n;p,q)$ is compact, the minimum $\ep\geq 0$ is achieved in the inclusions from (\ref{dfn:tolerant_metric}b) for some isometries $f\in\Or(\R^n;p,q)$ and $g\in\Or(\R^n;q,p)$.
Then, for any clusters $C(S,\ti p;\al)$ and $C(Q,\ti q;\al)$ isometric to $C(S,p;\al)$ and $C(Q,q;\al)$ via $f_S\in\Or(\R^n;\ti p,p)$ and $g_Q\in\Or(\R^n;\ti q,q)$, respectively, the same minimum $\ep\geq 0$ is achieved
in the following inclusions (and vice versa), which proves the independence of $\BT(\si,\xi)$ under a choice of clusters.
\smallskip

\noindent
$C(Q,\ti q;\al-\ep)\subseteq \ti  f (C(S,\ti p;\al))+\bar B(0;\ep)$ for 
$\ti f=g_Q^{-1}\circ f\circ f_S\in\Or(\R^n;\ti p,\ti q)$, and 
\smallskip

\noindent
$C(S,\ti p;\al-\ep)\subseteq \ti g(C(Q,\ti q;\al))+\bar B(0;\ep)$ for 
$\ti g=f_S^{-1}\circ g\circ g_Q\in\Or(\R^n;\ti p,\ti q)$.
\medskip

Now we prove the coincidence axiom.
By Definition~\ref{dfn:tolerant_metric}, $\BT(\si,\xi)=0$ means that
some representatives of given classes $\si,\xi$ satisfy 
$C(Q,q;\al)\subseteq f(C(S,p;\al))$ for some $f\in\Or(\R^n;p,q)$ and
$C(S,p;\al)\subseteq g(C(Q,q;\al))$ for some $g\in\Or(\R^n;q,p)$.
Combining these inclusions, we get $C(Q,q;\al)\subseteq f\circ g(C(Q,q;\al))$.
Since $f\circ g\in \Or(\R^n;q)$ is an isometry fixing $q$ and both clusters in the inclusion above consist of the same number of points the surjection $a\mapsto f\circ g(a)$ for $a\in C(Q,q;\al)$ should bijective, so $C(Q,q;\al)=f\circ g(C(Q,q;\al))$.
Then the initial inclusions are equalities.
Hence $C(S,p;\al),C(Q,a;\al)$ are related by the isometry $f\in\Or(\R^n;p,q)$, so $\si=\xi$.
The symmetry axiom holds because the inclusions in condition (\ref{dfn:tolerant_metric}b) are symmetric to each other under swapping the arguments.
\medskip

To prove the triangle inequality, let clusters $C(S,p;\al)$, $C(Q,f(p);\al)$, $C(T,g\circ f(p);\al)$ represent $\si,\xi,\zeta$, respectively, so that $\ep_1=\BT(\si,\xi)$ and $\ep_2=\BT(\xi,\zeta)$ are achieved for inclusions
$C(Q,f(p);\al-\ep_1)\subseteq f(C(S,p;\al))+\bar B(0;\ep_1)$ and 
$C(T,g\circ f(p);\al-\ep_2)\subseteq g(C(Q,f(p);\al))+\bar B(0;\ep_2)$ for isometries $f,g$ of $\R^n$.
\medskip

The last inclusion gives $C(T,g\circ f(p);\al-\ep_1-\ep_2)\subseteq g(C(Q,f(p);\al-\ep_1))+\bar B(0;\ep_2)$ because we can reduce the radius $\al$ in the cluster $g(C(Q,f(p);\al))$ to $\al-\ep_1$.
Indeed, if a point $t\in C(T,g\circ f(p);\al-\ep_1-\ep_2)$ is covered by a closed ball $\bar B(q;\ep_2)$ for some $q\in g(C(Q,f(p);\al))$, then $|q-t|\leq\ep_2$ and
$$|q-g\circ f(p)|\leq |q-t|+|t-g\circ f(p)|\leq\ep_2+(\al-\ep_1-\ep_2)=\al-\ep_1.$$
Hence $q$ belongs to the smaller cluster $g(C(Q,f(p);\al-\ep_1)$ as required.
Now we apply the isometry $g$ to the inclusion $C(Q,f(p);\al-\ep_1)\subseteq f(C(S,p;\al))+\bar B(0;\ep_1)$ to get  
$C(T,g\circ f(p);\al-\ep_1-\ep_2)\subseteq 
g(C(Q,f(p);\al-\ep_1))+\bar B(0;\ep_2)\subseteq 
g\circ f(C(S,p;\al))+\bar B(0;\ep_1+\ep_2)$ as $(q+\bar B(0;\ep_1))+\bar B(0;\ep_2)=q+\bar B(0;\ep_1+\ep_2)$.
\medskip

Swapping the roles of $S,T$ in the arguments above, we similarly prove that if 
$C(S,p;\al-\ep_1)\subseteq f(C(Q,f^{-1}(p);\al))+\bar B(0;\ep_1)$ and 
$C(Q,f^{-1}(p);\al-\ep_2)\subseteq g(C(T,g^{-1}\circ f^{-1}(p);\al))+\bar B(0;\ep_2)$ for some isometries $f,g$ of $\R^n$, then
$$C(S,p;\al-\ep_1-\ep_2)\subseteq f\circ g(C(T,g^{-1}\circ f^{-1}(p);\al))+\bar B(0;\ep_1+\ep_2).$$

Definition~\ref{dfn:tolerant_metric} implies that $\BT(\si,\zeta)\leq \ep_1+\ep_2=\BT(\si,\xi)+\BT(\xi,\zeta)$.
\end{proof}

\begin{exa}[detailed computations for Example~\ref{exa:square_vs_hexagon}]
\label{exa:square_vs_hexagon_details}
Fig.~\ref{fig:square_vs_hexagon} shows the stable 2-clusters $C(\La_4,0;2)$ and $C(\La_6,0;2)$ of the square ($\La_4$) and hexagonal ($\La_6$) lattices.
Without rotations, the 1st picture of Fig.~\ref{fig:square_vs_hexagon} shows the directed Hausdorff distance $d_{\vec H}=\sqrt{(1-\frac{\sqrt{3}}{2})^2+(\frac{1}{2})^2}=\sqrt{2-\sqrt{3}}\approx 0.52$ between clusters with the added boundary circle $\bd B(0;2)$.
Due to high symmetry, it suffices to consider rotations 
of the square vertex $(1,1)$ 
for angles $\ga\in[45^{\circ},60^{\circ}]$ because all other ranges can be isometrically mapped to this range for another vertex of the square.
We find the squared distances $s_1(\ga)$ and $s_2(\ga)$ from the vertex $(\sqrt{2}\cos\ga,\sqrt{2}\sin\ga)$ rotated from $(1,1)$ at $\ga=45^\circ$ through the angle $\ga-45^{\circ}$ to its closest neighbors $(\frac{1}{2},\frac{\sqrt{3}}{2})$ and $(\frac{3}{2},\frac{\sqrt{3}}{2})$ in $C(\La_6,0;2)$.
\medskip

\noindent
$s_1(\gamma)=\left|(\sqrt{2}\cos\ga,\sqrt{2}\sin\ga)-\Big(\frac{1}{2},\frac{\sqrt{3}}{2}\Big)\right|^2
=
\left(\sqrt{2}\cos\ga-\frac{1}{2}\right)^2+\left(\sqrt{2}\sin\ga-\frac{\sqrt{3}}{2}\right)^2=3-\sqrt{2}\cos\ga-\sqrt{6}\sin\ga,\; 
\frac{ds_1}{d\ga}=\sqrt{2}\sin\ga-\sqrt{6}\cos\ga=0,\;
\tan\ga=\sqrt{3}, \; 
\ga=60^{\circ}, 
s_1=(\sqrt{2}-1)^2 \mbox{ is minimal for the points in } y=\sqrt{3}x \mbox{ at distances } 1,\sqrt{2} \mbox{ from } 0.$
\medskip

\noindent
$s_2(\gamma)=\left|(\sqrt{2}\cos\ga,\sqrt{2}\sin\ga)-\Big(\frac{3}{2},\frac{\sqrt{3}}{2}\Big)\right|^2=
\left(\sqrt{2}\cos\ga-\frac{3}{2}\right)^2+\left(\sqrt{2}\sin\ga-\frac{\sqrt{3}}{2}\right)^2
=5-3\sqrt{2}\cos\ga-\sqrt{6}\sin\ga,\;
\frac{ds_2}{d\ga}=3\sqrt{2}\sin\ga-\sqrt{6}\cos\ga=0,
\ga=30^{\circ}, 
s_2=(\sqrt{3}-\sqrt{2})^2 \mbox{ is minimal for the points in } y=\frac{x}{\sqrt{3}} \mbox{ at distances } \sqrt{2},\sqrt{3} \mbox{ from } 0.$ 
\medskip

It might look that the second minimum is smaller.
However, for the angle $\ga=30^{\circ}$, another vertex $(-1,1)$ rotated through $\ga-45^\circ=-15^{\circ}$ has distance $\sqrt{2}-1$ to its closest neighbor $(-\frac{1}{2},\frac{\sqrt{3}}{2})\in C(\La_6,0;2)$.
For any angle $\ga\in[45^\circ,60^\circ]$, the second function has the minimum $s_2(45^{\circ})=2-\sqrt{3}=d_{\vec H}^2$ in the 1st picture of Fig.~\ref{fig:square_vs_hexagon}.
Hence the vertex $(1,1)$ has the minimum distance $\sqrt{2}-1\approx 0.41<\sqrt{2-\sqrt{3}}\approx 0.52$ in the 3rd picture of Fig.~\ref{fig:square_vs_hexagon}.
All other points of the square cluster $C(\La_4,0;2)$ are even closer to their neighbors in $C(\La_6,0;2)$.
For example, the point $(1,0)$ rotated by $15^\circ$ has the distance to $(1,0)$ equal to $\sqrt{(\cos 15^\circ-1)^2+\sin^2 15^\circ}\approx 0.26$.
The final picture in Fig.~\ref{fig:square_vs_hexagon} confirms that all points of the hexagonal cluster $C(\La_6,0;2)$ are covered by the $(\sqrt{2}-1)$-offset of 
$C(\La_4,0;2)$ and the boundary circle.
So 
$\BT=\sqrt{2}-1\approx 0.41$.
\bs
\end{exa}

\begin{proof}[Proof of Lemma~\ref{lem:EMD_metric}]
The symmetry axiom holds because Definition~\ref{dfn:EMD_isosets} is symmetric under swapping $S,Q$.
The triangle axiom was proved for any weighted distributions in \cite[Appendix~A]{rubner2000earth}.
We will prove that if $\EMD( I(S;\al), I(Q;\al) )=0$ then $I(S;\al)=I(Q;\al)$ as isosets.
Indeed, $\sum\limits_{i=1}^{m(S)} \sum\limits_{j=1}^{m(Q)} f_{ij} \BT(\si_i,\xi_j)=0$ means that, for any $i,j$, if $f_{ij}>0$ then $\BT(\si_i,\xi_j)=0$, so $\si_i=\xi_j$ by the coincidence axiom of $\BT$ from Lemma~\ref{lem:tolerant_metric}(a).
Hence any flow $f_{ij}>0$ is always between equal isometry classes.
The conditions on weights of $\si_i,\xi_j$ in Definition~\ref{dfn:EMD_isosets} imply that every class $\si_i$ should `flow' to its equal class $\xi_j$ of the same weight.
These flows define a bijection $I(S;\al)\to I(Q;\al)$ respecting all weights.
\end{proof}


\begin{proof}[Proof of Lemma~\ref{lem:common_lattice}]
Choose the origin $0\in\R^n$ at a point of $S$.
Applying translations, we can assume that primitive unit cells $U(S),U(Q)$ of the given periodic sets $S,Q$ have a vertex at the origin $0$.
Then $S=\La(S)+(U(S)\cap S)$ and $Q=\La(Q)+(U(Q)\cap Q)$, where 
$\La(S),\La(Q)$ are lattices of $S,Q$, respectively.
\medskip

We are given that every point of $Q$ is $d_B(S,Q)$-close to a point of $S$, where the bottleneck distance $d_B(S,Q)$ is strictly less than the packing radius $r(Q)$.  
\medskip

Assume towards contradiction that $S,Q$ have no common lattice.
Then there is a point $p\in\La(S)$ whose all integer multiples $kp\in\La(S)$ do not belong to $\La(Q)$ for $k\in\Z-\{0\}$.
Any such multiple $kp$ can be translated by a vector of $\La(Q)$ to a point $q(k)$ in the unit cell $U(Q)$ so that $kp\equiv q(k)\pmod{\La(Q)}$.
Since the cell $U(Q)$ contains infinitely many points $q(k)$,
one can find a pair $q(i)\neq q(j)$ at a distance less than $\de=r(Q)-d_B(S,Q)>0$.
For any $m\in\Z$, the following points are equivalent modulo (translations along the vectors of) the lattice $\La(Q)$.
$$q(i+m(j-i))\equiv (i+m(j-i))p= 
 ip + m(jp-ip)\equiv 
 q(i) + m(q(j)-q(i)).$$
These points for $m\in\Z$ lie in a straight line with gaps $|q(j)-q(i)|<\de$.
The open balls with the packing radius $r(Q)$ and centers at all points of $Q$ do not overlap.
Hence all closed balls with the radius $d_B(S,Q)<r(Q)$ and the same centers are at least $2\de$ away from each other.
Due to $|q(j)-q(i)|<\de=r(Q)-d_B(S,Q)$, there is $m\in\Z$ such that $q(i) + m(q(j)-q(i))$ is outside the union $Q+\bar B(0;d_B(S,Q))$ of all these smaller balls.
Then $q(i) + m(q(j)-q(i))$ has a distance more than $d_B(S,Q)$ from any point of $Q$. 
The translations along all vectors of the lattice $\La(Q)$ preserve the union of balls $Q+\bar B(0;d_B(S,Q))$.
Then the point $(i+m(j-i))p\in S$, which is equivalent to $q(i) + m(q(j)-q(i))$ modulo $\La(Q)$, has a distance more than $d_B(S,Q)$ from any point of $Q$.
This conclusion contradicts Definition~\ref{dfn:Hausdorff+bottleneck}(b) of the bottleneck distance $d_B(S,Q)$. 
\end{proof}

\begin{proof}[Proof of Lemma~\ref{lem:compute_cluster}]
To find all points in $C(S,p;\al)$, we will extend $U$ by adding adjacent cells in`spherical' shells around $U$.
After considering the initial cell $U$ with a basis $\vec v_1,\dots,\vec v_n$, we take $3^n-1$ cells $U+\vec v$ for vectors $\vec v=\sum\limits_{i=1}^n c_i \vec v_i\in\La-\{0\}$ with integer coordinates $c_i\in\{-1,0,1\}$.
The next `spherical' shell consists of $5^n-3^n$ cells $U+\vec v$ and so on.
For any shifted cell $U+\vec v$ with $v\in\La$, if all vertices have distances more than $\al$ to $p$, this cell is discarded.
Otherwise, we check if any translated points $M+\vec v$ are within the closed ball $\bar B(p;\al)$ of radius $\al$.
The upper union $\bar U=\bigcup \{(U+\vec v) :  v\in\La, (U+\vec v)\cap \bar B(p;\al)\neq\emptyset\}$
consists of $\frac{\vol[\bar U]}{\vol[U]}$ cells and is contained in the larger ball $B(p;\al+d)$, because any shifted cell $U+\vec v$ within $\bar U$ has the longest diagonal $d$ and intersects $B(p;\al)$.
Since each $U+\vec v$ contains $m$ points of $S$, we check at most 
$m\frac{\vol[\bar U]}{\vol[U]}$ points.
So
$$|C(S,p;\al)|\leq 
m\frac{\vol[\bar U]}{\vol[U]}\leq 
m\frac{\vol[B(p;\al+d)]}{\vol[U]}=
m\frac{(\al+d)^n V_n}{\vol[U]}=\nu(U,\al,n)m,$$
where $\nu(U,\al,n)=\frac{(\al+d)^n V_n}{\vol[U]}$.
We will estimate $\nu(U,\al,n)$ using the upper bound $\al\leq (L+m+1)2R(S)$ from Lemma~\ref{lem:upper_bounds}(c).
Since the longest diagonal has the upper bound $2R(S)\geq d$ because the closed balls with the radius $\frac{d}{2}$ and centers at the vertices of a unit cell $U$ cover $U$, so $\al+d\leq (L+m+2)2R(S)$.
\medskip

Since $\Ga(\frac{n}{2}+1)
=\sqrt{\pi}\frac{(2n-1)(2n-3)\dots 1}{2^n}
=\sqrt{\pi}\frac{(2n)(2n-1)(2n-2)(2n-3)\dots 1}{2^{n}(2n)(2n-2)\dots 2}=\sqrt{\pi}\frac{(2n)!}{2^{2n}n!}$,
the volume of the unit ball becomes
$V_n=\frac{\pi^{n/2}}{\Ga(\frac{n}{2}+1)}=(\sqrt{\pi})^{n-1}\frac{2^{2n}n!}{(2n)!}$.
The bounds 
$\sqrt{2\pi n}\left(\frac{n}{e}\right)^n \exp(\frac{1}{12n+1}) < n!< \sqrt{2\pi n}\left(\frac{n}{e}\right)^n \exp(\frac{1}{12n})$, 
imply that
$V_n=(\sqrt{\pi})^{n-1}\frac{2^{2n}n!}{(2n)!}\leq 
\frac{(\sqrt{\pi})^{n-1}}{\sqrt{2}}2^{2n}
 \left(\frac{n}{e}\right)^n \left(\frac{e}{2n}\right)^{2n} \exp(\frac{1}{12n}-\frac{1}{24n+1})\leq
\frac{\exp(\frac{1}{22})}{\sqrt{2\pi}} \left(\frac{e\sqrt{\pi}}{n}\right)^n$ because
$\frac{1}{12n}-\frac{1}{24n+1}=\frac{(24n+1)-12n}{12n(24n+1)}\leq\frac{12n+1}{12\times 24 n}\leq\frac{13n}{12\times 24n}<\frac{1}{22}$ for $n\geq 1$.
Then 
$V_n\leq \frac{\exp(\frac{1}{22})}{\sqrt{2\pi}} \left(\frac{e\sqrt{\pi}}{n}\right)^n$ implies that
$\nu(U,\al,n)=\frac{(\al+d)^n V_n}{\vol[U]}
\leq \frac{((L+m+2)2R(S))^n V_n}{\vol[U]}
\leq\frac{\exp(\frac{1}{22})}{\sqrt{2\pi}}\frac{((L+m+2)2R(S)e\sqrt{\pi}/n)^n}{\vol[U]}
\leq \frac{(10(L+m+2)R(S)/n)^n}{2\vol[U]}=\GC(S)$ as required.
\end{proof}

\begin{proof}[Proof of Lemma~\ref{lem:max-min_formula}]
The directed distance $d_{\vec R}(C\cup\bd\bar B(0;\al),D\cup\bd\bar B(0;\al))$ is the minimum $\ep\in[0,\al]$ such that, for some $f\in\Or(\R^n)$, all points of $C\cap B(0;\al-\ep)$ are covered by $f(D)+\bar B(0;\ep)$ as all points of $C\setminus B(0;\al-\ep)$ are $\ep$-close to the boundary $\bd\bar B(0;\al)$.
Let $j\in\{1,\dots,k\}$ be the largest index so that $|p_j|<\al-\ep$.
Then $C\cap B(0;\al-\ep)=\{p_1,\dots,p_j\}$ and $d_{\vec R}(\{p_1,\dots,p_i\}, D)\leq d_{\vec R}(C\cap B(0;\al-\ep),D)\leq\ep$ for all $i=1,\dots,j$. 
By the above choice of $j$, if $j<i\leq k$ then $\al-|p_i|\leq\ep$.
Hence, for all $i=1,\dots,k$, both terms in the minimum $\min\{\; \al-|p_i|,\; d_{\vec R}(\{p_1,\dots,p_i\}, D)\; \}$ are at most $\ep$.
Then $d_{\vec M}(C,D)\leq\ep$.
\medskip

Using the brief notation $d_{\vec M}=d_{\vec M}(C,D)$, to prove the converse inequality $d_{\vec R}(C\cup\bd\bar B(0;\al),D\cup\bd\bar B(0;\al))\leq d_{\vec M}$, we check below that $C\cap\bar B(0;\al-d_{\vec M})$ is covered by $f(D)+\bar B(0;d_{\vec M})$ for some $f\in\Or(\R^n)$ or, equivalently, the inequality $d_{\vec R}(C\cap\bar B(0;\al-d_{\vec M}),D)\leq d_{\vec M}$ holds.
Let $j\in\{1,\dots,k\}$ be the largest index so that $|p_j|\leq\al-d_{\vec M}$.
Since $\al-|p_j|\geq d_{\vec M}$ and 
$\min\{\; \al-|p_j|,\; d_{\vec R}(\{p_1,\dots,p_j\}, D)\; \}$ $\leq\max\limits_{i=1,\dots,k}\min\{ \al-|p_i|,\; d_{\vec R}(\{p_1,\dots,p_i\}, D) \}=d_{\vec M}$, the term $d_{\vec R}(\{p_1,\dots,p_j\}, D)$ in the minimum above is at most $d_{\vec M}$. 
Due to $C\cap\bar B(0;\al-d_{\vec M})=\{p_1,\dots,p_j\}$, we get 
$d_{\vec R}(C\cap\bar B(0;\al-d_{\vec M}),D)\leq d_{\vec M}$, which proves the required equality.
\end{proof}

\begin{proof}[Proof of Lemma~\ref{lem:rot-inv_distance}]
Let $q_1\in D$ be a point that has a maximum distance to the origin $0\in\R^n$.
If there are several points at the same maximum distance, choose any of them.
Similar choices below do not affect the estimates.
For any $1<i<n$, let $q_i$ be a point of $D$ that has a maximum perpendicular distance to the linear subspace spanned by the previously defined vectors $\vec q_1,\dots,\vec q_{i-1}$.
\medskip

The key idea is to replace the minimization in $d_{\vec R}(C_j,D)$ over infinitely many $f\in\Or(\R^n)$ by a finite minimization over compositions $f_{n-1}[s_{n-1}]\circ\dots\circ f_1[s_1]\in\Or(\R^n)$ depending on finitely many unknown points $s_1,\dots,s_{n-1}\in C_j$, which will be exhaustively checked in time $O(\ar{C}^{n-1})$.
\medskip

If the point $q_1$ belongs to the infinite straight line $L(s_1)$ through the points $s_1$ and $0$, then set $f_1[s_1]$ to be the identity map.
Otherwise, let $f_1[s_1]\in\SO(\R^n)$ fix the linear subspace orthogonal to the plane spanned by $\vec s_1,\vec q_1$, and then rotate the point $q_1$ to $L(s_1)$ through the smallest possible angle.
Since $q_1$ is a furthest point of $D$ from the origin $0$ and $|s_1-q_1|\leq d_j$, the rotation $f_1[s_1]$ moves $q_1$ and hence any other point of $D$ by at most $d_j$.
Then any point in $f_1[s_1](D)$ is at most $2d_j$ away from its closest neighbor in $C_j$, so $d_{\vec H}(C_j,f_1[s_1](D))\leq 2d_j$.
\medskip

For any $1<i<n$, if $\vec q_i$ belongs to the linear subspace $L(q_1,\dots,q_{i-1},s_i)$ spanned by $\vec q_1,\dots,\vec q_{i-1},\vec s_i$, set $f_i[s_i]$ to be the identity map.
Else let the rotation $f_i[s_i]\in\SO(\R^n)$ fix the linear subspace orthogonal to $\vec q_1,\dots,\vec q_{i},\vec s_i$, and rotate $q_i$ to $L(q_1,\dots,q_{i-1},s_i)$ through the smallest possible angle.
Since $f_1[s_1](q_2)$ is at most $2d_j$ away from $s_2\in C_j$, the map $f_2[s_2]$ moves $f_1[s_1](q_2)$ and hence any other point of $f_1[s_1](D)$, by at most $2d_j$.
Since $q_2$ had a maximum perpendicular distance from the line through $\vec q_1$, the composition $f_2[s_2]\circ f_1[s_1]$ moves any point of $D$ by at most $d_j+2d_j=3d_j$.
For $2<i<n$, the composition $f_{n-1}[s_{n-1}]\circ \dots\circ f_1[s_1]$ moves any point of $D$ by the maximum distance $d_j+2d_j+\dots+(n-1)d_j=\frac{n(n-1)}{2}d_j=(\omega-1) d_j$.
Define the rotated image $D'=f_{n-1}[s_{n-1}]\circ\cdots\circ f_1[s_1](D)$ based on the neighbors $s_1,\dots,s_{n-1}\in C_j$ of $q_1,\dots,q_{n-1}\in D$, respectively.
Since the subcloud $C_j$ is covered by the $d_j$-offset of $D$ and hence by the $\om d_j$-offset of $D'$, the approximation $d_{\vec H}(C_j, D')$ is non-strictly between the exact distance $d_j=d_{\vec H}(C_j, D)$ and its upper bound $\omega d_j$.
\medskip

The algorithm starts by finding $n-1$ `farthest' points $q_1,\dots,q_{n-1}\in D$, which are independent of $j$, in time $O(n|D|)$.
Here $q_1$ is a point of $D$ with a maximum distance $|q_1|$ to the origin, $q_2$ is a next `farthest' point of $D$ from $0$ and so on.
If some of these points have equal distances to $0$, they can be chosen in any order.
Though we do not know which points $s_1,\dots,s_{n-1}\in C_j$ become nearest neighbors of $q_1,\dots,q_{n-1}\in D$ after an optimal rotation, we will consider all (unit vectors of) points $s_1,\dots,s_{n-1}\in C_j$ to compute the approximate distance
$d'_j=\min\limits_{s_1,\dots,s_{n-1}\in C_j} d_{\vec H}(C_j,D')$ for $j=1,\dots,k$.
The minimization for all such points keeps both bounds: $d_j\leq d'_j\leq \omega d_j$.
\medskip

To minimize choices for $s_1,\dots,s_{n-1}$, we remove from the ordered list $p_1,\dots,p_k$ all points $p_i$ whose unit vectors $\vec p_i/|\vec p_i|$ appear earlier with smaller indices.
Then we consider each variable point $s_i$ from the remaining list $C'$ in increasing order of distances from the origin.
For any chosen points $s_1,\dots,s_{n-1}$, we compute the rotated image $D'=f_{n-1}[s_{n-1}]\circ\cdots\circ f_1[s_1](D)$ by computing matrix products in time $O(|D|)$, where we skip polynomial factors of the fixed dimension $n$.
\medskip

To compute the approximation $d'_j=d'_{\vec R}(C_j, D)=\min\limits_{s_1,\dots,s_{n-1}\in C_j} d_{\vec H}(C_j, D')$, we keep the current minimum of $d'_j$, which will be updated after getting the distance $d_{\vec H}(C_j, D')$ for every new choice of $s_1,\dots,s_{n-1}\in C_j$.
For each rotated image $D'=f_{n-1}[s_{n-1}]\circ\cdots\circ f_1[s_1](D)$, we run the internal loop for $j=1,\dots,k$.
For each point $p_j$ from the ordered full cloud $C$, we compute the distance $d(p_j,D')=\min\limits_{q\in D'}|p_j-q|$ to its nearest neighbor in $D'$ in time $O(|D|)$.
We use the previous iteration for $j-1$ to get $d_{\vec H}(C_j,D')=\max\{ d_{\vec H}(C_{j-1},D'), d(p_j,D')\}$, where we set $d_{\vec H}(C_0,D')=0$.
If all $s_1,\dots,s_{n-1}\in C_j$ and the current value of $d'_j$ 
is larger than $d_{\vec H}(C_j,D')$, we update $d'_j:=d_{\vec H}(C_j,D')$.
Since $C$ has $O(\ar{C}^{n-1})$ points $s_1,\dots,s_{n-1}$ with distinct normalized vectors, which determine $D'=f_{n-1}[s_{n-1}]\circ\cdots\circ f_1[s_1](D)$, the total time is $O(|C|\ar{C}^{n-1}|D|)$.
\end{proof}

Lemma~\ref{lem:re-ordering} implies that ordering lists of pairs of $\ep$-perturbations will keep $\ep$-closeness of corresponding ordered values.
This result will help prove Lemma~\ref{lem:lower_bound}.

\begin{lem}[re-ordering of $\ep$-close values]
\label{lem:re-ordering}
For any $\ep\geq 0$, let $C=\{c_{1},\dots,c_{k}\}$ and $D=\{d_{1},\dots,d_{k}\}$ satisfy $|c_{i}-d_{i}|\leq\ep$, $i=1,\dots,k$.
For any $i=1,\dots,k$, let $c_{(i)},d_{(i)}$ be the $i$-th largest values in $C,D$, respectively.
Then $|c_{(i)}-d_{(i)}|\leq\ep$.
\bs
\end{lem}
\begin{proof}
We consider any real $d_i$ an $\ep$-perturbation of the corresponding value $c_i$ for $i=1,\dots,k$.
Assume towards contradiction that $c_{(i)}<d_{(i)}-\ep$, so all $i$ smallest values of $C$ are less than $d_{(i)}-\ep$.
Then all $\ep$-perturbations of these $i$ values in $D$ are less than $d_{(i)}$, so $D$ has $i$ values that are strictly smaller than $d_{(i)}$.
This conclusion contradicts that $d_{(i)}$ is the $i$-th largest value in $D$.
The assumption $c_{(i)}>d_{(i)}+\ep$ similarly leads to contradiction.
Hence, after writing both $C,D$ in increasing order, their $i$-th largest values remain $\ep$-close. 
\end{proof}

\begin{proof}[Proof of Lemma~\ref{lem:lower_bound}]
Definition~\ref{dfn:tolerant_metric}c of $\ep=\BT([C(S,p;\al)],[C(Q,q;\al)])$ implies that, for a suitable isometry $f\in\Or(\R^n)$, the image $f(C(S,p;\al-\ep)-\vec p)$ is covered by the $\ep$-offset of $C(Q;q;\al)-\vec q$ shifted by $q$ to the origin. 
Since $\ep$ is smaller than a minimum half-distance between points of $S,Q$, the above covering establishes a bijection $g$ with all (at least k) neighbors of $p$ and $q$ in their $(\al-\ep)$-clusters.
\smallskip

The covering condition above means that the corresponding neighbors are at a maximum distance $\ep$ from each other.
The triangle inequality implies that the distances from corresponding neighbors to their centers $p,q$ differ by at most $\ep$.
The ordered distances from $p,q$ to their $k$ neighbors in the $(\al-\ep)$-clusters form the rows of $p,q$ in $\PDD(S;k),\PDD(Q;k)$. 
The bijection $g$ may not respect their order.
By Lemma~\ref{lem:re-ordering} the ordered distances with the same indices are $\ep$-close.
So the $L_{\infty}$ distance between the rows of $p,q$ is at most $\ep$. 
\end{proof}

The proof of Lemma~\ref{lem:upper_bounds}(c) referred to Lemma~\ref{lem:finite_symmetry}, which was briefly proved in the 2nd paragraph in \cite[p.~20]{delone1976local} without a formal statement.
We stated and prove this auxiliary result below to make all arguments complete.

\begin{lem}[finite symmetry group]
\label{lem:finite_symmetry}
Let a periodic point set $S\subset\R^n$ be $n$-dimensional, i.e. $S$ is not contained a lower-dimensional affine subspace of $\R^n$.
Then the symmetry group $\sym(S,p;2R(S))$ is finite for any point $p\in S$.
\end{lem}
\begin{proof}
Recall that the covering radius $R(S)$ is the largest radius $R$ of an open ball $B(q;R)$ within the complement $\R^n\setminus S$ for $q\in\R^n$.
Consider any such ball $\bar B(q;R(S))$ whose boundary sphere passes through the given point $p\in S$ and whose interior contains no points of $S$.
Then the closed ball $\bar B(q;R(S))$ should include at least one more point $p_1\in S\setminus\{p\}$ with $|p-p_1|\leq 2R(S)$.
Otherwise, the ball $\bar B(q;R(S))$ can be slightly expanded from $p\in S$ without including any points of $S\setminus\{p\}$, which contradicts the definition of the covering radius $R(S)$.
\smallskip

If $n\geq 2$, consider another open ball $B(q_1;R(S))\subset\R^n\setminus S$ that touches at $p$ the straight line $L(p,p_1)$ through $p,p_1$. 
Then the closed ball $\bar B(q_1;R(S))$ should include at least one more point $p_2\in S$ outside the line $L(p,p_1)$, so $|p-p_2|\leq 2R(S)$ and $p,p_1,p_2\in S$ span the 2-dimensional plane $L(p,p_1,p_2)$.
Otherwise the closed ball $\bar B(q_1;R(S))$ can be slightly expanded from $p\in S$ on the boundary $\partial B(q_1;R(S))$ without including any points of $S\setminus\{p,p_1\}$.
\smallskip

If $n\geq 3$, consider another open ball $B(q_2;R(S))\subset\R^n\setminus S$ that touches at $p$ the plane $L(p,p_1,p_2)$.
Then the closed ball $\bar B(q_2;R(S))$ should include at least one more point $p_3\in S$ outside the plane $L(p,p_1,p_2)$.
Then $p,p_1,p_2,p_3\in S$ span the 3-dimensional subspace $L(p,p_1,p_2,p_3)$ and so on until we find $n+1$ affinely independent points $p,p_1,\dots,p_n\in S$ such that each $p_i$ is at a maximum distance $2R(S)$ from $p$ for $i=1,\dots,n$.
Since the cluster $C(S,p;2R(S))$ contains $n+1$ affinely independent points, its symmetry group $\sym(S,p;2R(S))$ is finite.
\end{proof}

\begin{exa}[Earth Mover's Distance for lattices with bottleneck distance $d_B=+\infty$]
\label{exa:EMD_PDD}
The 1D lattices $S=\Z$ and $Q=(1+\de)\Z$ with the bottleneck distance $d_B(S,Q)=+\infty$ have PDD consisting of a single row (as for any lattice).
For instance, $\PDD(S;4)=(1,1,2,2)$ and $\PDD(Q;4)=(1+\de,1+\de,2+2\de,2+2\de)$. 
For the common stable radius $\al=2+2\de$, Example~\ref{exa:EMD_cluster} computed $\EMD(I(S;\al),I(Q;\al))=2\de$.
Theorem~\ref{thm:lower_bound} considers the maximum number $k$ of points in clusters of $S,Q$ with the radius $\al-2\de=2$, so $k=2$.
\medskip

Then $\EMD(\PDD(S;2),\PDD(Q;2))$ equals the $L_{\infty}$ distance $\de$ between the short rows $(1,1)$ and $(1+\de,1+\de)$.
The above computations illustrate the lower bound 
$\EMD(\PDD(S;2),\PDD(Q;2))=\de\leq\EMD(I(S;\al),I(Q;\al))=2\de.$ 
This inequality becomes equality for the larger stable radius $\al=2+4\de$, because the clusters of $S,Q$ with the radius $\al-2\de=2+2\de$ contain $k=4$ points.
The $L_{\infty}$ distance between $(1,1,2,2)$ and $(1+\de,1+\de,2+2\de,2+2\de)$ is $2\de$  for $\de<\frac{1}{8}$, so $\EMD(\PDD(S;4),\PDD(Q;4))=2\de=\EMD(I(S;2+4\de),I(Q;2+4\de))$.
\bs
\end{exa}

\begin{exa}[lower bound for a distance between square and hexagonal lattices]
\label{exa:lower_bound}
The square lattice $\La_4$ and hexagonal lattice $\La_6$ with minimum inter-point distance $1$ have a common stable radius $\al=2$ as shown in Fig.~\ref{fig:square_vs_hexagon}.
The maximum number of points in the stable 2-clusters is $k=12$.
The rows 
$\PDD(\La_4;12)=(1,1,1,1,\sqrt{2},\sqrt{2},\sqrt{2},\sqrt{2},2,2,2,2)$ and $\PDD(\La_6;12)=(1,1,1,1,1,1,\sqrt{3},\sqrt{3},\sqrt{3},\sqrt{3},\sqrt{3},\sqrt{3})$ have the $L_{\infty}$ distance $\max\{\sqrt{2}-1,2-\sqrt{3}\}=\sqrt{2}-1$, which coincides with $\EMD(I(\La_4;2),I(\La_6;2))$ in Example~\ref{exa:square_vs_hexagon}.
\bs
\end{exa}

The latest version includes Algorithms~{\ref{algorithm:isoset_EMD}} and~{\ref{algorithm:d_M}}, which implement the approximations from section~{\ref{sec:algorithms}} of new metrics defined in section~{\ref{sec:metrics}}.
Tables~A.1-A.2 with near-duplicates and run times on the CSD and GNoME are new.  
\medskip

\begin{algorithm}[H]
    \DontPrintSemicolon
    \caption{{Pseudocode for the boundary tolerant metric BT between $\al$-clusters and EMD between isosets in Theorem{~\ref{thm:approximate_BT}} and Corollary~{\ref{cor:approximate_EMD}}, where the max-min distance $d_{\vec M}$ is approximated by Algorithm~{\ref{algorithm:d_M}}.}}
    \label{algorithm:isoset_EMD}
    \medskip
    \SetKwInOut{Input}{Input}
    \SetKwInOut{Output}{Output}
    \Input{Isosets \texttt{is1} = $I(S;\al)$, \texttt{is2} = $I(Q;\al)$, with weights \texttt{w1}, \texttt{w2}}
    \Output{$\EMD(I(S;\al),I(Q;\al))$}

    \medskip

    \texttt{distance\_matrix = zeros(len(is1), len(is2))}\\

    \For{\texttt{\emph{i in range(len(is1))}}}{
        \For{\texttt{\emph{j in range(len(is2))}}}{
            \texttt{BT = max($d_{\vec M}$(C, D), $d_{\vec M}$(D, C))}\\
            \texttt{distance\_matrix[i, j] = BT}\\
        }
    }

    \texttt{emd = EMD(w1, w2, distance\_matrix)}\\

    \Return{} \texttt{emd}
\end{algorithm}

To confirm near-duplicates in the CSD and GNoME database by complete isosets,
we first filtered out pairs of crystals that have $L_\infty\geq 10^{-4}\angstrom$ on faster invariants ADA100.
The full table of the resulting 4385 pairs with all distances and running times is in the supplementary materials.
In most pairs, crystals belong to the same 6-letter code family because their structures are either polymorphs (different phases with the same composition) or slightly different versions determined under different temperatures or pressures.
However, 398 pairs consist of geometric near-duplicates that were assigned to (unexpectedly) different families.
In almost all these cases, the EMD metric on isosets was 0 after rounding to $10^{-10}$ (floating point error) in Angstroms. 
Table~A.1 shows all 25 crystals where the $\EMD$ metric was only slightly above 0.

\begin{algorithm}[H]
    \DontPrintSemicolon
    \caption{{Pseudocode for the directed max-min distance $d_{\vec M}$ in Definition~{\ref{dfn:directed_distances}}(b) by using an approximation of $d_{\vec R}$ in Lemma~{\ref{lem:rot-inv_distance}}.}}
    \label{algorithm:d_M}
    \medskip
    \SetKwInOut{Input}{Input}
    \SetKwInOut{Output}{Output}
    \Input{finite clouds C, D (ordered by distance to the origin $0$)}
    \Output{$d_{\vec M}(C, D)$}

    \medskip

     \texttt{alpha = max(norm(C[-1]), norm(D[-1]))}\\
    \texttt{q1, d\_R, max\_d, res = D[-1], [inf] * len(C), -inf, -inf}\\
    
    \For{\texttt{\emph{q in D}}}{
        \lIf{\texttt{\emph{d := perp\_dist(q, q1) > max\_d}}}{
            \texttt{max\_d = d; q2 = q}
        }
    }

    \For{\texttt{\emph{i1, p1 in enumerate(C)}}}{

        \texttt{R1b = rotation\_to\_align(q1, p1)} \\

        \texttt{q1\_, q2\_ = dot(R1b, q1), dot(R1b, q2)}

        \For{\texttt{\emph{i2, p2 in enumerate(C)}}}{

            \lIf{\texttt{\emph{i1 == i2}}}{\Continue}

            \texttt{src\_normal = cross(q1\_, q2\_)} \\ 

            \texttt{R2b = identity(3)}

            \If{\texttt{\emph{norm(src\_normal) != 0}}}{
                \texttt{tar\_normal = cross(p1, p2)}

                \If{\texttt{\emph{norm(tar\_normal) == 0}}}{
                
                    \If{\texttt{\emph{dot(p1, p2) < 0}}}{
                        \texttt{R2b = rotation\_to\_align(p1, p2)}
                    }
                }
                \Else{
                    \texttt{R2b = rotation\_to\_align(src\_normal, tar\_normal)}
                }
            }

            \texttt{d\_H = Hausdorff\_dist(C, dot(D, R2b $\times$ R1b))}
            \texttt{}

            \lIf{\texttt{\emph{v := min(d\_H, alpha - norm(p)) > res}}}{
                \texttt{res = v}
            }
        }
    }
    
    \Return{} \texttt{res}

\end{algorithm}

\begin{table}[H]
\label{tab:CSD}
\caption{The first pair consists of rigidly different mirror images from Fig.~{\ref{fig:mirror_images}}~(right).
All others are geometric near-duplicates from (surprisingly) different families in the CSD, confirmed by tiny values of the EMD metric on isosets.
The distance units are in \emph{attometers}: 1 am $=10^{-8}\angstrom=10^{-18}$ meter.
The run times in milliseconds (ms) depend on the cluster size (maximum number of atoms in $\al$-clusters) according to Theorem~\ref{thm:compute_isoset} and Corollary~\ref{cor:approximate_EMD}.}
\hspace*{-20mm}
\begin{tabular}{llrrrr}
CSD id1  & CSD id2  & EMD\_isosets, am & isosets time, ms & EMD\_isosets time, ms & cluster size \\
WODLOS   & XAWGAE   & 85856.22         & 129.619         & 1204.58               & 7            \\
TAFQIA   & VAVQIS   & 952.96           & 1690.949        & 321603.86             & 20           \\
FIJKIU   & IPEQUR   & 728.43           & 407.579         & 77455.47              & 16           \\
JIZMIR01 & JIZNAK   & 496.08           & 40.454          & 634.73                & 5            \\
HIYVUG01 & MASPIF   & 334.62           & 35.518          & 543.45                & 7            \\
KIVXEW10 & KIWCEC   & 125.03           & 22.456          & 32.32                 & 5            \\
XAYZOP   & ZEMDAZ   & 89.47            & 301.217         & 1697.26               & 4            \\
KIVXEW10 & KIWCEC28 & 83.07            & 21.701          & 32.49                 & 5            \\
AFIBOH   & NENCUF   & 31.67            & 126.582         & 1160.58               & 5            \\
KIVXEW07 & KIWCEC09 & 31.11            & 22.287          & 36.95                 & 5            \\
KIVXEW07 & KIWCEC11 & 31.11            & 22.434          & 36.75                 & 5            \\
KIVXEW11 & KIWCEC26 & 26.11            & 21.646          & 32.33                 & 5            \\
SERKIL   & SERKOR   & 23.78            & 2444.885        & 18485.57              & 6            \\
ADESAG   & REWPOB   & 5.81             & 54.675          & 5689.27               & 15           \\
GEQRAX   & IFOQOL   & 0.05             & 265.4           & 2090.42               & 6            \\
BUKYEN   & UYOCES   & 0.03             & 398.129         & 15739.11              & 11           \\
GOHYOT   & VIHCEY   & 0.01             & 100.031         & 940.16                & 5            \\
JUMCUP   & QAHBOT   & 0.01             & 179.367         & 4234.6                & 5            \\
CALMOV   & CALNAI   & 0.01             & 128.437         & 3913.32               & 4            \\
NABKOT   & ZIVSEF   & 0.01             & 75.401          & 796.14                & 5            \\
LIBGAE   & VESJUY   & 0.01             & 41.535          & 403.87                & 3            \\
AMEVEV   & OLERON   & 0                & 70.172          & 558.78                & 4            \\
SIHFIZ   & TEZBUV   & 0                & 207.984         & 1761.39               & 5            \\
XATCAA   & ZAQMEN   & 0                & 60.254          & 394.74                & 4            \\
PIDREA   & XIZNOL   & 0                & 94.135          & 243.46                & 5              
\end{tabular}
\end{table}

\begin{table}[H]
\label{tab:GNoME}
\caption{{After excluding 3248 exact numerical duplicates from} \cite[Table~1]{anosova2024importance}, {the next 25 pairs of closest near-duplicates in the GNoME database are confirmed by tiny values of EMD on isosets, see the first pair in Fig.~\ref{fig:exact_duplicates}.
The distance units are \emph{attometers}: 1 am $=10^{-8}\angstrom=10^{-18}$ meter.
The run times are in milliseconds (ms).
The cluster size is the maximum number of atoms in $\al$-clusters.
The full table of $2858\time 2$ pairs is in the supplementary materials.}}
\hspace*{-20mm}
\begin{tabular}{lllrrr}
GNoME id1  & GNoME id2  & EMD\_isosets, am & isosets time, ms & EMD\_isosets time, ms & cluster size \\
\hline
1547d30046 & ddc216e80c & 1                & 1.659            & 434.362               & 14           \\
b4065a4798 & e78d3559e6 & 1.7              & 3.034            & 13.271                & 6            \\
98ab164895 & df1252bc44 & 2                & 1.002            & 419.142               & 14           \\
0de9d25713 & b1733941a7 & 2.7              & 1.971            & 49.816                & 6            \\
0e79f7c053 & 6cf951ac6f & 3                & 1.035            & 429.487               & 14           \\
07ece241f0 & 45cacc8d45 & 3.2              & 0.618            & 14.374                & 6            \\
a58dc74a92 & c16bf63220 & 4.1              & 2.532            & 641.086               & 14           \\
5023e3a4b8 & 8f7ffb4d4a & 4.6              & 2.776            & 10.02                 & 6            \\
3198d1a3ea & 35f67abe6d & 5                & 1.031            & 403.398               & 14           \\
6826b81efb & 76ee112799 & 5                & 0.985            & 407.618               & 14           \\
6826b81efb & e9be17f0ee & 5                & 1.008            & 404.306               & 14           \\
2cff5f2fa0 & f470a5f6fa & 5.3              & 0.635            & 169.911               & 17           \\
2ce912f039 & 9de239ee0c & 5.5              & 0.632            & 3.456                 & 2            \\
c9f5a7a51b & fd9f40e0e1 & 6                & 1.14             & 195.261               & 10           \\
18078e002b & aca2a892a5 & 6                & 1.028            & 421.009               & 14           \\
18078e002b & b9722429b1 & 6                & 1.18             & 445.453               & 14           \\
18078e002b & b702e73db3 & 6                & 1.035            & 414.325               & 14           \\
34b4204eee & adee17535b & 6                & 1.017            & 396.855               & 14           \\
506b8b5646 & 60d266db80 & 6                & 1.174            & 413.254               & 14           \\
506b8b5646 & ec7b789cb3 & 6                & 1.174            & 403.014               & 14           \\
780741962f & a19688f106 & 6.5              & 1.804            & 870.629               & 15           \\
780741962f & c6af1fc763 & 6.5              & 2.731            & 921.496               & 15           \\
780741962f & c64c3e245c & 6.5              & 1.792            & 829.979               & 15           \\
b06353561c & b6d2341d32 & 6.6              & 2.387            & 259.359               & 12           \\
ebb33e044c & ebc9a4db61 & 6.8              & 1.232            & 450.351               & 14           \end{tabular}
\end{table}

\begin{figure}[h!]
\includegraphics[height=100mm]{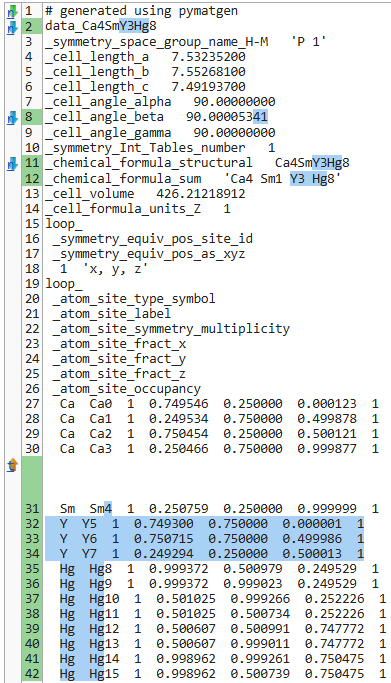}
\includegraphics[height=100mm]{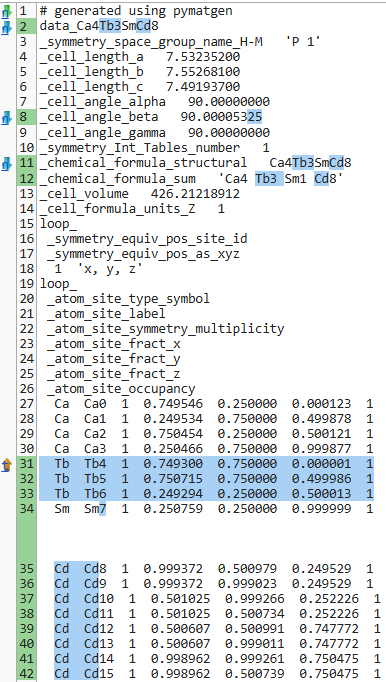}
\vspace*{-5mm}
\caption{The GNoME crystals 1547d30046 and ddc216e80c in the first row of Table~A.2 are compared as texts by https://text-compare.com.
All differences are highlighted in blue.
}
\label{fig:exact_duplicates}
\end{figure}
\vspace*{-2mm}

\begin{figure}[h!]
\includegraphics[height=30mm]{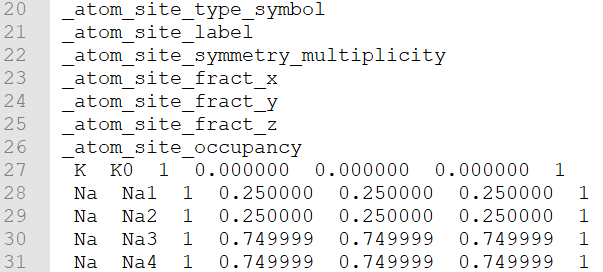}
\includegraphics[height=30mm]{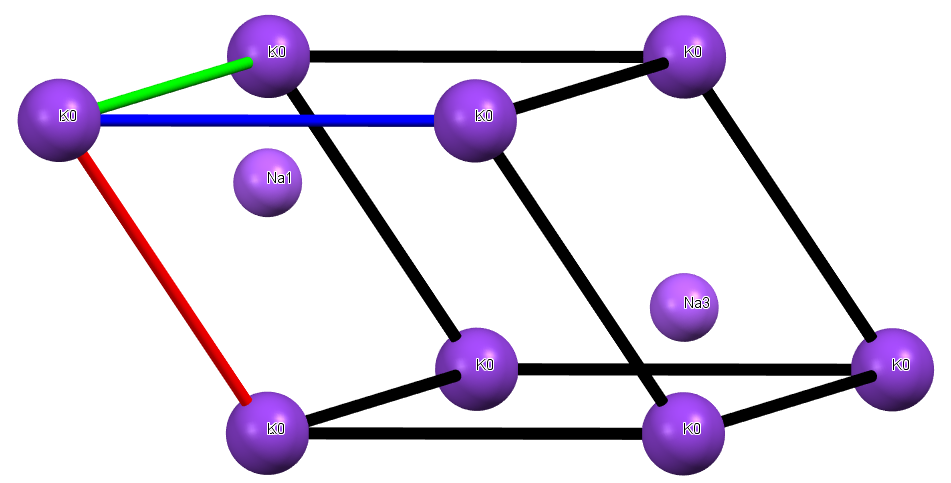}
\vspace*{-8mm}
\caption{Different entries cdc06a1a2a and 0e2d8f26d6 in the GNoME database are not only identical symbol by symbol but also contain two pairs of atoms (Na1=Na2 and Na3=Na4) at the same positions.
\textbf{Left}: a screenshot from the CIF.
\textbf{Right}: Mercury visualization can show only one atom in each pair of coinciding atoms, e.g. only Na1 and not Na2 from the CIF.
}
\label{fig:repeated_atoms}
\end{figure}

Fig.~\ref{fig:repeated_atoms} shows the most striking pair of exact duplicates in the GNoME is cdc06a1a2a and 0e2d8f26d6, whose CIFs are identical symbol by symbol in addition to two pairs of atoms at the same positions (Na1=Na2 and Na3=Na4).

\end{document}